\documentclass[12pt,a4paper]{article}
\usepackage{epsfig,latexsym,amsfonts,amssymb,amsmath,amscd,epic,graphics,theorem,epstopdf}
\usepackage{tikz-cd}
\theoremstyle{plain}
\newtheorem{lemma}{Lemma}[section]

\newtheorem{remark}[lemma]{Remark}
\newtheorem{example}[lemma]{Example}
\newtheorem{theorem}[lemma]{Theorem}
\newtheorem{definition}[lemma]{Definition}

{\theorembodyfont{\rmfamily}  \font\ncsc=cmcsc10
 \font\ntt=cmtt12

\begin{document}
\newcommand{\val}{\operatorname{val}}\newcommand{\Val}{\operatorname{Val}}
\newcommand{\pperp}{\hbox{$\perp\hskip-6pt\perp$}}
\newcommand{\ssim}{\hbox{$\hskip-2pt\sim$}}\newcommand{\ini}{\operatorname{ini}}
\newcommand{\aleq}{{\ \stackrel{3}{\le}\ }}\newcommand{\Vers}{\mathrm{Vers}}
\newcommand{\ageq}{{\ \stackrel{3}{\ge}\ }}
\newcommand{\aeq}{{\ \stackrel{3}{=}\ }}
\newcommand{\bleq}{{\ \stackrel{n}{\le}\ }}
\newcommand{\bgeq}{{\ \stackrel{n}{\ge}\ }}
\newcommand{\beq}{{\ \stackrel{n}{=}\ }}
\newcommand{\cleq}{{\ \stackrel{2}{\le}\ }}
\newcommand{\cgeq}{{\ \stackrel{2}{\ge}\ }}
\newcommand{\ceq}{{\ \stackrel{2}{=}\ }}
\newcommand{\fm}{\mathfrak{m}}
\newcommand{\N}{{\mathbb N}}\newcommand{\T}{{\mathbb T}}
\newcommand{\A}{{\mathbb A}}
\newcommand{\K}{{\mathbb K}}
\newcommand{\Z}{{\mathbb Z}}\newcommand{\F}{{\mathbf F}}
\newcommand{\R}{{\mathbb R}}
\newcommand{\C}{{\mathbb C}}
\newcommand{\Q}{{\mathbb Q}}
\newcommand{\PP}{{\mathbb P}}
\newcommand{\cA}{{\mathcal A}}
\newcommand{\cB}{{\mathcal B}}
\newcommand{\cC}{{\mathcal C}}
\newcommand{\cD}{{\mathcal D}}
\newcommand{\cF}{{\mathcal F}}
\newcommand{\cI}{{\mathcal I}}
\newcommand{\cL}{{\mathcal L}}
\newcommand{\cM}{{\mathcal M}}
\newcommand{\cR}{{\mathcal R}}
\newcommand{\cO}{{\mathcal O}}
\newcommand{\cP}{{\mathcal P}}
\newcommand{\cS}{{\Sigma}}
\newcommand{\cT}{{\mathfrak T}}
\newcommand{\mcT}{{\mathcal T}}
\newcommand{\cU}{{\mathcal U}}
\newcommand{\cZ}{{\mathcal Z}}
\newcommand{\cOK}{\mathcal{OK}}
\newcommand{\cEA}{{\mathcal EA}}
\newcommand{\cNDF}{\mathcal{NDF}}
\newcommand{\mkfO}{{\mathfrak O}}
\newcommand{\go}{{\mathfrak o}}
\newcommand{\gc}{{\mathfrak c}}
\newcommand{\bbz}{{x_0}}
\newcommand{\Al}{\operatorname{Ad}}\newcommand{\Men}{\operatorname{Men}}
\newcommand{\EAl}{\operatorname{EAd}}\newcommand{\Ann}{\operatorname{Ann}}
\newcommand{\red}{{\operatorname{red}}}
\newcommand{\Pic}{{\operatorname{Pic}}}\newcommand{\Sym}{{\operatorname{Sym}}}
\newcommand{\QI}{{\operatorname{QI}}}\newcommand{\Div}{{\operatorname{Div}}}
\newcommand{\oDel}{{\widetilde\Del}}
\newcommand{\real}{{\operatorname{Re}}}\newcommand{\Aut}{{\operatorname{Aut}}}
\newcommand{\conv}{{\operatorname{conv}}}\newcommand{\Ima}{{\operatorname{Im}}}
\newcommand{\Span}{{\operatorname{Span}}}\newcommand{\Trop}{{\operatorname{Trop}}}
\newcommand{\Ker}{{\operatorname{Ker}}}\newcommand{\qind}{{\operatorname{qind}}}
\newcommand{\Cycle}{{\operatorname{Cycle}}}\newcommand{\OG}{{\operatorname{OG}}}
\newcommand{\Fix}{{\operatorname{Fix}}}\newcommand{\ina}{{\operatorname{in}}}
\newcommand{\sign}{{\operatorname{sign}}}
\newcommand{\even}{{\operatorname{even}}}
\newcommand{\odd}{{\operatorname{odd}}}
\newcommand{\com}{{\operatorname{com}}}
\newcommand{\ncom}{{\operatorname{ncom}}}
\newcommand{\nmob}{{\operatorname{nmob}}}
\newcommand{\bound}{{\operatorname{bound}}}
\newcommand{\nbound}{{\operatorname{ends}}}
\newcommand{\Inn}{{\operatorname{In}}}
\newcommand{\Ex}{{\operatorname{Ex}}}
\newcommand{\opp}{{\operatorname{opp}}}\newcommand{\Par}{{\operatorname{Par}}}
\newcommand{\Cheb}{{\operatorname{Cheb}}}\newcommand{\arccosh}{{\operatorname{arccosh}}}
\newcommand{\Card}{{\operatorname{Card}}}
\newcommand{\alg}{{\operatorname{alg}}}
\newcommand{\ord}{{\operatorname{ord}}}
\newcommand{\mt}{{\operatorname{mult}}}
\newcommand{\cheb}{{\operatorname{Cheb}}}
\newcommand{\oi}{{\overline i}}\newcommand{\oGamma}{{\overline\Gamma}}
\newcommand{\oj}{{\overline j}}\newcommand{\oh}{{\overline h}}
\newcommand{\ob}{{\overline b}}
\newcommand{\os}{{\overline s}}
\newcommand{\oa}{{\overline a}}
\newcommand{\oy}{{\overline y}}
\newcommand{\ow}{{\overline w}}
\newcommand{\ot}{{\overline t}}
\newcommand{\oz}{{\overline z}}
\newcommand{\eps}{{\varepsilon}}
\newcommand{\proofend}{\hfill$\Box$\bigskip}
\newcommand{\Int}{{\operatorname{Int}}}
\newcommand{\pr}{{\operatorname{pr}}}
\newcommand{\Hom}{{\operatorname{Hom}}}
\newcommand{\Ev}{{\operatorname{Ev}}}
\newcommand{\im}{{\operatorname{Im}}}\newcommand{\br}{{\operatorname{br}}}
\newcommand{\sk}{{\operatorname{sk}}}\newcommand{\DP}{{\operatorname{DP}}}
\newcommand{\const}{{\operatorname{const}}}
\newcommand{\Sing}{{\operatorname{Sing}}\hskip0.06cm}
\newcommand{\conj}{{\operatorname{Conj}}}
\newcommand{\Cl}{{\operatorname{Cl}}}
\newcommand{\Crit}{{\operatorname{Crit}}}
\newcommand{\Ch}{{\operatorname{Ch}}}
\newcommand{\discr}{{\operatorname{discr}}}
\newcommand{\Tor}{{\operatorname{Tor}}}
\newcommand{\Conj}{{\operatorname{Conj}}}
\newcommand{\Log}{{\operatorname{Log}}}
\newcommand{\vol}{{\operatorname{vol}}}
\newcommand{\defect}{{\operatorname{def}}}
\newcommand{\codim}{{\operatorname{codim}}}
\newcommand{\tmu}{{\C\mu}}
\newcommand{\wt}{{\operatorname{wt}}}
\newcommand{\ov}{{\overline v}}
\newcommand{\ox}{{\overline{x}}}
\newcommand{\bw}{{\boldsymbol w}}
\newcommand{\hbw}{{\widehat\bw}}
\newcommand{\mfw}{{\mathfrak{w}}}
\newcommand{\bv}{{\boldsymbol v}}
\newcommand{\bn}{{\Phi}}
\newcommand{\bx}{{\boldsymbol x}}
\newcommand{\bd}{{\boldsymbol d}}
\newcommand{\bz}{{\boldsymbol z}}
\newcommand{\bL}{{\boldsymbol L}}
\newcommand{\bP}{{\boldsymbol P}}
\newcommand{\bp}{{\boldsymbol p}}
\newcommand{\bq}{{\boldsymbol p_{tr}}}
\newcommand{\be}{{\boldsymbol e}}
\newcommand{\bc}{{\boldsymbol c}}
\newcommand{\ba}{{\boldsymbol a}}
\newcommand{\bb}{{\boldsymbol b}}
\newcommand{\tet}{{\theta}}
\newcommand{\Del}{{\Delta}}
\newcommand{\bet}{{\beta}}
\newcommand{\kap}{{\kappa}}
\newcommand{\del}{{\delta}}
\newcommand{\sig}{{\sigma}}
\newcommand{\alp}{{\alpha}}
\newcommand{\Sig}{{\Sigma}}
\newcommand{\Gam}{{\Gamma}}
\newcommand{\gam}{{\gamma}}\newcommand{\idim}{{\operatorname{idim}}}
\newcommand{\Lam}{{\Lambda}}
\newcommand{\lam}{{\lambda}}
\newcommand{\SC}{{SC}}
\newcommand{\MC}{{MC}}
\newcommand{\nek}{{,...,}}
\newcommand{\cim}{{c_{\mbox{\rm im}}}}
\newcommand{\clM}{\tilde{M}}
\newcommand{\clV}{\bar{V}}
\newcommand{\rtm}{{\mu}}

\title{Real enumerative invariants relative to the
toric boundary and their
refinement}
\author{Ilia Itenberg
\and Eugenii Shustin}
\date{}
\maketitle
\begin{abstract}
We introduce 
new invariants of a class of toric surfaces (including the projective plane)
that arise from appropriate enumeration of real curves of genus
one and two.
These invariants admit a refinement
similar to the one introduced by Grigory Mikhalkin in the rational case.
\end{abstract}

\medskip

{\bf MSC-2010 classification:} Primary 14N10, Secondary 14J26, 14P05 

\medskip

\section{Introduction}\label{rel-intro}

Refined enumerative geometry, initiated in \cite{BG,GSh}, became one of the central topics in enumerative geometry with important links to closed and open Gromov-Witten invariants and to Donaldson-Thomas invariants.
In a big part of
known examples, refined invariants appear as one-parameter deformations
of complex enumerative invariants (see, for example, \cite{BG,FS,GSc,GSh}).
In his groundbreaking paper \cite{Mir}, G. Mikhalkin
has proposed a refined invariant
provided by enumeration of real rational curves
and has related this invariant to the refined tropical invariants of F. Block and L. G\"ottsche \cite{BG}.
Namely, he
introduced an integer-valued {\it quantum index}
for real algebraic curves in toric surfaces.
The quantum index is defined in \cite{Mir} for real curves satisfying the following assumptions:
the curve
has to intersect toric divisors only at real points and to be irreducible and {\it separating}; the latter condition means
that, in the complex point set of the normalization of the curve, the complement of the real part
is disconnected, {\it i.e.}, formed by two halves exchanged by the complex conjugation
(in fact, the quantum index is associated to a half of a separating real curve,
while the other half has the opposite quantum index; for detailed definitions, see Section \ref{sec2.3}).
Mikhalkin \cite{Mir} showed that, for an appropriate
kind of constraints,
a Welschinger-type enumeration of
real rational curves ({\it cf}. \cite{W1})
in a given divisor class and with a given quantum index
produces an invariant
which
can be directly related to the numerator
of a Block-G\"ottsche refined tropical invariant (represented as a fraction with the standard denominator).

The main goal of the present
paper is to introduce refined invariants of Mikhalkin's type
in the case of curves of genus $1$ and $2$.
We follow the ideas of \cite{Sh} and choose constraints so that every counted real curve of genus $g=1$ or $2$ appears
to be a {\it maximal} one
({\it i.e.}, it has $g+1$ global real branches), and hence is
separating.
More precisely, given a toric surface with the
tautological
real structure and a very ample divisor class, we fix maximally many
real points in a generic Menelaus position (see Section \ref{sec-td-m})
on the toric boundary of the positive quadrant, where
genus $g$
curves from the given linear system must
be tangent to toric divisors with prescribed even intersection multiplicities,
and we fix $g$ more generic real points inside different non-positive quadrants as
extra
constraints (through which the curves under consideration should pass).
There are finitely many real curves of genus $g$
matching the constraints,
and all these curves are separating.
Their halves have
quantum index, and
we equip each curve
with a certain Welschinger-type sign.
This gives rise to
a signed enumeration of halves of real curves of genus $g = 1$ and $2$
that match given constraints, belong to a given linear system,
and have a prescribed quantum index $\kappa$
(we also treat the case $g = 0$, slightly generalizing Mikhalkin's enumeration
to the situation with arbitrary even intersection multiplicities).

We prove that,
for some toric surfaces (including the projective plane;
in the case of genus $2$, in order to avoid cumbersome technicalities,
this is the only toric surface considered in the paper),
the result $W_g^\kappa(\Delta, \widehat\bw)$
of such enumeration,
where the sequence $\Delta$ encodes the orders of tangencies
and the sequence $\widehat\bw$ encodes the position of the constraints,
does not depend on the choice of a (generic)
sequence $\widehat\bw$;
precise statements can be found in
Section \ref{sec2.3} (rational case, Theorem \ref{at1}),
Section \ref{asec3} (elliptic case, Theorem \ref{at2}),
and Section \ref{sec-g2} (genus two case, Theorem \ref{g2t1}).
The obtained invariants are said to be {\it refined}.
In particular, we get new real enumerative invariants (without prescribing values for quantum index)
in genus one
and two.
Tropical counterparts of the above refined invariants, their computation
and comparison with the Block-G\"ottsche refined tropical invariants
are presented in the forthcoming papers \cite{IS, IS2}.

The
proofs of invariance are decomposed into four steps
(these steps are common for all genera $g=0$, $1$, and $2$):
\begin{itemize}
\item[(1)] identification of {\it walls}, {\it i.e.},
loci of codimension one in the space of constraints,
corresponding
to degenerations in the set of counted curves,
\item[(2)] analysis of (local) deformations of germs of degenerate curves at their singular points,
\item[(3)] establishing
{\it transversality conditions}
that ensure
the existence of one-dimensional families
in the space of constraints transversally crossing
the walls and realizing prescribed local deformations
of degenerate counted curves,
\item[(4)] verification of the persistence of the count in all wall-crossing events.
\end{itemize}

In
addition to the signed enumeration mentioned above,
we suggest another choice for a Welschinger-type sign
in the enumerative problems considered.
This choice also yields
an invariant enumeration.
The proofs of the invariance for these two signed enumerations
are literally the same,
so we treat them together.
It
turns out that the results
$W_g^\kappa(\Delta, \widehat\bw)$
and ${\widetilde W}_g^\kappa(\Delta, \widehat\bw)$
of two signed enumerations
(see Section \ref{sec-inv} for the notations)
coincide in the case $g = 0$ and in the case $g = 1$
(and seem to coincide in the case $g = 2$),
up to a global sign depending only on $\Delta$ and $\kappa$.
This statement is proved in the forthcoming paper \cite{IS}
using the tropical calculation of these invariants.
In fact, for the constraints that are sufficiently close to the tropical limit, the two introduced signs coincide
(up to the above mentioned global sign)
for individual curves under enumeration.
Probably, it is also the case in general, without the additional assumption on the constraints.

\vskip5pt

The
paper is organized as follows.
Section \ref{asec1} contains a preparatory material.
In Section \ref{sec-inv}, we introduce refined invariants arising in enumeration
of real curves of genus
zero, one and two
(the genus zero version slightly generalizes Mikhalkin's refined rational invariants \cite{Mir}).
Sections \ref{asec2},
\ref{asec5}, and \ref{asec6}
are devoted to the proofs of the invariance statements
in genus zero, one, and two, respectively.

\vskip5pt

{\bf Acknowledgements}.
We started this work during our stay at the Mittag-Leffler Institute,
Stockholm, in April 2018 and during the visit of the second author to the \'Ecole Normale Sup\'erieure, Paris, in June 2019, and we completed the work during our research stays at the Mathematisches Forschungsinstitut Oberwolfach in July 2021 and
March 2022 (in the framework of Research-in-Pairs program) and in February - March 2023
(in the framework of Oberwolfach Research Fellows program).
We are very grateful to these institutions for the support and excellent working conditions.
The first author was supported in part by the ANR grants ANR-18-CE40-0009 ENUMGEOM
and ANR-22-CE40-0014 SINTROP.
The second author has been partly supported by the Israeli Science
Foundation grant no. 501/18, and by the Bauer-Neuman Chair in Real and Complex Geometry.

\section{Preparation}\label{asec1}
\subsection{Convex lattice polygons and
real toric surfaces}\label{sec2.1.1}
Consider the lattice $\Z^2$
and its ambient plane $\R^2
=\Z^2\otimes\R$. Let $P\subset\R^2$
be a non-degenerate convex
lattice polygon.
For a lattice segment $\sigma\subset\R^2$
(respectively, a vector $\ba\in\Z^2\setminus\{(0, 0)\}$),
denote by $\|\sigma\|_\Z$ (respectively, $\|\ba\|_\Z$) its {\it lattice length}, {\it i.e.},
the ratio of the Euclidean length and the minimal length of a non-zero parallel integral vector.

Denote by $\Tor(P)$
the complex toric surface associated with $P$.
Let $\Tor(P)^\times \simeq (\C^*)^2$ (respectively, $\Tor(\partial P)$) be the dense orbit
(respectively, the union of all toric divisors)
of $\Tor(P)$.
The toric surface $\Tor(P)$ has the tautological real structure,
and the real part $\Tor_\R(P)$ of $\Tor(P)$
contains the positive
quadrant $\Tor_\R(P)^+\simeq(\R_{>0})^2$.
The closure
$\Tor_\R^+(P)$ of $\Tor_\R(P)^+$
(with respect to the Euclidean topology)
is diffeomorphically taken onto $P$ by the moment map.
We pull back the standard orientation and the metric
of $P\subset\R^2
$ to $\Tor_\R^+(P)$ and
induce an orientation and a metric
on the boundary $\partial\Tor_\R^+(P)$ of $\Tor_\R^+(P)$;
in particular, we get a cyclic order on
the sides of $P$; we call this order {\it positive}.
Denote by ${\mathcal L}_P$
the tautological line bundle over $\Tor(P)$; the global sections of ${\mathcal L}_P$ are spanned by the monomials
$z^\omega$, $z=(z_1,z_2)$, $\omega\in P\cap\Z^2$.

For each edge $\sigma\subset\partial P$, we consider the toric curve $\Tor(\sigma)$,
its dense orbit $\Tor(\sigma)^\times$,
the real part $\Tor_\R(\sigma)$ of $\Tor(\sigma)$,
and the positive half-axis $\Tor_\R(\sigma)^+ \subset \Tor(\sigma)^\times \cap \Tor_\R(\sigma)$.
The closure
$\Tor_\R^+(\sigma)$ of $\Tor_\R(\sigma)^+$ coincides with
$\Tor(\sigma)\cap\partial\Tor^+_\R(P)$.

\subsection{Toric degree and Menelaus condition}\label{sec-td-m}
A multi-set $\Delta \subset \Z^2 \setminus \{(0, 0)\}$ is said to be
\begin{itemize}
\item {\it balanced} if the sum of vectors in $\Delta$ is equal to $0$,
\item {\it non-degenerate} if the vectors of $\Delta$ span $\R^2$,
\item {\it even} if $\Delta \subset (2\Z)^2$.
\end{itemize}
A balanced non-degenerate multi-set $\Delta \subset \Z^2 \setminus \{(0, 0)\}$ is called a {\em {\rm (}toric{\rm )} degree}.

Let $\Delta \subset \Z^2 \setminus \{(0, 0)\}$ be a non-degenerate balanced
multi-set.
For each vector $\ba\in\Delta$, denote by
$\check\ba$ the vector obtained from $\ba$ by the counter-clockwise
rotation by $\pi/2$.
The vectors $\check\ba$ can be
attached to each other so that the next vector starts at the end of the preceding one in order to form a simple
broken line bounding a convex lattice polygon
whose boundary gets the counter-clockwise orientation.
The latter polygon is
denoted by $P_\Delta$; it is determined by $\Delta$ up to translation.
For a convex lattice polygon $\delta$,
we set $\cA(\delta)$ to be the Euclidean area of $\delta$
and $\cI(\delta)$ to be the number of integer points in the interior of $\delta$.
In the case of $P_\Delta$, we shortly write $\cA(\Delta)$ and $\cI(\Delta)$.
The elements of
$\Delta$ are denoted by
$$\ba^\sigma_i,\quad i=1, \ldots, n^\sigma,$$
where $\sigma$ ranges over the set $P^1_\Delta$ of
sides of the polygon $P_\Delta$, and $n^\sigma$
is the number of vectors in $\Delta$ that are outer normals to $\sigma$.

From now on, assume that $\Delta$ is even;
for each $\sigma \in P^1_\Delta$ and each $i = 1, \ldots, n^\sigma$, put
$2k^\sigma_i = \|\ba^\sigma_i\|_\Z$. 

\begin{definition}\label{ad1}
For each $\sigma \in P^1_\Delta$, consider a sequence $\bz^\sigma$ of $n^\sigma$ points in $\Tor(\sigma)^\times$,
and denote by $\bz$ the
double index sequence
{\rm (}the upper index being $\sigma \in P^1_\Delta$ and the lower index being $i = 1, \ldots, n^\sigma${\rm )}
formed by the sequences $\bz^\sigma$.
We say that the
sequence $\bz$
satisfies the Menelaus $\Delta$-condition {\rm (}cf. \cite[Section 5.1]{Mir}\/{\rm )}
if there exists a curve $C\in|{\mathcal L}_{P_\Delta}|$
such that $C$ does not contain toric divisors as
components and,
for each $\sigma\in P^1_\Delta$,
the scheme-theoretic intersection of $C$ and $\Tor(\sigma)$
coincides with $\sum_{i=1}^{n^\sigma}2k^\sigma_iz^\sigma_i$.
\end{definition}

\begin{lemma}\label{al1}
The
sequences $\bz$ satisfying the Menelaus $\Delta$-condition form
an algebraic
hypersurface
$\Men(\Delta)\subset \prod_{\sigma\in P^1_\Delta}\Tor(\sigma)^{n^\sigma}$.
\end{lemma}

{\bf Proof.}
We present an explicit equation
of $\Men(\Delta)$.
Consider a linear functional $\lambda:\R^2\to\R$ which is injective on $\Z^2$.
The maximal and the minimal points of $\lambda$ on $P_\Delta$ divide the boundary $\partial P_\Delta$
of $P_\Delta$
into two broken lines $B', B''$, and the maximal and minimal points of $\lambda$ on each side $\sigma$ define its orientation. An automorphism of $\Z^2$,
which takes an oriented side $\sigma$ of $P$ onto the naturally oriented segment
$[0,\|\sigma\|_\Z]$ of the first axis of $\R^2$,
defines an isomorphism of $\Tor(\sigma)^\times$ with
$\C^\times$.
Denote by $\xi^\sigma=(\xi_i^\sigma)_{i=1}^{n^\sigma}$ the sequence of images in
$\C^\times$
of the points of $\bz^\sigma$. Then, $\Men(\Delta)$
is given by the equation
 $$\prod_{\sigma\subset B'}\prod_{i=1}^{n^\sigma}(\xi_i^\sigma)^{2k_i^\sigma}=
\prod_{\sigma\subset B''}\prod_{i=1}^{n^\sigma}(\xi_i^\sigma)^{2k_i^\sigma}\ .$$
\proofend

The hypersurface $\Men(\Delta)$ is called {\it Menelaus hypersurface}.
Notice that
it is reducible. More precisely,
$\Men(\Delta)$ splits into $2k_0$ components, where
$$k_0=\gcd\{k_i^\sigma\ :\ i=1,...,n^\sigma,\ \sigma\subset\partial P_\Delta\}.$$
Denote by $\Men(\Delta)_{red}$ the component given by the equation
\begin{equation}
\prod_{\sigma\subset B'}\prod_{i=1}^{n^\sigma}(\xi_i^\sigma)^{k_i^\sigma/k_0}=
\prod_{\sigma\subset B''}\prod_{i=1}^{n^\sigma}(\xi_i^\sigma)^{k_i^\sigma/k_0}\ .
\label{ae1}
\end{equation}

Note that the introduced above isomorphism $\Tor(\sigma)^\times \simeq \C^\times$ takes $\Tor_\R(\sigma)^+$ onto
$\R_{>0}\subset\C^\times$.
For 
$\rho>0$, denote by $\Men_\R^\rho(\Delta)$ the part of $\Men(\Delta)_{red}$
specified by the condition that, for each $\sigma\in P^1_\Delta$,
all points of the subsequence $\bw^\sigma$
belong to $\Tor_\R^+(\sigma)$ and lie
at the distance $\ge\rho$ (with respect to the metric fixed in Section \ref{sec2.1.1})
from the endpoints of $\Tor_\R^+(\sigma)$. Here, we always assume that, in each interval $\Tor_\R^+(\sigma)$,
the complement to $\rho$-neighborhoods of the endpoints
is non-empty.

\begin{lemma}\label{al2}

The {\rm (}metric{\rm )} closure
$\overline \Men_\R^\rho(\Delta)\subset\prod_{\sigma\subset\partial P}(\Tor_\R^+(\sigma))^{n^\sigma}$
is diffeomorphic to a convex polytope.
\end{lemma}

{\bf Proof.}
Since $\Men_\R^\rho(\Delta)$ is given by equation (\ref{ae1}) with positive variables $\xi_i^\sigma$ satisfying restrictions of the form
$0<\const_1
\le\xi_i^\sigma\le\const_2
<\infty$,
the coordinate-wise logarithm takes $\overline \Men_\R^\rho(\Delta)$
onto a hyperplane section of a convex polytope.
\proofend

\subsection{Deformations of plane curve germs}\label{sec-def}
Here we introduce notations and recall some known facts on deformations of plane curve germs.

Let $(C,z)\subset(\C^2,z)$ be the germ of a reduced 
curve $C$ with an isolated singular point $z$, and let $L\supset\{z\}$ be a smooth curve intersecting $(C,z)$ only at $z$.
The space $\Vers(C,z)$, defined as
the germ at zero of ${\mathcal O}_{C,z}/{\mathfrak m}_{C,z}^s$, $s\gg0$, where ${\mathfrak m}_{C,z}$ is the maximal ideal in the local ring ${\mathcal O}_{C,z}$,
serves as the base of a versal deformation of the germ $(C,z)$.
Denote by $Q_1,...,Q_r$ the irreducible components
of $(C, z)$, {\it i.e.}, the local branches of $C$ centered at $z$,
and consider the following loci in $\Vers(C,z)$.
\begin{enumerate}
\item[(a)] The equisingular locus $ES(C,z)\subset\Vers(C,z)$ parameterizes curve germs
with a singular point topologically equivalent to $(C,z)$.

It is smooth (see \cite[Item (ii) on page 435]{DH}),
and by \cite[Theorem 1 and formula (19)]{GuS}, its tangent space $T(ES(C,z))$
is contained in
$$\left\{\varphi\in \Vers(C,z)\ \Big|\ \begin{array}{c}{\ord\varphi\big|_{Q_i}\ge2\delta(Q_i)
+\sum_{j\ne i}(Q_i\cdot Q_j)}\\ {\qquad\qquad+\mt Q_i-1,\quad
i=1,...,r}\end{array}\right\},$$
where $(Q_i\cdot Q_j)$ stands for
the intersection multiplicity of the local branches $Q_i$ and $Q_j$, the $\delta$-invariant $\delta(Q_i)$
denotes the maximal number of nodes that may appear in a deformation of the local branch $Q_i$,
the multiplicity of $Q_i$ at $z$ is denoted by $\mt Q_i$,
and at last, $\ord\varphi\big|_{Q_i}$
is the vanishing order of $\varphi(x(t),y(t))$, where $(x(t),y(t))$ is a parameterization of $Q_i$.
\item[(b)] The equisingular loci
$ES_L(C,z)\subset ES(C,z)$ and $ES^{fix}_L(C,z)\subset ES_L(C,z)$
additionally constrained by the requirement to intersect $L$ at a point close to $z$
with multiplicity $(C\cdot L)_z$, or to intersect $L$ at $z$ with multiplicity $(C\cdot L)_z$,
respectively.
Both loci are smooth,
and by \cite[Theorem 2 and formulas (21), (22)]{GuS},
their tangent spaces satisfy
$$T(ES_L(C,z))\subset\left\{\varphi\in \Vers(C,z)\ \Big|\ \begin{array}{c}{\ord\varphi\big|_{Q_i}\ge2\delta(Q_i)+\sum_{j\ne i}(Q_i\cdot Q_j)}\\ {\qquad\qquad+(Q_i\cdot L)-1,\quad i=1,...,r}\end{array}\right\},$$
$$T(ES^{fix}_L(C,z))\subset\left\{\varphi\in \Vers(C,z)\ \Big|\ \begin{array}{c}{\ord\varphi\big|_{Q_i}\ge2\delta(Q_i)+\sum_{j\ne i}(Q_i\cdot Q_j)}\\ {+\mt Q_i+(Q_i\cdot L)-1,\quad i=1, \ldots, r}\end{array}\right\}.$$
\item[(c)] The equigeneric locus $EG(C,z)\subset\Vers(C,z)$
parameterizing equigeneric (or $\delta=\const$)
deformations.
It is smooth if and only if
all local branches $Q_1$, \ldots,
$Q_r$
are smooth.
In any case, its tangent cone is the linear space
$$
J^{cond}(C,z)/{\mathfrak m}_{C,z}^s=\left\{\varphi\in\Vers(C,z)\ \Big|\
\begin{array}{c}{\ord\varphi\big|_{Q_i}\ge2\delta(Q_i)+\sum_{j\ne i}(Q_i\cdot Q_j)}\\
{i=1, \ldots, r}
\end{array}\right\}
$$
of codimension $\codim EG(C,z)=\delta(C,z)$,
where
$$J^{cond}(C,z)=\Ann(\bn_*{\mathcal O}_{\widehat C}/{\mathcal O}_C)\subset{\mathcal O}_{C,z}$$
is the conductor ideal (see details in \cite[Item (iii) on page 435 and Theorem 4.15]{DH}),
and $\delta(C,z)=\dim{\mathcal O}_{C,z}/J^{cond}(C,z)$.
\item[(d)] The equigeneric loci $EG_L(C,z)\subset EG(C,z)$ and
$EG_L^{fix}(C,z)\subset EG_L(C,z)$ additionally constrained by
the requirement to intersect $L$ at a point close to $z$ with multiplicity $(C\cdot L)_z$, or to intersect $L$ at $z$
with multiplicity $(C\cdot L)_z$, respectively.
If $(C,z)$ is irreducible, {\it i.e.}, $r = 1$,
then these loci are smooth \cite[Lemma 3(1)]{Sh2}, and their tangent spaces satisfy
        $$T(EG_L(C,z))=\left\{\varphi\in\Vers(C,z)\
        \Big| \ \ord\varphi\big|_{(C,z)}\ge2\delta(C,z)+(C\cdot L)_z-1\right\},$$
        $$T(EG^{fix}_L(C,z))=\left\{\varphi\in\Vers(C,z)\
        \Big| \ \ord\varphi\big|_{(C,z)}\ge2\delta(C,z)+(C\cdot L)_z\right\}.$$
\end{enumerate}

\subsection{Curves on toric surfaces}\label{sec-cts}
Given a morphism $\bn:\widehat C\to\Tor(P_\Delta)$ of a smooth curve $\widehat C$
to $\Tor(P_\Delta)$, we denote by
\begin{itemize}
\item $\bn_*(\widehat C)$ the scheme-theoretic push-forward, {\it i.e.},
the one-dimensional part of the image, whose components are taken with the corresponding multiplicities;
\item $C:=\bn(\widehat C)$ the reduced model of $\bn_*(\widehat C)$,
where all components are taken with multiplicity one; we split $\bn$ into $\bn:\widehat C\twoheadrightarrow C\overset{\ina}{\hookrightarrow}\Tor(P_\Delta)$;
\item $\bn^*(D)$, where $D\subset\Tor(P_\Delta)$ is a divisor intersecting $C$ in finitely many points,
the divisor on $\widehat C$ which is the pull-back of the scheme-theoretic
intersection $D\cap\bn_*(\widehat C)$.
\end{itemize}

Denote by ${\mathcal M}_{g,m}(\Tor(P_\Delta),{\mathcal L}_{P_\Delta})$
the moduli space of stable maps $\bn:(\widehat C,\bp)\to\Tor(P_\Delta)$,
where $\widehat C$ is a smooth connected
curve of genus $g$ and
$\bp$ is a sequence of $m$ distinct points in $\widehat C$,
such that $\bn_*(\widehat C)\in|{\mathcal L}_{P_\Delta}|$.
Correspondingly,
we denote
by $\overline{\mathcal M}_{g,m}(\Tor(P_\Delta),{\mathcal L}_{P_\Delta})$
the Kontsevich
compactification of the above moduli space;
this compactification is obtained by adding isomorphism classes of stable maps $\bn:(\widehat C,\bp)\to \Tor(P_\Delta)$,
where $\widehat C$ is a connected nodal curve of arithmetic genus $g$ and
$\bp$
is a sequence of $m$ distinct points
in ${\widehat C} \setminus \Sing({\widehat C})$,
such that $\bn_*(\widehat C)\in|{\mathcal L}_{P_\Delta}|$
(here $\Sing({\widehat C})$ stands for the set of singularities of $\widehat C$).
In what follows, we work with certain subspaces of these moduli spaces specified for genus $g=0$
in Section \ref{sec2.3}, for genus $g=1$ in Section \ref{asec3}, and for genus $2$ in Section \ref{sec-g2}.
For any subset ${\mathcal M} \subset {\mathcal M}_{g,m}(\Tor(P_\Delta),{\mathcal L}_{P_\Delta})$,
denote by $\overline{\mathcal M}$ the closure of $\mathcal M$ in the
compactified moduli space $\overline{\mathcal M}_{g,m}(\Tor(P_\Delta),{\mathcal L}_{P_\Delta})$.

We denote by $(C_1\cdot C_2)_z$ the intersection multiplicity of
curves $C_1,C_2\subset\Tor(P_\Delta)$ at a smooth point $z$ of $\Tor(P_\Delta)$.
By $C_1C_2$
we mean
the total intersection multiplicity of the curves $C_1,C_2\subset\Tor(P_\Delta)$
provided that all intersection points of these curves are smooth in $\Tor(P_\Delta)$.

Recall
that deformations of a morphism $\bn:\widehat C\to \Tor(P_\Delta)$
of a smooth curve $\widehat C$ are encoded by the normal sheaf on $\widehat C$ (see, for instance, \cite{GK}):
$${\mathcal N}_\bn:=\bn^*{\mathcal T}\Tor(P_\Delta)/{\mathcal T}\widehat C$$
(where ${\mathcal T}$ denotes the tangent bundle).
In the case of an immersion, ${\mathcal N}_\bn$ is a line bundle of degree
\begin{equation}\deg{\mathcal N}_\bn=c_1(\Tor(P_\Delta))c_1({\mathcal L}_{P_\Delta})-2.\label{eis4}\end{equation}
In case of $\bn:\widehat C\to C\hookrightarrow\Tor(P_\Delta)$ birational onto its image $C$, but not necessarily immersion,
we
consider deformations realizing specific local deformations of $C$
at its singular points (see details in Section \ref{sec-def}).

Let $\Delta=(\ba^\sigma_i,\ i=1,...,n^\sigma,\ \sigma\in P_\Delta^1)$ be a toric degree as introduced in Section \ref{sec-td-m}.
Fix a non-negative integer $g\le\cI(\Delta)$ and a non-negative integer $n(\sigma) \leq n^\sigma$
for each edge $\sigma \in P^1_\Delta$
such that $\sum_{\sigma\in P^1_\Delta} n(\sigma) < \sum_{\sigma\in P^1_\Delta} n^\sigma$.
Put $n_\partial = \sum_{\sigma\in P^1_\Delta} n(\sigma)$ and $n_{\ina} = n - n_\partial$,
where
$$n = \sum_{\sigma\in P^1_\Delta}n^\sigma + g - 1.$$
Choose a sequence $\bw$ of $n$ distinct points in $\Tor(P_\Delta)$ splitting into two subsequences:
\begin{itemize}
\item $\bw_\partial$ consisting of $n_\partial$ points in general position on $\Tor(\partial P_\Delta)$ so that,
for each side $\sigma\in P^1_\Delta$, the sequence $\bw_\partial$ contains $n(\sigma)$ points
$\bw^\sigma_i$, $1 \leq i \leq n(\sigma)$, lying on $\Tor(\sigma)$;
\item
$\bw_{\ina}$ consisting of $n_{\ina}$ points $w_j$,
$j = 1$, $\ldots$, $n_{\ina}$ in general position in $(\C^\times)^2\subset\Tor(P_\Delta)$.
\end{itemize}

Introduce the subset ${\mathcal M}_{g,n}(\Delta,\bw)\subset{\mathcal M}_{g,n}(\Tor(P_\Delta),{\mathcal L}_{P_\Delta})$
consisting of the elements $\left[\bn:(\widehat C,\bp)\to\Tor(P_\Delta)\right]$
satisfying the following assumptions:
\begin{enumerate}
\item[(M1)]
the sequence $\bp$ splits into two disjoint subsequences
$\bp_\partial$ and $\bp_{\ina}$ containing $n_\partial$ points $p^\sigma_i$, $\sigma \in P^1_\Delta$, $1 \leq i \leq n(\sigma)$,
and $n_{\ina}$ points $p_j$, 
$j = 1$, $\ldots$, $n_{\ina}$, respectively,
such that
$\bn(\bp_\partial)=\bw_\partial$
(more precisely, for each $\sigma \in P^1_\Delta$, one has $\bn(p^\sigma_i) = w^\sigma_i$, $i = 1$, $\ldots$, $n(\sigma)$)
and
$\bn(p_j) = w_j$, $j = 1$, $\ldots$, $n_{\ina}$;
\item[(M2)]
for each $\sigma \in P^1_\Delta$, one has
$$\bn^*(\Tor(\sigma))=\sum_{i=1}^{n^\sigma}k_i^\sigma p_i^\sigma\in\Div(\widehat C),$$
where
$\bn(p_i^\sigma)=w_i^\sigma$ for all $p^\sigma_i\in\bp_\partial$ and $w^\sigma_i\in\bw_\partial$,
$\sigma\in P^1_\Delta$, $1\le i\le n(\sigma)$,
and
$p^\sigma_i$, $n(\sigma)<i\le n^\sigma$, are some points of $\widehat C$ not counted in the sequence $\bp_\partial$.; at last, $\Div(\widehat C)$ denotes the 
group of Weil divisors on $\widehat C$.
\end{enumerate}

\begin{lemma}\label{l3-1a}
Suppose that the sequence $\bw$ is in general position subject to the location with respect to the toric divisors as indicated above. Then,
the space ${\mathcal M}_{g,n}(\Delta,\bw)$ is finite.
Moreover, for each element $\left[\bn:(\widehat C,\bp)\to\Tor(P_\Delta)\right]\in
{\mathcal M}_{g,n}(\Delta,\bw)$, the map $\bn$ takes $\widehat C$ birationally onto
an immersed curve $C=\bn(\widehat C)$ that,
for each $\sigma \in P^1_\Delta$,
intersects the toric divisor $\Tor(\sigma)$ at $n^\sigma$ distinct points and $C$ is smooth at each of these intersection points. Furthermore, $\bn(\bp)\cap\Sing(C)=\emptyset$.
\end{lemma}

{\bf Proof.}
Consider the family
$${\mathcal M}'_{g,n}(\Delta)=\big\{[\bn:(\widehat C,\bp_\partial)\to\Tor(P_\Delta)]\in{\mathcal M}_{g,n}(\Tor(P_\Delta),{\mathcal L}_{P_\Delta})$$
such that
$$\bn:\widehat C\to C=\bn(\widehat C)\ \text{is birational},\qquad\qquad$$
\begin{equation}\qquad\quad[\bn:(\widehat C,\bp_\partial)\to\Tor(P_\Delta)]\ \text{satisfies condition (M2)}\big\}.
\label{efamily}\end{equation}
We claim that
\begin{equation}\dim{\mathcal M}'_{g,n}(\Delta)\le n.\label{edop-alg1}\end{equation}
Assuming (\ref{edop-alg1}) and imposing the conditions $\bn(p^\sigma_i)=w^\sigma_i\in\bw_\partial$,
$1\le i\le n(\sigma)$, for every edge $\sigma\in P^1_\Delta$,
as well as
imposing
extra $n_{\ina}$ point constraints $\bw_{\ina}$ in general position in $\Tor(P_\Delta)^\times$, we obtain a finite set.

To prove (\ref{edop-alg1}), take a generic element $\xi=[\bn:(\widehat C,\bp)\to\Tor(P_\Delta)]$ in an irreducible component
of ${\mathcal M}'_{g,n}(\Delta)$. It is a smooth point of that component, and the germ of ${\mathcal M}'_{g,n}(\Delta)$ at $\xi$
realizes an equisingular deformation. We
estimate the dimension of the tangent space to this germ.
Put $C=\bn(\widehat C)$. Taking into account the formulas in Section \ref{sec-def},
we derive that
\begin{equation}\dim{\mathcal M}'_{g,n}(\Delta)\le h^0\left(C,{\mathcal J}_{Z/C}\otimes_{{\mathcal O}_C}\ina^*{\mathcal L}_{P_\Delta}\right)
=h^0\left(\widehat C,{\mathcal O}_{\widehat C}(-\bd)\otimes_{{\mathcal O}_{\widehat C}}\bn^*{\mathcal L}_{P_\Delta}\right),\label{edop-alg2}\end{equation}
where \begin{itemize}\item $Z\subset C$ is the zero-dimensional scheme
supported at $\Sing(C)\cup(C\cap\Tor(\partial P_\Delta))$,
defined at each point $z\in\Sing(C)\cap\Tor(P_\Delta)^\times$ by the conductor ideal $I^{cond}(C,z)$, and defined
at each point $w\in C\cap\Tor(\sigma)$, $\sigma\in P^1_\Delta$, by the ideal\
described by the right-hand side
of the first formula of item (b) in Section \ref{sec-def};
\item ${\mathcal J}_{Z/C}\subset{\mathcal O}_C$ is the ideal sheaf of the scheme $Z\subset C$;
\item the divisor $\bd$ satisfies $$\deg \bd=2\sum_{z\in\Sing(C)}\delta(C,z)
+\sum_{\sigma\in P^1_\Delta}\sum_{w\in C\cap\Tor(\sigma)}\sum_{Q\subset(C,w)}((Q\cdot\Tor(\sigma))_w-1)$$
$$+\sum_{z\in\Sing(C)\cap\Tor(P_\Delta)^\times}\sum_{Q\subset(C,z)}(\mt Q-1).$$
\end{itemize} It follows that
$$\deg\left({\mathcal O}_{\widehat C}(-\bd)\otimes_{{\mathcal O}_{\widehat C}}\bn^*{\mathcal L}_{P_\Delta}\right)=C^2-\deg \bd$$
$$=C^2-(C^2+CK_{\Tor(P_\Delta)}+2-2g)$$
$$+CK_{\Tor(P_\Delta)}+\#\{Q\subset (C,w),\ w\in C\cap\Tor(\sigma)\}$$
$$-\sum_{z\in\Sing(C)\cap\Tor(P_\Delta)^\times}\sum_{Q\subset(C,z)}(\mt Q-1)$$
$$=(2g-2)+\#\{Q\subset (C,w),\ w\in C\cap\Tor(\sigma)\}$$
\begin{equation}-\sum_{z\in\Sing(C)\cap\Tor(P_\Delta)^\times}\sum_{Q\subset(C,z)}(\mt Q-1),\label{edop-alg4}\end{equation}
where the symbol $\#$ stands for the number of elements of a given set.
Since
\begin{equation}3\le\#\{Q\subset (C,w),\ w\in C\cap\Tor(\sigma)\}\le\#\Delta,
\label{edop-alg3}\end{equation}
by (\ref{edop-alg4}) and by the Riemann-Roch formula,
the right-hand side of (\ref{edop-alg2}) does not exceed
$$2g-2+\#\big\{Q\subset (C,w),\ w\in C\cap\Tor(\sigma)\big\}-g+1$$
$$=g-1+\#\{Q\subset (C,w),\ w\in C\cap\Tor(\sigma)\}\le g-1+\#\Delta=n.$$
Furthermore, if either $\#\{Q\subset (C,w),\ w\in C\cap\Tor(\sigma)\}<\#\Delta$,
or there is a singular local branch $Q$ of $C$ centered in $\Tor(P_\Delta)^\times$, then we get ${\mathcal M}'_{g,n}(\Delta,\bw)=\emptyset$.
Thus, we are left with the following case:
$\dim{\mathcal M}'_{g,n}(\Delta)=n$,
the curve $C$ is immersed in $\Tor(P_\Delta)^\times$,
and $\#\{Q\subset (C,w),\ w\in C\cap\Tor(\sigma)\}=\#\Delta$.
Now we show that $C$ is smooth along the toric divisors.
Suppose that $C$ has a local singular branch $Q_0$ centered at some point $w\in\Tor(\sigma)$, $\sigma\in P^1_\Delta$. Fixing the position of the point $w$, we obtain a subfamily ${\mathcal V}_1\subset{\mathcal M}'_{g,n}(\Delta)$ of dimension $\ge n-1$. On the other hand, the computations as above yield that
$$\dim{\mathcal V}_1\le h^0\left(\widehat C,{\mathcal O}_{\widehat C}(-\bd_1)\otimes_{{\mathcal O}_{\widehat C}}\bn^*{\mathcal L}_{P_\Delta}\right),$$
where by the second formula of item (b) in Section \ref{sec-def},
$$\deg \bd_1=2\sum_{z\in\Sing(C)}\delta(C,z)
+\sum_{\sigma\in P^1_\Delta}\sum_{w\in C\cap\Tor(\sigma)}\sum_{Q\subset(C,w)}((Q\cdot\Tor(\sigma))_w-1)+\mt Q_0.$$
This implies that
$$\deg\left(\widehat C,{\mathcal O}_{\widehat C}(-\bd_1)\otimes_{{\mathcal O}_{\widehat C}}\bn^*{\mathcal L}_{P_\Delta}\right)=2g-2+\#\Delta-\mt Q_0\le2g-2+\#\Delta-2,$$
and by the Riemann-Roch formula,
we obtain
$$\dim{\mathcal V}_1\le\#\Delta+g-1-2=n-2,$$
which is a contradiction. Thus, $C$ must be immersed.
Suppose now that some point $w\in C\cap\Tor(\sigma)$, $\sigma\in P^1_\Delta$, is a center of local branches
$Q_1$, \ldots, $Q_r$,
$r\ge2$.
Again, fixing the position of the point $w$, we obtain a subfamily ${\mathcal V}_2\subset{\mathcal M}'_{g,n}(\Delta)$
of dimension $\ge n-1$. On the other hand, the computations as above yield that
$$\dim{\mathcal V}_2\le h^0\left(\widehat C,{\mathcal O}_{\widehat C}(-\bd_2)\otimes_{{\mathcal O}_{\widehat C}}\bn^*{\mathcal L}_{P_\Delta}\right),$$
where
$$\deg \bd_2=2\sum_{z\in\Sing(C)}\delta(C,z)
+\sum_{\sigma\in P^1_\Delta}\sum_{w\in C\cap\Tor(\sigma)}\sum_{Q\subset(C,w)}((Q\cdot\Tor(\sigma))_w-1) + r,$$
which implies
$$\deg\left(\widehat C,{\mathcal O}_{\widehat C}(-\bd_2)\otimes_{{\mathcal O}_{\widehat C}}\bn^*{\mathcal L}_{P_\Delta}\right)=2g-2+\#\Delta-r\le2g-2+\#\Delta-2,$$
and, by the Riemann-Roch formula,
we obtain
$$\dim{\mathcal V}_2\le\#\Delta+g-1-2=n-2,$$
which is a contradiction. Hence, $C$ is smooth along the toric divisors.

It remains to show that ${\mathcal M}_{g,n}(\Delta,\bw)$ does not contain elements $[\bn:(\widehat C,\bp)\to\Tor(P_\Delta)]$ such that $\bn:\widehat C\to\Tor(P_\Delta)$ is a multiple covering of the image $C'=\bn(\widehat C)$. Following the preceding argument, we assume that $\dim{\mathcal M}_{g,n}(\Delta)=n$ and that a generic element of some component ${\mathcal V}_0\subset{\mathcal M}_{g,n}(\Delta)$ is represented by a multiple covering $\bn:\widehat C\to C'\hookrightarrow\Tor(P_\Delta)$,
and then come to a contradiction.
Denote by $\bn':\widehat C'\to C'$ the normalization.
Notice that the genus $g'$ of $\widehat C'$ is at most $g$,
and that the number of local branches of $C'$ centered on the toric divisors is bounded from above by $\#\Delta$. Further on, we can assume that the variation of $\bn:\widehat C\to\Tor(P_\Delta)$ along ${\mathcal V}_0$ induces an equisingular deformation of $C'$ in the linear system $|C'|$.
The computations as in the preceding paragraph imply that
$$n={\mathcal V}_0\le\#\{Q\subset (C',w),\ w\in C'\cap\bigcup_\sigma\Tor(\sigma\}+g'-1,$$
which, in turn, yields that
$$\#\{Q\subset (C',w),\ w\in C'\cap\bigcup_\sigma\Tor(\sigma\}=\#\Delta\quad\text{and}\quad g'=g.$$
Moreover, the argument of the preceding paragraph would imply that $C'$ is immersed
and smooth along the toric divisors.
That means that $\bn:\widehat C\to \widehat C'$ is ramified at each of the $\#\Delta\ge3$ points
of $\widehat C'$ mapped to the toric divisors, and the ramification index at each of them
is maximal, but this contradicts the Riemann-Hurwitz formula.
\proofend

A map $\bn:(\widehat C,\bp)\to\Tor(P_\Delta)$ such that $[\bn:(\widehat C,\bp)\to\Tor(P_\Delta)]\in{\mathcal M}_{g,n}(\Delta,\bw)$
is said to be {\it real} if
\begin{enumerate}
\item[(i)] the sequence $\bw$ is invariant with respect to the tautological real structure on $\Tor(P_\Delta)$,
\item[(ii)] $(\widehat C,\bp)$ is equipped with a real structure, and
\item[(iii)] $\bn:(\widehat C,\bp)\to\Tor(P_\Delta)$ commutes with the real structures in the source and in the target.
\end{enumerate}
The set of equivalence classes of such real
maps $\bn:(\widehat C,\bp)\to\Tor(P_\Delta)$
considered
up to equivariant isomorphism is denoted by ${\mathcal M}^\R_{g,n}(\Delta,\bw)$.
We say that an element $[\bn:(\widehat C,\bp)\to\Tor(P_\Delta)]\in{\mathcal M}^\R_{g,n}(\Delta,\bw)$ is {\it separating} if the complement in $\widehat C$
to the real point set $\R\widehat C$ is disconnected, {\it i.e.}, consists of two connected components.
The choice of one of these
halves induces the so-called {\it complex orientation} on $\R\widehat C$,
as well as on $\bn(\R\widehat C)$
if $\bn$ is birational onto its image.
Denote by
$\overrightarrow{\mathcal M}_{g,n}^{\R}(\Delta,\bw)$ the set
of separating elements
$[\bn: (\widehat C, \bp) \to \Tor(P_\Delta)] \in {\mathcal M}^\R_{g,n}(\Delta,\bw)$ equipped with
a choice of a half
$\widehat C_+$ of $\widehat C\setminus\R\widehat C$.

Following \cite{Mir},
we assign to each element
$\xi=([\bn:(\widehat C, \bp)\to\Tor(P_\Delta)], \widehat C_+)\in\overrightarrow{\mathcal M}_{g,n}^{\R}(\Delta,\bw)$
its {\it quantum index}
\begin{equation}
\QI(\xi)
=\frac{1}{\pi^2}\int_{\widehat C_+}\bn^*\left(\frac{dx_1\wedge dx_2}{x_1x_2}\right),\quad x_1=|z_1|,\ x_2=|z_2|,\label{eis1}
\end{equation}
with $z_1,z_2$ coordinates in the torus $(\C^\times)^2\simeq
\Tor(P_\Delta)^\times$ such that the form $dx_1\wedge dx_2$ agrees with the orientation of $\Tor_\R^+(P_\Delta)$ defined in Section \ref{sec2.1.1}.
By \cite[Theorem 3.1]{Mir}, if all intersection points of $\bn(\widehat C)$ with the toric divisors are real,
then
$\QI(\xi)\in\frac{1}{2}\Z$ and $|\QI(\xi)| \leq \cA(\Delta)$.

\section{Real refined enumerative invariants in genus $0$, $1$, and $2$}\label{sec-inv}

\subsection{Preliminary remarks}
Recall
that the Welschinger signed count of real plane curves of a given degree and genus, passing through a generic conjugation-invariant configuration of points in the plane, becomes non-invariant as long as
the genus is positive \cite[Theorem 3.1]{IKS0}.
However, in the case of real rational surfaces
with a disconnected real
part,
Welschinger invariants survive for certain positive genera $g$, provided
that some conditions on constraints are fulfilled: in particular, each real connected component
of the ambient surface contains at most one {\it oval} ({\it i.e.}, global real branch) of any
curve under consideration, and the counted curves have the maximal possible number $g + 1$ of ovals; see  \cite{Sh}.
We exploit here a similar approach in order to extend the invariant count of
real rational curves \cite{Mir} to positive genera.
Our idea is to consider real curves that are tangent
to toric divisors with even intersection multiplicity at each intersection point.
In such a case, every quadrant of a real toric surface plays the role of a separate component
of the real part of a real surface.
The lack of the ``bad" wall in the space of constraints indicated in \cite[Section 3]{IKS0}
is expressed as a cohomology vanishing condition which we prove below
in Lemmas \ref{lis2}, \ref{lis3}, and \ref{lis5}.

\begin{remark}\label{ris}
The above idea does not work for genus $g\ge3$: if we fix constraints on the boundary
of the positive quadrant {\rm (}see Sections \ref{sec2.3}, \ref{asec3}, and \ref{sec-g2}{\rm )} and $g$
more extra points in the big real torus, the set of counted curves will necessarily
include non-dividing real curves with at most $g$ ovals, which does not allow one to perform Mikhalkin's quantization.
\end{remark}

\subsection{Refined rational invariants}\label{sec2.3}
In the notation of Section \ref{sec-cts}, assume that
$$g=0,\quad n_{in}=0,\quad n_\partial=n=\sum_{\sigma\in P^1_\Delta}n^\sigma,\quad\text{and}\quad\bw=\bw_\partial\subset\partial\Tor_\R^+(P_\Delta).$$
In particular, for all $\sigma\in P^1_\Delta$ but one,
we have $n(\sigma)=n^{\sigma}$, and for the remaining edge $\sigma'$,
we have $n(\sigma')=n^{\sigma'}-1$. Note that for a given $\bw$, there exists a unique point $w^{\sigma'}_{n^{\sigma'}}\in \partial\Tor_\R^+(P_\Delta)$ such that the sequence $\widehat\bw:=\bw\cup\{w^{\sigma'}_{n^{\sigma'}}\}$
belongs to the hypersurface $\Men(\Delta)$.
We denote by ${\mathcal M}_{0,n}(\Delta,\widehat\bw)\subset{\mathcal M}_{0,n}(\Delta,\bw)$ the subset formed by the curves passing through $\widehat\bw$.

Denote by $\Men_0(\Delta)$ the subset of $\Men(\Delta)$ formed by the above sequences $\widehat\bw$ such that ${\mathcal M}_{0,n}(\Delta,\widehat\bw)$ satisfies the conclusions of Lemma \ref{l3-1a}.
The closure
$$\overline{\bigcup_{\widehat\bw\in \Men_0(\Delta)}{\mathcal M}^\R_{0,n}(\Delta,\widehat\bw)}
\subset\overline{\mathcal M}^\R_{0,n}(\Tor(P_\Delta),{\mathcal L}_{P_\Delta})$$ is naturally fibered over the space of sequences $\widehat\bw\in\Men(\Delta)$,
and for the sake of notation,
we denote the fibers by $\overline{\mathcal M}^\R_{0,n}(\Delta,\widehat\bw)$. 

Consider the subset ${\mathcal M}^{\R,+}_{0,n}(\Delta,\widehat\bw)\subset{\mathcal M}^\R_{0,n}(\Delta,\widehat\bw)$,
formed by the elements $[\bn:(\PP^1,\bp)\to\Tor(P_\Delta)]$ such that $\bn(\PP_\R^1)\subset
\Tor_\R^+(P_\Delta)$\ \footnote{Note that ${\mathcal M}^{\R,+}_{0,n}(\Delta,\widehat\bw)\subsetneq{\mathcal M}^\R_{0,n}(\Delta,\widehat\bw)$ when there exists a connected component of
$\Tor_\R(P_\Delta)^\times\setminus\Tor_\R(P_\Delta)^+$
such that the boundary of the component coincides with $\partial\Tor_\R^+(P_\Delta)$.}.
For every $[\bn:(\PP^1,\bp)\to\Tor(P_\Delta)]\in{\mathcal M}^{\R,+}_{0,n}(\Delta,\widehat\bw)$,
the curve $C=\bn(\PP^1)$ intersects each toric divisor $\Tor(\sigma)$ in points $w^\sigma_i$, $1\le i\le n^\sigma$.
Hence, each element $\xi$ of the set
\begin{align*}\overrightarrow{\mathcal M}_{0,n}^{\R,+}(\Delta,\widehat\bw)=&\{([\bn:(\PP^1,\bp)\to\Tor(P_\Delta)], \PP^1_+) \in\overrightarrow{\mathcal M}_{0,n}^{\R}(\Delta,\widehat\bw)\ :\\ &\qquad\qquad[\bn:(\PP^1,\bp)\to\Tor(P_\Del)]\in{\mathcal M}_{0,n}^{\R,+}(\Delta,\widehat\bw)\}\end{align*}
possesses a quantum index $\QI(\xi)\in\frac{1}{2}\Z$.

Furthermore, if $\widehat\bw\in \Men_0(\Delta)$, then
there are well-defined
{\it Welschinger signs}
\begin{equation}
W_0(\xi)=(-1)^{e_+(C)}\cdot
\prod_{\renewcommand{\arraystretch}{0.6}
\begin{array}{c}
\scriptstyle{\sigma\in P^1_\Delta,\ 1\le i\le n^\sigma}\\
\scriptstyle{k_i^\sigma \ \equiv \ 0  \mod 2}\end{array}}\sign(\xi,w_i^\sigma)\ ,
\label{ae13}
\end{equation}
\begin{equation}{\widetilde W}_0(\xi)=(-1)^{e(C)}\prod_{\sigma\in P^1_\Delta,\ 1\le i\le n^\sigma}\sign(\xi,w_i^\sigma)\ ,\label{ae9i}\end{equation}
where $e_+(C)$, respectively, $e(C)$, is the number of elliptic nodes ({\it i.e.}, real intersection points of a pair of smooth transversal complex conjugate local branches)
of $C = \bn(\PP^1)$ 
that arise from a real nodal equigeneric deformation of all singular points of $C$
in the positive quadrant $\Tor^+_\R(P_\Delta)$, respectively, in the torus $\Tor_\R(P_\Delta)^\times$
({\it cf}. \cite[Section 1.1]{IKS4}),
and
$\sign(\xi,w_i^\sigma)$ equals $1$ or $-1$ according as the complex orientation of $\R C = \bn(\PP_\R^1)$
at $w_i^\sigma$ agrees or not with the
fixed orientation of $\partial \Tor^+_\R(P_\Delta)$.

For each $\kappa \in \frac{1}{2}\Z$ such that $|\kappa| \leq \cA(\Delta)$, put
$$W_0^\kappa(\Delta,\widehat\bw)=
\sum_{\renewcommand{\arraystretch}{0.6}
\begin{array}{c}
\scriptstyle{\xi\in\overrightarrow{\mathcal M}_{0,n}^{\R,+}(\Delta,\widehat\bw)}\\
\scriptstyle{\QI(\xi)=\kappa}\end{array}}W_0(\xi),\quad \widetilde W_0^\kappa(\Delta,\widehat\bw)=
\sum_{\renewcommand{\arraystretch}{0.6}
\begin{array}{c}
\scriptstyle{\xi\in\overrightarrow{\mathcal M}_{0,n}^{\R,+}(\Delta,\widehat\bw)}\\
\scriptstyle{\QI(\xi)=\kappa}\end{array}}\widetilde W_0(\xi)\ .$$

\begin{theorem}\label{at1}
For each $\kappa \in \frac{1}{2}\Z$ such that $|\kappa| \leq \cA(\Delta)$,
the values $W_0^\kappa(\Delta,\widehat\bw)$ and $\widetilde W_0^\kappa(\Delta,\widehat\bw)$ do not depend on
the choice of a
sequence $\widehat\bw\in\Men_0(\Delta)$.
\end{theorem}

We prove this theorem in Section \ref{asec2}.
Since the values $W_0^\kappa(\Delta,\widehat\bw)$ and $\widetilde W_0^\kappa(\Delta,\widehat\bw)$ do not depend on
the choice of a
sequence $\widehat\bw\in\Men_0(\Delta)$, we use for these values the notations
$W_0^\kappa(\Delta)$ and $\widetilde W_0^\kappa(\Delta)$, respectively.

\begin{remark}\label{ar1}
(1) Theorem \ref{at1} slightly generalizes the statement of \cite[Theorem 5]{Mir},
because we allow tangencies of any even intersection multiplicity.

(2) One can show {\rm (}\cite{Mir, IS}{\rm )}
that the invariants $W_0^\kappa(\Delta)$ and $\widetilde W^\kappa_0(\Delta)$
vanish for all $\kappa\not\in2\Z$ and for all $\kappa\in2\Z$ such that $\kappa\not\equiv\cA(\Delta)\mod4$. In general, the invariants $W_0^\kappa(\Delta)$
can be computed via the formula in \cite[Theorem 3.4]{Blomme}.
\end{remark}

\subsection{Refined elliptic invariants}\label{asec3}
In the notation of Section \ref{sec-cts}, assume that
$$g=1,\quad n_{in}=1,\quad n_\partial=\sum_{\sigma\in P^1_\Delta}n^\sigma-1,\quad\bw_\partial\subset\partial\Tor_\R^+(P_\Delta),\quad
\bw_{in}=\{w_0\}\subset {\mathcal Q},$$
where ${\mathcal Q}$ is an open quadrant in $\Tor_\R(P_\Delta)^\times\setminus\Tor^+_\R(P_\Delta)$.
As in the preceding section, for all $\sigma\in P^1_\Delta$ but one, we have $n(\sigma)=n^{\sigma}$,
and for the remaining edge $\sigma'$, we have $n(\sigma')=n^{\sigma'}-1$. As in Section \ref{sec2.3},
we define the completion $\widehat\bw_\partial:=\bw_\partial\cup\{w^{\sigma'}_{n^{\sigma'}}\}\in\Men(\Delta)$,
where $w^{\sigma'}_{n^{\sigma'}}\in \partial\Tor_\R^+(P_\Delta)$, and put $\widehat\bw=\widehat\bw_\partial\cup\bw_{in}$.
We denote by ${\mathcal M}_{1,n}(\Delta,\widehat\bw)\subset{\mathcal M}_{1,n}(\Delta,\bw)$
the subset formed by the curves passing through $\widehat\bw$.

Denote by $\Men_1(\Delta)$ the set of sequences $\widehat\bw\in \Men(\Delta)\times {\mathcal Q}$ such that
the set ${\mathcal M}_{1,n}(\Delta,\widehat\bw)$ satisfies the conclusions of Lemma \ref{l3-1a}.
The closure
$$\overline{\bigcup_{\widehat\bw\in \Men_1(\Delta)}{\mathcal M}^\R_{1,n}(\Delta,\widehat\bw)}
\subset\overline{\mathcal M}^\R_{1,n}(\Tor(P_\Delta),{\mathcal L}_{P_\Delta})$$ is naturally fibered over the space of sequences $\widehat\bw\subset\Men(\Delta)\times {\mathcal Q}$, and
we denote the fibers by $\overline{\mathcal M}^\R_{1,n}(\Delta,\widehat\bw)$.

Consider the subset ${\mathcal M}^{\R,+}_{1,n}(\Delta,\widehat\bw)\subset{\mathcal M}^\R_{1,n}(\Delta,\widehat\bw)$,
formed by the elements
$$[\bn:(\widehat C,\bp)\to\Tor(P_\Delta)], \;\;\; \text{where} \; \widehat C \; \text{is a smooth elliptic curve},$$
defining real curves $C=\bn(\widehat C)\subset\Tor(P_\Delta)$
that have a one-dimensional real branch in
$\Tor_\R^+(P_\Delta)$ and that intersect
the toric divisor $\Tor(\sigma)$ in $w^\sigma_i$, $1\le i\le n^\sigma$, for all $\sigma\in P^1_\Delta$.
It follows from Lemma \ref{l3-1a} that for $\widehat\bw\in \Men_1(\Delta)$, the curve $C$ has two one-dimensional branches, one in $\Tor^+_\R(P_\Delta)$ and the other in the closed quadrant $\overline{\mathcal Q}$
containing the point $w_0$. This defines a unique real structure on $\widehat C$ whose fixed point set $\R\widehat C$ consists of two ovals, making $\bn:\widehat C\to C$
a separating real curve.
The choice of a component of $\widehat C\setminus
\R{\widehat C}$
defines a complex orientation
on the one-dimensional branches of $\R C$.

We impose a further restriction
on the choice of the quadrant ${\mathcal Q} \ni w_0$ dictated by our wish to avoid specific wall-crossing events
that may occur in the proof of the invariance of the count to be specified below:
first, the escape of a real one-dimensional branch of counted curves out of the positive quadrant
and, second, degenerations in which the union of the toric divisors splits off.

\begin{definition}\label{ad2} We say that the quadrant ${\mathcal Q}\subset \Tor_\R(P_\Delta)^\times$
satisfies the admissible quadrant condition {\rm (}AQC{\rm )} if the closure $\overline {\mathcal Q}$ of ${\mathcal Q}$
shares with $\partial\Tor^+_\R(P_\Delta)$ at most one side.
\end{definition}

\begin{example}\label{ar2}
Let $d$, $d_1$ and $d_2$ be positive integers.
Condition {\rm (}AQC{\rm )}
imposes no restriction on the choice of ${\mathcal Q}$ in the case of the triangle $\conv\{(0,0),(d,0),(0,d)\}$
{\rm (}the convex hull of the points $(0, 0)$, $(d, 0)$, $(0, d)${\rm )}
having the projective plane $\PP^2$ as associated toric surface, while for the rectangle
$\conv\{(0,0), (d_1,0),(d_1,d_2), (0,d_2)\}$
which has
$\PP^1 \times \PP^1$ as associated toric surface,
the condition authorizes only the quadrant ${\mathcal Q}$,
where the both coordinates are negative.
\end{example}

Denote by
$\overrightarrow{\mathcal M}_{1,n}^{\R,+}(\Delta,\widehat\bw)$ the subset
of $\overrightarrow{\mathcal M}_{1,n}^{\R}(\Delta,\widehat\bw)$ formed by the elements
$$([\bn:(\widehat C,\bp)\to\Tor(P_\Delta)], \widehat C_+)$$
such that $[\bn:(\widehat C,\bp)\to\Tor(P_\Delta)]\in{\mathcal M}_{1,n}^{\R,+}(\Delta,\widehat\bw)$ and $\widehat C_+$
is a half of $\widehat C\setminus\R\widehat C$.
According to \cite[Theorem 1]{Mir}, each element $\xi$ of $\overrightarrow{\mathcal M}_{1,n}^{\R,+}(\Delta,\widehat\bw)$
has a quantum index $\QI(\xi)\in\frac{1}{2}\Z$, $|\QI(\xi)|\le\cA(\Delta)$ (defined by formula (\ref{eis1})
with the integration domain $\widehat C_+$),
and for $\bw$ such that $\widehat\bw\in \Men_1(\Delta)$, according to the argument in \cite[Section 1.1]{IKS4},
there is a well-defined {\it Welschinger sign}
\begin{equation}W_1(\xi)=(-1)^{e_+(C)}
\cdot(-1)^{h(C,{\mathcal Q})}\cdot\prod_{\renewcommand{\arraystretch}{0.6}
\begin{array}{c}
\scriptstyle{\sigma\in P^1_\Delta,\ 1\le i\le n^\sigma}\\
\scriptstyle{k_i^\sigma \ \equiv \ 0  \mod 2}\end{array}}
\sign(\xi,w_i^\sigma)\ ,\label{ae21}\end{equation}
where $h(C,{\mathcal Q})$ is the number of hyperbolic nodes of $C = \bn(\widehat C)$
in the quadrant ${\mathcal Q}$ that arise in any real nodal equigeneric deformation of all singularities of $C$ in ${\mathcal Q}$
(here $e_0(C)$ and $\sign(\xi, w_i^\sigma)$ are straightforward analogs of the corresponding ingredients
in formula (\ref{ae13})).

For each $\kappa \in \frac{1}{2}\Z$ such that $|\kappa| \leq \cA(\Delta)$, put
\begin{equation}W_1^\kappa(\Delta,\widehat\bw)=
\sum_{\renewcommand{\arraystretch}{0.6}
\begin{array}{c}
\scriptstyle{\xi\in\overrightarrow{\mathcal M}_{1,n}^{\R,+}(\Delta,\widehat\bw)}\\
\scriptstyle{\QI(\xi)=\kappa}\end{array}}W_1(\xi),\quad \widetilde W_1^\kappa(\Delta,\widehat\bw)=
\sum_{\renewcommand{\arraystretch}{0.6}
\begin{array}{c}
\scriptstyle{\xi\in\overrightarrow{\mathcal M}_{1,n}^{\R,+}(\Delta,\widehat\bw)}\\
\scriptstyle{\QI(\xi)=\kappa}\end{array}}\widetilde W_1(\xi)\ ,\label{f-elliptic}\end{equation}
where ({\it cf}. formula (\ref{ae9i}))
$$\widetilde W_1(\xi)=(-1)^{e(C)}\prod_{\sigma\in P^1_\Delta,\ 1\le i\le n^\sigma}\sign(\xi,w_i^\sigma)\ .$$

\begin{theorem}\label{at2} Suppose that $n\ge4$.
Then,
the values $W_1^\kappa(\Delta,\widehat\bw)$ and $\widetilde W^\kappa_1(\Delta,\widehat\bw)$ do not depend
on the choice of a
sequence
$\widehat\bw\in\Men_1(\Delta)\times{\mathcal Q}$
subject to the following conditions:
$\widehat\bw_\partial\subset\partial\Tor_\R^+(P_\Delta)$ and $w_0\in {\mathcal Q}$, where ${\mathcal Q}$ is a fixed non-positive quadrant satisfying {\rm (}AQC{\rm )}.
\end{theorem}

We prove this theorem in Section \ref{asec5}.
Since the values $W_1^\kappa(\Delta,\widehat\bw)$ and $\widetilde W_1^\kappa(\Delta,\widehat\bw)$ do not depend on
the choice of a
sequence $\widehat\bw\in\Men_1(\Delta)\times{\mathcal Q}$, we use for these values the notations
$W_1^\kappa(\Delta, {\mathcal Q})$ and $\widetilde W_1^\kappa(\Delta, {\mathcal Q})$, respectively.

\begin{remark}\label{rat2a} (1) We expect that the statement of Theorem \ref{at2} holds for the remaining value $n=3$ too. It corresponds to $\#\Delta=3$ and
$P_\Delta$
a lattice triangle, while the counted elliptic curves in $\Tor(P_\Delta)$ intersect each toric divisor at one point, where they are unibranch.

(2) One can show \cite{IS}
that the invariants $W_1^\kappa(\Delta,{\mathcal Q})$
and $\widetilde W_1^\kappa(\Delta,{\mathcal Q})$
vanish for all $\kappa\not\in2\Z$ and for all $\kappa\in2\Z$
such that $\kappa\not\equiv\cA(\Delta)\mod4$.
\end{remark}

\subsection{Refined invariants of genus two}\label{sec-g2}

For genus $g=2$, we restrict our attention to the case $\Tor(P_\Delta)\simeq\PP^2$. That is, $\Delta$ consists of positive multiples of the vectors $(1,1),(-1,0),(0,-1)$ and defines plane curves of an even degree $2d$. The boundary of the positive quadrant $(\PP^2_\R)^+$ consists of the three segments $D_s^+$, $s=0,1,2$, where $D_s^+=D_s\cap\Cl(\PP_\R^2)^+$, $D_s=\{x_s=0\}\subset\PP^2$ corresponding the sides
$$[(0,d),(d,0)],\ [(0,0),(0,d)],\ [(0,0),(d,0)]\in P^1_\Delta,$$ respectively.
Introduce the completions $\widehat\bw_\partial=\bw_\partial\cup\{w^{(0)}_{n^0}\}\in\Men(\Delta)$ and $\widehat\bw=\widehat\bw_\partial\cup\bw_{in}$, where
$\bw_{in}=\{w_1,w_2\}$,
$$
\displaylines{
w_1\in(\PP_\R^2)^{(-,+)}=\left\{(x_0 : x_1 : x_2) \in \PP_\R^2 \; : \; \frac{x_1}{x_0}<0,\ \frac{x_2}{x_0}>0\right\}, \cr
w_2\in(\PP_\R^2)^{(+,-)}=\left\{(x_0 : x_1 : x_2) \in \PP_\R^2 \; : \; \frac{x_1}{x_0}>0,\ \frac{x_2}{x_0}<0\right\}.
}
$$

Denote by $\Men_2(\Delta)$ the set of sequences $\widehat\bw\in \Men(\Delta)\times(\PP_\R^2)^{(-,+)}\times(\PP_\R^2)^{(+,-)}$ such that
the set ${\mathcal M}_{2,n}(\Delta,\widehat\bw)$ satisfies the conclusions of Lemma \ref{l3-1a}.
The closure
$$\overline{\bigcup_{\widehat\bw\in \Men_2(\Delta)}{\mathcal M}^\R_{2,n}(\Delta,\widehat\bw)}
\subset\overline{\mathcal M}^\R_{2,n}(\Tor(P_\Delta),{\mathcal L}_{P_\Delta})$$ is naturally fibered over the space of sequences $\widehat\bw\subset\Men(\Delta)\times(\PP_\R^2)^{(-,+)}\times(\PP_\R^2)^{(+,-)}$, and
we denote the fibers by $\overline{\mathcal M}^\R_{2,n}(\Delta,\widehat\bw)$. 

Denote by ${\mathcal M}^\R_{2,n}(\Delta,\widehat\bw)$ the space of equivariant isomorphism classes of real stable maps of curves of genus $2$. Due to the choice of the points of $\widehat\bw$, each element $\xi=[\bn:(\widehat C,\bp)\to\PP^2]\in{\mathcal M}^\R_{2,n}(\Delta,\widehat\bw)$, where $\widehat C$ stands for a real Riemann surface of genus $2$, and $\xi$
defines an $M$-curve, {\it i.e.}, $\R\widehat C$ consists of three ovals that are mapped to the quadrants
$\Cl(\PP_\R^2)^+$, $\Cl(\PP_\R^2)^{(-,+)}$, $\Cl(\PP_\R^2)^{(+,-)}$, respectively. The complement $\widehat C\setminus\R\widehat C$ consists of two components, and the choice of one of them, $\widehat C^+$, defines
a complex orientation on $\R\widehat C$ and on $\bn(\R\widehat C)$. Denote by $\overrightarrow{\mathcal M}^\R_{2,n}(\Delta,\widehat\bw)$ the set of elements of ${\mathcal M}^\R_{2,n}(\Delta,\bw)$ equipped with a complex orientation. Similarly to formula (\ref{ae21}), we define for each element $\xi=[\bn:(\widehat C,\bp)\to\PP^2]\in\overrightarrow{\mathcal M}^\R_{2,n}(\Delta,\widehat\bw)$,
$$W_2(\xi)=(-1)^{e_+(C)}\cdot(-1)^{h_-(C)}\cdot\prod_{\sigma\in P^1_\Delta}\prod_{k_i^\sigma\equiv0\mod2}
\sign(\xi,w_i^\sigma),$$ where $C=\bn(\widehat C)$ and
$h_-(C)$ is the number of hyperbolic nodes which appear in real local equigeneric nodal deformations of singular points of $C$ in $(\PP_\R^2)^{(-,+)}\cup(\PP_\R^2)^{(+,-)}$. Put
$$W^\kappa_2(\Delta,\widehat\bw)=\sum_{\renewcommand{\arraystretch}{0.6}
\begin{array}{c}
\scriptstyle{\xi\in\overrightarrow{\mathcal M}^\R_{2,n}(\Delta,\widehat\bw)}\\
\scriptstyle{\QI(\xi)=\kappa}\end{array}}W_2(\xi),\quad \widetilde W^\kappa_2(\Delta,\widehat\bw)=\sum_{\renewcommand{\arraystretch}{0.6}
\begin{array}{c}
\scriptstyle{\xi\in\overrightarrow{\mathcal M}^\R_{2,n}(\Delta,\widehat\bw)}\\
\scriptstyle{\QI(\xi)=\kappa}\end{array}}\widetilde W_2(\xi),$$
where ({\it cf}. formula (\ref{ae9i})
$$\widetilde W_2(\xi)=(-1)^{e(C)}\prod_{\sigma\in P^1_\Delta,\ 1\le i\le n^\sigma}\sign(\xi,w_i^\sigma)\ .$$

\begin{theorem}\label{g2t1} Suppose that $n\ge5$.
For any toric degree $\Delta$ as above and any $\kappa\in\frac{1}{2}\Z$ such that $|\kappa|\le2d^2$, the values $W^\kappa_2(\Delta,\widehat \bw)$ and $\widetilde W^\kappa_2(\Delta,\widehat\bw)$
do not depend on the choice of
a
sequence $\widehat\bw\in\Men_2(\Delta)\times(\PP_\R^2)^{(-,+)}\times(\PP_\R^2)^{(+,-)}$.
\end{theorem}

We prove this theorem in Section \ref{asec6}.
Since the values $W_2^\kappa(\Delta,\widehat\bw)$ and $\widetilde W_2^\kappa(\Delta,\widehat\bw)$ do not depend on
the choice of a
sequence $\widehat\bw\in\Men_2(\Delta)\times(\PP_\R^2)^{(-,+)}\times(\PP_\R^2)^{(+,-)}$,
we use for these values the notations
$W_2^\kappa(\Delta)$ and $\widetilde W_2^\kappa(\Delta)$, respectively.

\begin{remark}\label{rat2b}
We expect that the statement of Theorem \ref{g2t1} holds for the remaining value $n=4$ too. It corresponds to $\#\Delta=3$ and $T_\Delta$ a lattice triangle, while the counted curves of genus $2$ in $\Tor(P_\Delta)$ intersect each toric divisor at one point, where they are unibranch.
\end{remark}

\section{Invariance in genus zero: proof of Theorem \ref{at1}}\label{asec2}

Note that $\Men_0(\Delta)$ is an open dense semialgebraic subset of $\Men^\rho_\R(\Delta)$. We call its connected components the {\it chambers}. The complement to this subset is the union of finitely many semialgebraic strata of codimension $\ge1$. We call the strata of codimension $1$ the {\it walls}.

Now we pick two
sequences $\widehat\bw(0)$ and $\widehat\bw(1)$
satisfying the hypotheses of Theorem \ref{at1}.
We can join
these sequences by a generic path $\{\widehat\bw(t)\}_{0\le t\le 1}\subset\Men_\R^\rho(\Delta)$ which avoids the loci of codimension $>1$.
To prove Theorem \ref{at1}, we need to verify the constancy of $W_0^\kappa(\Delta,\widehat\bw(t))$ and $\widetilde W_0^\kappa(\Delta,\widehat\bw(t))$ inside any given chamber and in all possible wall-crossing events, for all $\kappa$.

\subsection{Invariance inside chambers}\label{chamb0}
Here,
we verify that $W_0^\kappa(\Delta,\widehat\bw(t))$ does not change along intervals $t'<t<t''$
for which $\widehat\bw(t)$ lies in a given chamber.
To this end, we show that the projection
of ${\mathcal M}_{0,n}^{\R,+}(\Delta,\widehat\bw(t))_{t'<t<t''}$
to the interval $(t',t'')$
does not have critical points. In such a case, the family ${\mathcal M}_{0,n}^{\R,+}(\Delta,\widehat\bw(t))_{t'<t<t''}$ is the union of intervals trivially covering $(t',t'')$, and hence it lifts to a trivial family of complex oriented curves and the quantum index persists along each of the components of the latter family.

It is convenient to
look at the elements $\widehat\bw\in \Men_\R^\rho(\Delta)$ as follows:
pick some pair $(\sigma_0,i_0)$, where $\sigma_0\in P^1_\Delta$ and $1\le i_0\le n^{\sigma_0}$,
and think of $w^{\sigma_0}_{i_0}$ as a mobile point, whose position is determined
by the other points $w^\sigma_i\in\widehat\bw$ via the Menelaus condition.
Put $C = \bn_*(\PP^1)$.
Observe that the tangent space at $C$ to the family of curves $C'\in|{\mathcal L}_{P_\Delta}|$ intersecting $\Tor(\sigma)$ at $w^\sigma_i$ (where $(\sigma,i)\ne(\sigma_0,i_0)$) with multiplicity at least $2k^\sigma_i$ can be identified with the linear system $\{C'\in|{\mathcal L}_{P_\Delta}|\ :\ (C'\cdot C)_{w^\sigma_i}\ge 2k^\sigma_i\}$. In turn, the tangent space at $C$ to the family of curves $C'\in|{\mathcal L}_{P_\Delta}|$ intersecting $\Tor(\sigma_0)$ at a point close to $w^{\sigma_0}_{i_0}$ with multiplicity at least $2k^{\sigma_0}_{i_0}$ can be identified with the linear system $\{C'\in|{\mathcal L}_{P_\Delta}|\ :\ (C'\cdot C)_{w^{\sigma_0}_{i_0}}\ge 2k^\sigma_i-1\}$.
To ensure the behavior of ${\mathcal M}_{0,n}^{\R,+}(\Delta,\widehat\bw(t))_{t'<t<t''}$
as indicated in the previous paragraph, we have to show the transversality of intersection
of the tangent space at $C$ to the family of rational curves and the tangent spaces to families of curves matching
the above tangency conditions.
The required transversality is reformulated in the cohomology language (and proved) in the following lemma.

\begin{lemma}\label{lis2}
Let $\xi=[\bn:(\PP^1,\bp)\to\Tor(P_\Delta)]\in{\mathcal M}_{0,n}^{\R,+}(\Delta,\widehat\bw)$ with $\widehat\bw\in\Men_0(\Delta)$. 
Then,
$$H^0\left(\PP^1,{\mathcal N}_\bn\left(-\sum_{\sigma\ne\sigma_0,i\ne i_0}2k^\sigma_ip^\sigma_i-(2k^{\sigma_0}_{i_0}-1)p^{\sigma_0}_{i_0}\right)\right)$$
$$=H^1\left(\PP^1,{\mathcal N}_\bn\left(-\sum_{\sigma\ne\sigma_0,i\ne i_0}2k^\sigma_ip^\sigma_i-(2k^{\sigma_0}_{i_0}-1)p^{\sigma_0}_{i_0}\right)\right)=0.$$
\end{lemma}

{\bf Proof.}
We have
$$\deg{\mathcal N}_\bn\left(-\sum_{\sigma\ne\sigma_0,i\ne i_0}2k^\sigma_ip^\sigma_i-(2k^{\sigma_0}_{i_0}-1)p^{\sigma_0}_{i_0}\right)
\overset{\text{(\ref{eis4})}}{=}-1>-2.$$
Thus, the claim of the lemma follows
from the Riemann-Roch theorem.
\proofend

\subsection{Walls in the space of constraints}\label{walls0}

Here, we characterize elements $\xi\in\overline{\mathcal M}_{0,n}^{\R,+}(\Delta,\widehat\bw(t))_{0\le t\le1}$ such that
$\widehat\bw(t)$ is a generic member of a stratum of dimension $n-1$ 
in $\Men_\R^\rho(\Delta) \setminus \Men_0(\Delta)$ ({\it i.e.,} of codimension $1$ in $\Men_\R^\rho(\Delta)$).
In view of the results of Section \ref{chamb0}, we can ignore the walls along which
$\overline{\mathcal M}^\R_{0,n}(\Delta,\widehat\bw)={\mathcal M}^\R_{0,n}(\Delta,\widehat\bw)$.
The following lemma classifies other walls.

\begin{lemma}\label{al4} (1) The following elements $\xi=[\bn:(\widehat C,\bp)\to\Tor(P_\Delta)]$ cannot occur in $\overline{\mathcal M}_{0,n}^{\R,+}(\Delta,\widehat\bw(t))_{0\le t\le1}\setminus{\mathcal M}_{0,n}^{\R,+}(\Delta,\widehat\bw(t))_{0\le t\le1}$:
\begin{enumerate}\item[(1i)] $\widehat C$ is a reducible connected curve of arithmetic genus $0$ with a component mapped onto a toric divisor;
\item[(1ii)] $\widehat C$ is a reducible connected curve of arithmetic genus $0$ with at least three irreducible components.
\end{enumerate}

(2) If $\xi=[\bn:(\widehat C,\bp)\to\Tor(P_\Delta)]\in\overline{\mathcal M}_{0,n}(\Delta,\widehat\bw(t^*))$,
where $\widehat\bw(t^*)$ is a generic element in an $(n-1)$-dimensional stratum in $\Men_\R^\rho(\Delta) \setminus \Men_0(\Delta)$,
then $\xi$ is of one of the following types:
\begin{enumerate}
\item[(2i)] $\widehat\bw(t^*)$ consists of $n+1$ distinct points, $\widehat C\simeq\PP^1$,
the map $\bn$ is birational onto its image and satisfies the following:
either it is smooth at $\widehat\bw(t^*)$, but has singular branches in $\Tor(P_\Delta)^\times$,
or $C$ is immersed in $\Tor(P_\Delta)^\times$ and at some points $w_i^\sigma\in\widehat\bw(t^*)$
it has singularity of type
$A_{2m}$, $m\ge k_i^\sigma$,
whereas $C$ is smooth at the remaining points of $\widehat\bw(t^*)$;
\item[(2ii)] two points of the sequence $\widehat\bw(t^*)$ coincide
{\rm (}$w_i^\sigma(t^*)=w_j^\sigma(t^*)$
for some $\sigma\in P^1_\Delta$
and $i\ne j${\rm )}
and $\widehat C\simeq\PP^1$, the map $\bn$ being an immersion such that the point $w_i^\sigma(t^*)=w_j^\sigma(t^*)$
is a center of one or two smooth branches;
\item[(2iii)] $\widehat\bw(t^*)$ consists of $n+1$ distinct points, $\widehat C=\widehat C_1\cup\widehat C_2$, where $\widehat C_1\simeq\widehat C_2\simeq\PP^1$ and
$\widehat C_1\cap\widehat C_2$ is one point $p$, each map $\bn:\widehat C_i\to\Tor(P_\Delta)$, $i = 1, 2$,
is either an immersion smooth along
$\Tor(\partial P_\Delta)$, or a multiple covering of a line intersecting only two toric divisors, while these divisors correspond to opposite parallel sides of $P_\Delta$ and the intersection points with these divisors are ramification points of the covering with the maximal ramification index;
furthermore, either $\bn(p)\in\Tor(P_\Delta)^\times$
and then all intersection points of the curves $C_1=\bn(\widehat C_1)$, $C_2=\bn(\widehat C_2)$ are ordinary nodes,
or $w=\bn(p)\in\bw(t^*)$, and in the latter case at least one of the maps $\bn:\widehat C_i\to\Tor(P_\Delta)$ is birational
onto its image,
the curves $C_1=\bn(\widehat C_1)$, $C_2=\bn(\widehat C_2)$ intersect $\Tor(\sigma)$ at $w$ with even multiplicities, and they do not have common point in $\widehat\bw(t^*)\setminus\{w\}$;
\item[(2iv)] $P_\Delta$ is a triangle, $n=3$, two points of the sequence $\widehat\bw(t^*)$ coincide,
the curve $C=\bn(\PP^1)$ is rational and smooth at $\widehat\bw(t^*)$, and
$\bn:\PP^1\to C$ is a double covering ramified at two distinct points of $\widehat\bw(t^*)$.
\end{enumerate}
\end{lemma}

{\bf Proof.}
(1i) A toric divisor $\Tor(\sigma)$ cannot split off alone, since, otherwise, in the deformation along the path $\{\widehat\bw(t)\}_{0\le t\le1}$,
the intersection points with the neighboring toric divisors would yield points on these toric divisors
on the distance less than $\rho$
from the corners of $\Tor_\R^+(P_\Delta)$.
For the same reason we obtain that only all toric divisors together may split off,
while their intersection points must smooth out in the deformation.
However, this contradicts the rationality of the considered curves.

\smallskip
(1ii) First, note that $\bn(\widehat C)$ intersects $\Tor(\partial P_\Delta)$ only at $\widehat\bw(t^*)$. Furthermore, if the images of
two
irreducible components $\widehat C_1,\widehat C_2$ of $\widehat C$ contain the same point $w_i^\sigma(t^*)$ that is different from any other point
$w_j^\sigma(t^*)$, $j\ne i$, then $w_i^\sigma(t^*)=\bn(\widehat C_1\cap\widehat C_2)$ (cf. Lemma \ref{l3-1a}). Due to the genus restriction,
it follows that the sequence
$\widehat\bw(t^*)$ satisfies at least three Menelaus conditions, which cuts off $\Men_\R^\rho(\Delta)$ a polytope of dimension $\le n-2$, contrary to the
dimension $n - 1$ assumption.

\smallskip
Let $\xi$ satisfy the hypotheses of item (2).

\smallskip
(2i) Let $\widehat C\simeq\PP^1$, $\widehat\bw(t^*)$ consist of $n+1$ distinct points, and let $\bn$ take $\widehat C$ birationally onto $C=\bn(\widehat C)\subset\Tor(P_\Delta)$.

First, we consider the minimal value $n=2$. In such a case, $P_\Delta$ is a triangle and, for each $\widehat\bw\in \Men_\R^\rho(\Delta)$, there exists a unique curve $C=\bn(\PP^1)$ (where $\xi=[\bn:(\PP^1,\bp)\to\Tor(P_\Delta)]\in{\mathcal M}_{0,n}^{\R,+}(\Delta,\widehat\bw)$ which is smooth along the toric divisors and is nodal, while the number of elliptic nodes in $\Tor^+_\R(P_\Delta)$ does not depend on the choice of $\widehat\bw\in\Men_\R^\rho(\Delta)$, see \cite[Lemma 3.5]{Sh0}.

Suppose that $n\ge3$.
Note that $\xi$ must be unibranch at each point of $\widehat\bw(t^*)$.
Since $\xi$ is a generic element of an $(n-1)$-dimensional family in ${\mathcal M}_{0,n}(\Delta)$, fixing the position of a point $w_i^\sigma\in\widehat\bw(t^*)$, we obtain a family of dimension $\ge n-2\ge1$. Applying \cite[Inequality (5) in Lemma 2.1]{IKS4},
we get
$$-CK_{\Tor(P_\Delta)}\ge2+\left(-CK_{\Tor(P_\Delta)}
-n\right)+(n-3)+\sum_B(\mt Q-1)$$
\begin{equation}=-CK_{\Tor(P_\Delta)}-1+\sum_B(\mt Q-1)\quad\Longrightarrow\quad\sum_B(\mt Q-1)\le1\ ,\label{ae9f}\end{equation}
where $B$ is the set of all singular branches of $\xi$ in $\Tor(P_\Delta)^\times\cup\{w_i^\sigma\}$.
In particular, if $C$ is singular at some point $w_i^\sigma$, the singularity is of type
$A_{2m}$.
The requirements of item (2i) are fulfilled.

\smallskip
(2ii) Suppose that $\widehat C\simeq\PP^1$ is birationally taken by $\bn$ onto its image in $\Tor(P_\Delta)$,
and some of the points of the sequence $\widehat\bw(t^*)$ coincide. The latter is possible only if $n\ge3$.
For the dimension reason, the number $\#(\widehat\bw(t^*))$
of points in $\widehat\bw(t^*)$ is equal to $n$
({\it i.e.}, $w_i^\sigma(t^*)=w_j^\sigma(t^*)=w$
for some $\sigma\subset\partial P_\Delta$ and $i\ne j$), and all these $n\ge3$ points are in general position on $\Tor(\partial P_\Delta)$ subject
to the Menelaus relation (\ref{ae1}). If $C=\bn(\widehat C)$ is unibranch at each point of $\widehat\bw(t^*)$, then, by Lemma \ref{l3-1a},
the curve $C$ is immersed and smooth along the toric divisors. If $C$ is not unibranch at $w$,
then by Lemma \ref{l3-1a} it has two local branches at $w$; furthermore,
$C$ must be unibranch at each point $w_{i'}^{\sigma'}(t^*)\ne w$. So, assume that $C$ has two local branches at $w$.

If $n-1=2$, then $P_\Delta$ is a lattice triangle.
Using affine automorphisms of $\Z^2$,
we can take $P_\Delta=\conv\{(m_1,0),(m_1+k_1+k_2,0),(0,m_2)\}$
with $m_1,m_2,k_1,k_2$ even and write a parameterization of $C$ in the form
$$x=a\tau^{m_2},\quad y=b\tau^{m_1}(\tau-1)^{k_1}(\tau-\theta)^{k_2},\quad a,b,\theta\in\R^*.$$
Since $x(1)=x(\theta)$, we obtain $\theta=-1$. Note that $k_1\ne k_2$, since otherwise, $x$ and $y$ would be functions of $\tau^2$, that is, $\bn:\widehat C\to C$ would be a double covering against the initial assumption. Relations
$$x=a(1+m_2(\tau-1)+O((\tau-1)^2),\quad x=a(1-m_2(\tau+1)+O((\tau+1)^2),$$
show that the local branches of $C$ at $w$ are smooth. Furthermore,
$$x=a\tau^{m_2},\quad y=b\tau^{m_1}+b(k_2-k_1)\tau^{m_1+1}+O(\tau^{m_1+2});$$
hence, $C$ is smooth at the two remaining points of $\widehat\bw$. In addition,
$C$ is
immersed
everywhere, since $\frac{dx}{d\tau}\ne0$ for $\tau\ne0$.

Let $n-1\ge3$.
We
show that $C$ is immersed.
Fixing the position of $w$ and the position of one more point
$w'\in\bw(t^*)\setminus\{w\}$, we obtain a family of dimension at least $n-3\ge1$;
thus,
\cite[Inequality (5) in Lemma 2.1]{IKS4} applies:
$$-CK_{\Tor(P_\Delta)}\ge2+(-CK_{\Tor(P_\Delta)}-n+2)+\sum_B(\mt Q-1)+(n-4)$$
$$=-CK_{\Tor(P_\Delta)}+\sum_B(\mt Q-1)\ ,$$
where $B$ is the set of all singular local branches of $C$ in $\Tor(P_\Delta)^\times\cup\{w,w'\}$, and
hence $C$ is immersed.

\smallskip
(2iii)
Consider the case of $\widehat C=\widehat C_1\cup\widehat C_2$,
where $\widehat C_1\simeq\widehat C_2\simeq\PP^1$ and
$\widehat C_1\cap\widehat C_2=\{p\}$ is one point.
For the dimension reason, the points of
$\widehat\bw(t^*)$ are in general position subject to exactly two Menelaus conditions (induced by the components of $\widehat C$),
and they
are all distinct.
Each of the curves $C_1$, $C_2$, which passes through at least three points of $\widehat\bw(t^*)$,
is immersed and smooth along the toric divisors
and has even intersection multiplicity
with $\Tor(\partial P_\Delta)$ at any point, by Lemma \ref{l3-1a}. If $C_1$ or $C_2$ passes through exactly two points of $\widehat\bw(t^*)$, then it is a multiple covering of a line as described in
item (2iv). Let $\bn(p)\in\Tor(P_\Delta)^\times$. Then, by Lemma \ref{l3-1a}, the curves
$C_1$ and $C_2$ do not share points in $\widehat\bw(t^*)$, and we claim that all their intersection points are ordinary nodes. Due to the genericity assumptions for $\widehat\bw(t^*)$, we have to study the only case of both $C_1$ and $C_2$ immersed. If we freely move the points of $C_1\cap\widehat\bw(t^*)$ so that the corresponding Menelaus condition induced by $C_1$ retains, and fix the curve $C_2$,
then the persisting tangency condition of a germ $\bn:(\widehat C_1,q)\to\Tor(P_\Delta)$ to the curve $C_2$ would yield that the tangent space to the considered family of curves in the linear system $|C_1|$ would be contained in $H^0(\widehat C_1,{\mathcal O}_{\widehat C_1}(\bd_1))$, where
$$\deg\bd_1=C_1^2-(C_1^2+C_1K_{\Tor(P_\Delta)}
+2)-(-C_1K_{\Tor(P_\Delta)}
-\#(\bw(t^*)\cap C_1))-1$$ $$=
\#(\widehat\bw(t^*)\cap C_1)-3>-2\ ,$$ and hence, by the Riemann-Roch theorem,
$$h^0(\widehat C_1,{\mathcal O}_{\widehat C_1}(\bd_1))=\#(\widehat\bw(t^*)\cap C_1)-3+1=\#(\widehat\bw(t^*)\cap C_1)- 2
< \#\widehat\bw(t^*)\cap C_1)-1\ ,$$
which is a contradiction. If $\bn(p)=w_i^\sigma(t^*)\in\widehat\bw(t^*)$, then the intersection multiplicities of the germs $\bn(\widehat C_1,p)$ and $\bn(\widehat C_2,p)$ with $\Tor(\sigma)$ are even. Indeed, otherwise, $C_1$ and $C_2$ would have at least one additional local branch centered on some $\Tor(\sigma')$, $\sigma'\ne\sigma$, and intersecting $\Tor(\sigma')$ with an odd multiplicity, which contradicts our assumptions.

\smallskip
(2iv) Suppose that $\bn:\PP^1\to\Tor(P_\Delta)$ is an $s$-multiple covering onto its image, $s\ge2$.
For a generic $t\in[0,1]$, $\widehat\bw(t)$ consists of $n+1$ distinct points; hence the sequence $\widehat\bw(t^*)$ contains at least $n$ distinct points.
If $\widehat\bw(t^*)$ consists of $n+1$ distinct points, then
at each point of $\widehat\bw(t^*)$ we have the ramification index $s$;
thus, by the Riemann-Hurwitz formula, we get
\begin{equation}2\le2s-(s-1)(n+1),
\label{ae97}\end{equation}
which
gives $n\le 1$; the latter inequality
contradicts the fact that $n+1$ is
bounded from below by the number of sides of $P_\Delta$.
If the sequence $\widehat\bw(t^*)$ contains only $n$ distinct points,
we have $n\ge3$. Note that $\bn$ has an irreducible preimage
in at least $n-1$ points of $\widehat\bw(t^*)$ with ramification index $s$, and hence $n=3$ (cf. (\ref{ae97})).
It also follows that there are no other ramifications and that $s=2$, since the remaining point of $\widehat\bw(t^*)$,
where $\bn$ is not ramified,
lifts to at most two points in $\PP^1$.
Thus, we are left with the case
described in item (2iv) of the lemma.
\proofend

\subsection{Wall-crossing}\label{crossing0}
We complete the proof of Theorem \ref{at1} with the following lemma.

\begin{lemma}\label{al8}
Let $\{\widehat\bw(t)\}_{0\le t\le1}$ be a generic path in $\Men_\R^\rho(\Delta)$, and let $t^*\in(0,1)$ be such that $\overline{\mathcal M}_{0,n}^\R(\Delta,\widehat\bw(t^*))$ contains an element
$\xi$ as described in one of the items of Lemma \ref{al4}(2).
Then, for each $\kappa \in \frac{1}{2}\Z$ such that $|\kappa| \leq \cA(\Delta)$,
the numbers $W_0^\kappa(t)=W_0^\kappa(\Delta,\widehat\bw(t))$ and $\widetilde W_0^\kappa(t)=\widetilde W_0^\kappa(\Delta,\widehat\bw(t))$
do not change as $t$ varies in a neighborhood of $t^*$.
\end{lemma}

{\bf Proof.}
We always can assume that in a neighborhood of $t^*$, the path $\{\widehat\bw(t)\}_{0\le t\le1}$ is defined by fixing the position of some $n-1$ points of $\widehat\bw(t^*)$, while the other two points remain mobile (the choice of the two mobile points may depend on the considered degeneration). We also
notice that, in the degenerations as in Lemma \ref{al4}(2i,2ii),
the source curve and its real structure remain fixed,
which implies that the quantum index is constant in these wall-crossings. Except for the case of Lemma \ref{al4}(2iii) describing reducible degenerations, we work with families of curves which are trivially covered by families of complex oriented curves so that the quantum index persists along each component of the family of oriented curves.

The proof is separated into parts (1)-(6) in accordance with the classification of walls
in Lemma \ref{al4}(2)
and specific degenerations occurring on these walls.

\smallskip

{\bf(1)} Suppose that $\xi\in\overline{\mathcal M}^{\R,+}_{0,n}(\Delta,\widehat\bw(t^*))$ is as in Lemma \ref{al4}(2i), and $\xi$ has a singular local branch centered in $\Tor(P_\Delta)^\times$. Then the constancy of
$W_0^\kappa(t)$, $|t-t^*|<\eps$, follows from \cite[Lemmas 13 and 15]{IKS3}, provided, we establish the following transversality condition. Choose a sufficiently large integer $s$. For each point $z\in\Sing(C)\cap\Tor(P_\Delta)^\times$, we set $I_z=I^{cond}(C,z)/{\mathfrak m}_{C,z}^s\subset{\mathcal O}_{C,z}
/{\mathfrak m}_z^s$, the quotient of the conductor ideal by the power of the maximal ideal, which can be viewed as the tangent cone to the stratum parameterizing equigeneric deformations
(see \cite[Theorem 4.15]{DH}). The condition we need reads
$$h^1\left(C,{\mathcal J}_{Z/C}\otimes_{{\mathcal O}_{\widehat C}}\ina^*{\mathcal L}_{P_\Delta}
\right)=0,$$
where ${\mathcal J}_{Z/C}$ is the ideal sheaf of the zero-dimensional scheme $Z\subset C$ defined at each point $z\in\Sing(C)$ by the conductor ideal $I^{cond}(C,z)$ and at each point $w=w_i^\sigma\in\widehat\bw(t^*)$ by the ideal
$$I_w=\{\varphi\in{\mathcal O}_{C,w}\ :\ (\varphi\cdot C)_w\ge 2k_i^\sigma-1\}\subset{\mathcal O}_{C,w},$$
The above $h^1$-vanishing lifts to $\PP^1$ in the form
\begin{equation}h^1(\PP^1,{\mathcal O}_{\PP^1}(\bd))=0,\label{ae9a}\end{equation}
where
$$\deg\bd=C^2-(C^2+CK_{\Tor(P_\Delta)}+2)-(CK_{\Tor(P_\Delta)}-n-1)$$
\begin{equation}=n-1>-2=2g-2,\label{ae9cc}\end{equation}
and hence the required $h^1$-vanishing comes from the Riemann-Roch theorem.

\smallskip
{\bf(2)} Suppose that $\xi\in\overline{\mathcal M}^{\R,+}_{0,n}(\Delta,\widehat\bw(t^*))$ is as in Lemma \ref{al4}(2i)
with singularities of type
$A_{2m}$ at some points $w_i^\sigma\in\widehat\bw(t^*)$, $m\ge k_i^\sigma$.

The germ at $\xi$ of ${\mathcal M}_{0,n}(\Delta)$ (defined by (\ref{efamily})) can be identified with the germ at $0$ of $H^0(C,\ina^*{\mathcal L}_{P_\Delta})$, which in turn can be naturally embedded into $\prod_{\sigma,i}\Vers(C,w_i^\sigma)\times\prod_{z\in\Sing(C)\cap\Tor(P_\Delta)^\times}\Vers(C,z)$. We claim that the image of the germ at $0$ of $H^0(C,\ina^*{\mathcal L}_{P_\Delta})$ in $\prod_{i,\sigma}\Vers(C,w_i^\sigma)\times\prod_{z\in\Sing(C)\cap\Tor(P_\Delta)^\times}\Vers(C,z)$ intersects transversally with $\prod_{\sigma,i}EG_{\Tor(\sigma)}(C,w_i^\sigma)\times\prod_{z\in\Sing(C)\cap\Tor(P_\Delta)^\times}EG(C,z)$. Hence, this intersection ${\mathcal V}$ is smooth of the expected dimension $n$ since the latter product is smooth (see Section \ref{sec-def}, items (c) and (d)).
The required transversality amounts to the validation of the relation
\begin{equation}h^1(C,{\mathcal J}_{Z_1/C}\otimes_{{\mathcal O}_C}\ina^*{\mathcal L}_{P_\Delta})=0,
\label{ae9b}\end{equation}
where the zero-dimensional scheme $Z_1\subset C$ supported at $\widehat\bw(t^*)\cup\Sing(C)$ and defined by the ideals $I^{cond}(C,z)$ at all $z\in\Sing(C)\cap\Tor(P_\Delta)^\times$ and by the ideals
\begin{equation}\{\varphi\in{\mathcal O}_{C,w_i^\sigma}\ :\ \ord\varphi\big|_{C,w_i^\sigma}\ge2k_i^\sigma-1+2\delta(C,w_i^\sigma)\}\quad\text{for all}\ w_i^\sigma\in\widehat\bw(t^*)\label{ae9d}\end{equation}
(cf. Section \ref{sec-def}, items (c) and (d)). Relation (\ref{ae9b}) is equivalent to
$$h^1(\PP^1,{\mathcal O}_{\PP^1}(\bd_1))=0,$$
where $\deg \bd_1$ is given by (\ref{ae9cc}), and hence it holds by the Riemann-Roch theorem.

Note that $\xi$ is chosen as a generic member of the wall, that is, it represents a smooth point of the wall. Let us show that the germ ${\mathcal V}_1$ of the wall at $\xi$ has codimension one in ${\mathcal V}$. Pick a point $w_i^\sigma\in\widehat\bw(t^*)$, where $\bn(\widehat C)$ has a singularity of type $A_{2m}$, $m\ge k_i^\sigma$, and consider the intersection ${\mathcal V}_2$ of the image of the germ at $0$ of $H^0(C,\ina^*{\mathcal L}_{P_\Delta})$ in $\prod_{\sigma',i'}\Vers(C,w_{i'}^{\sigma'})\times\prod_{z\in\Sing(C)\cap\Tor(P_\Delta)^\times}\Vers(C,z)$ with \begin{equation}\prod_{(\sigma',i')\ne(\sigma,i)}ES_{\Tor(\sigma')}(C,w_{i'}^{\sigma'})\times ES^{fix}_{\Tor(\sigma)}(C,w_i^\sigma)
\times\prod_{z\in\Sing(C)\cap\Tor(P_\Delta)^\times}EG(C,z).\label{elocus}\end{equation} It has dimension $\dim{\mathcal V}_2\ge\dim{\mathcal V}_1-1$. We claim that $\dim{\mathcal V}_2\le\dim{\mathcal V}-2=n-2$. Indeed,
\begin{equation}h^1(C,{\mathcal J}_{Z_2/C}\otimes_{{\mathcal O}_C}\ina^*{\mathcal L}_{P_\Delta})=0,\label{ae9bb}\end{equation}
where the zero-dimensional scheme $Z_2$ supported at $\widehat\bw(t^*)\cup\Sing(C)$ coincides with $Z_1$ (appearing in (\ref{ae9b})) along $(\widehat\bw(t^*)\setminus\{w_i^\sigma\})\cup\Sing(C)$ and is defined at $w_i^\sigma$ by the ideal
$$\{\varphi\in{\mathcal O}_{C,w_i^\sigma}\ :\ \ord\varphi\big|_{C,w_i^\sigma}\ge2k_i^\sigma+1+2\delta(C,w_i^\sigma)\}$$
(cf. Section \ref{sec-def}, item (b)). To see this, as in the preceding paragraph we rewrite (\ref{ae9bb}) in the form
$$h^1(\PP^1,{\mathcal O}_{\PP^1}(\bd_2))=0,$$
where
$$\deg\bd_2=C^2-(C^2+CK_{\Tor(P_\Delta)}+2)-(CK_{\Tor(P_\Delta)}-n)-1=n-3>-2.$$
Thus, by the Riemann-Roch theorem, we have
$$\dim{\mathcal V}_2\le h^0(C,{\mathcal J}_{Z_2/C}\otimes_{{\mathcal O}_C}\ina^*{\mathcal L}_{P_\Delta})=(n-3)+1=n-2,$$
as required. Furthermore, relation (\ref{ae9bb}) yields the existence of the germ of a smooth path through $0$ in ${\mathcal V}$ transversally intersecting ${\mathcal V}_1$ and regularly projecting onto a smooth path through $\widehat\bw(t^*)$ in $\Men(\Delta)$.

\begin{figure}
\setlength{\unitlength}{1cm}
\begin{picture}(14,7)(-1.5,-1)
\includegraphics[width=0.8\textwidth, angle=0]{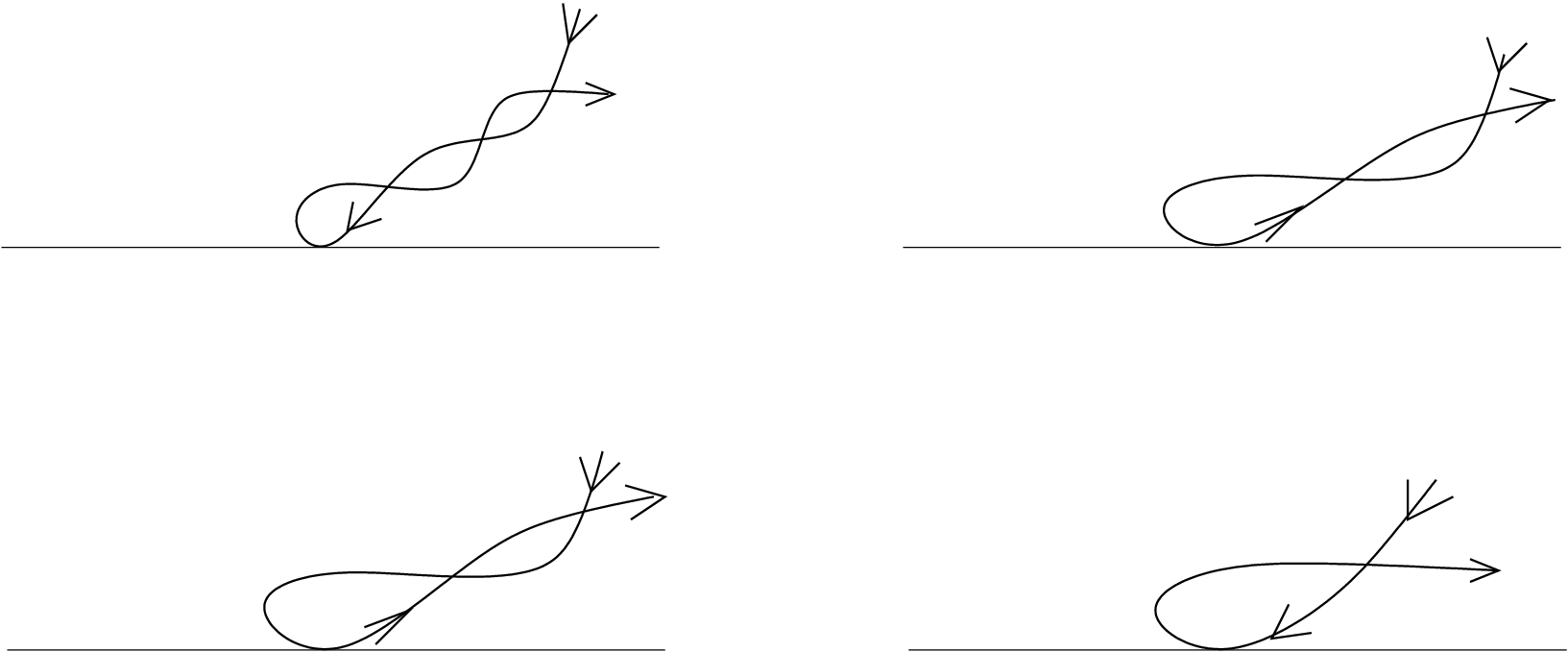}
\put(-9.7,2.7){$w_i^\sigma$}\put(-2.8,2.7){$w_i^\sigma$}\put(-9.7,-0.4){$w_i^\sigma$}\put(-2.8,-0.4){$w_i^\sigma$}
\put(-10.5,2.2){(a)}\put(-3.7,2.2){(b)}\put(-10.5,-1){(c)}\put(-3.7,-1){(d)}
\put(-11,3.4){$\bullet$}\put(-3.7,3.3){$\bullet$}\put(-4.2,3.4){$\bullet$}
\put(-11,-0.4){$\bullet$}\put(-3.7,-0.4){$\bullet$}\put(-4.6,-0.5){$\bullet$}
\put(-10,4.9){$t<t^*$}\put(-3,4.9){$t>t^*$}
\put(-10,1.3){$t<t^*$}\put(-3,1.3){$t>t^*$}
\end{picture}
\caption{Proof of Lemma \ref{al8}, part (2)}
\label{af6}
\end{figure}

Having a fragment of the path $\{\widehat\bw(t)\}_{|t-t^*|<\eps}$ whose existence we verified above, it remains to show that the value of $W^\kappa_0(\Delta,\widehat\bw(t))$ is the same for $t<t^*$ and for $t>t^*$. Given a point $w_i^\sigma\in\widehat\bw(t^*)$, where $C$ has a singularity of type $A_{2m}$, $m\ge k_i^\sigma$, take local coordinates $x,y$ in a neighborhood $U\subset\Tor(P_\Delta)$ of $w_i^\sigma$ so that $\Tor(\sigma)=\{y=0\}$.
Suppose that $C\cap U=\{f(x,y)=y^2-2yx^k+x^{2k}+\text{h.o.t.}=0\}$, where $k=k_i^\sigma$
and h.o.t. denotes terms above the Newton diagram.
In particular, the real part $\R C\cap U$ lies in $\{y\ge0\}$.
Without loss of generality, we can assume that $C$ deforms along the path $\{\widehat\bw(t)\}_{|t-t^*|<\eps}$ into nodal curves $C^{(t)}\cap U$, $0<|t-t^*|<\eps$.
We
show that all
hyperbolic nodes live in the domain $\{x>0,y>0\}$,
and the elliptic nodes live in the domain $\{y>0\}$
for an even $k$ and in the domain $\{y<0\}$ for an odd $k$.
Then, the constancy of $W_0^\kappa(\Delta,\widehat\bw(t))$, $0<|t-t^*|<\eps$,
follows from the fact that the total number of real nodes of $C^{(t)}\cap U$ equals $m$ modulo $2$
and that the orientation of $C^{(t)}\cap U$ at $w_i^\sigma$ correlates
with the number of the hyperbolic nodes,
see Figure \ref{af6}(a,b) for an even $k$ and Figure \ref{af6}(c,d) for an odd $k$.

Since the deformation $\{C^{(t)}\}_{|t-t^*|<\eps}$ is equigeneric,
we can describe it {\it via} deformation of local parameterizations,
see \cite[Theorem 1 on page 73]{Tei} and \cite[Theorems 4.1 and 4.2]{ChL}. Thus, $C\cap U=C^{(t^*)}\cap U$
admits a (real) parametrization
$$x=\tau^2+a_{2m+1}\tau^{2m+1}+\sum_{2\le r\le m}a_{2r}\tau^{2r}+\sum_{r\ge2m-2k+2}a_r\tau^r,\quad y=\tau^{2k},$$ where $\tau\in(\C,0)$ and $a_{2m-2k+1}\ne0$,
and the curves $C^{(t)}\cap U$ are parameterized by
$$x=\tau^2+a_{2m-2k+1}\tau^{2m-2k+1}+\sum_{2\le r\le m}a_{2r}\tau^{2r}+\sum_{r\ge2m-2k+2}a_r\tau^r+\sum_{r\ge1}b_r\tau^r,\quad y=\tau^{2k},$$ 
where $b_r=b_r(t)$ are analytic functions defined over $\R$ and vanishing at $t=t^*$.

The curves $C^{(t)}$, $t\ne t^*$, satisfy the conditions of Lemma \ref{l3-1a}, that is, $\alpha_1(t)>0$ for $0<|t-t^*|<\eps$. Since $y=\tau^{2k}>0$ as $\tau\in\R^*$, all hyperbolic nodes are located in the domain $\{y>0\}$.
Furthermore, a hyperbolic node corresponds to values $\tau_1,\tau_2\in\R$, $\tau_2=-\tau_1$, such that $x(\tau_)=x(\tau_2)$,
which yields $x(\tau_1)=\tau_1^2+O(\tau_1^4)>0$. Consider an elliptic node of $C^{(t)}$, $t\ne t^*$.
It corresponds to $\tau\in(\C,0)\setminus(\R,0)$ such tha $x(\tau),y(\tau)\in \R$.
Let $k$ be odd. Suppose that $y(\tau)>0=\tau^{2k}>0$, {\it i.e.},
$\tau=\theta\zeta_{2k}$, where $\theta>0$ and $(\zeta_{2k})^{2k}=1$, $\zeta_{2k}\ne\pm1$. Then,
$$x(\tau)=\alpha_1\theta\zeta_{2k}+(1+\alpha_2)\theta^2(\zeta_{2k})^2+O(\theta^3)\not\in\R,$$
since (recall that $k$ is odd)
$$(1+\alpha_2)\theta^2(\zeta_{2k})^2\not\in\R\quad\text{and}\quad\frac{\alpha_1\theta\zeta_{2k}}{(1+\alpha_2)
\theta^2(\zeta_{2k})^2}
\not\in\R.$$ Hence,
all elliptic nodes live in the domain $\{y<0\}$. Let $k$ be even. Suppose that $y(\tau)=\tau^{2k}<0$, {\it i.e.},
$\tau=\theta\eta\zeta_{2k}$, where $\theta>0$, $(\zeta_{2k})^{2k}=1$, $\eta=\exp\left(\frac{\pi\sqrt{-1}}{2k}\right)$. Thus,
$$x(\tau)=\alpha_1\theta\eta\zeta_{2k}+(1+\alpha_2)\theta^2\eta^2(\zeta_{2k})^2+O(\theta^3),$$
and again we have
$$(1+\alpha_2)\theta^2(\eta^2\zeta_{2k})^2\not\in\R\quad\text{and}\quad\frac{\alpha_1\theta\zeta_{2k}}{(1+\alpha_2)
\theta^2\eta^2(\zeta_{2k})^2}
\not\in\R.$$
Hence, all elliptic nodes live in the domain $\{y>0\}$.

\smallskip
{\bf(3)} Suppose that $\xi$ is as in Lemma \ref{al4}(2ii). Here we choose a path $\{\widehat\bw(t)\}_{0\le t\le1}$ defined in a neighborhood of $t^*$ so that the point $w_i^\sigma(t)$ is fixed and
the point $w_j(t^*)$ is mobile (together with some other point of $\widehat\bw(t)$). Since $C$ is immersed, we extract the required local constancy of $W_0^\kappa(t)$ and $\widetilde W_0^\kappa(t)$ from the transversality relation in the form of (\ref{ae9a}) which holds due to
\begin{equation}\deg\bd=C^2-(C^2+CK_{\Tor(P_\Delta)}+2)-(-CK_{\Tor(P_\Delta)}-2)=0>-2,\label{etrans1}\end{equation}
and from
the smoothness of the two following strata
in ${\mathcal O}_{C,w}/{\mathfrak m}_w^s$:
\begin{itemize}\item one stratum parameterizes deformations of a smooth branch intersecting
$\Tor(\sigma)$ at $w$ with multiplicity $2(k_i^\sigma+k_j^\sigma)$ into a smooth branch intersecting $\Tor(\sigma)$ at $w$ with multiplicity $2k_i^\sigma$ and in a nearby point with multiplicity $2k_j^\sigma$,
\item the other stratum parameterizes deformations of a couple of smooth branches intersecting $\Tor(\sigma)$ at $w$ with multiplicities $2k_i^\sigma$ and $2k_j^\sigma$, respectively, into a couple of smooth branches, one intersecting $\Tor(\sigma)$ at $w$ with multiplicity $2k_i^\sigma$ and the other intersecting $\Tor(\sigma)$ in a nearby point with multiplicity $2k_j^\sigma$.
\end{itemize}
Both smoothness statements are straightforward.

\begin{figure}
\setlength{\unitlength}{1cm}
\begin{picture}(12,7)(-3,-1)
\includegraphics[width=0.6\textwidth, angle=0]{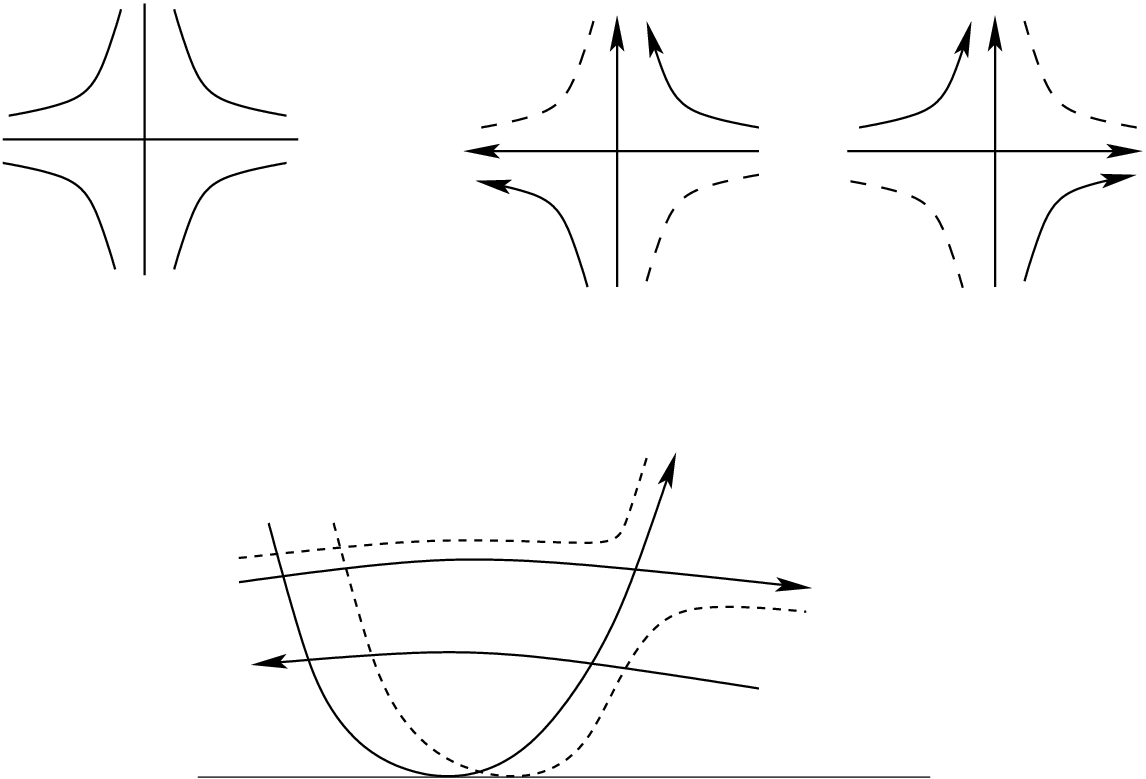}
\put(-8.15,3.2){(a)}\put(-3,3.2){(b)}\put(-4.5,-0.8){(c)}
\put(-9.5,4.2){$\lambda>0$}\put(-9.5,5.6){$\lambda<0$}
\put(-7.2,4.2){$\lambda<0$}\put(-7.2,5.6){$\lambda>0$}
\put(-6.2,-0.4){$w_i^\sigma(t^*)$}\put(-2,0.2){$\Tor(\sigma)$}
\put(-3.6,2.4){$C_1$}\put(-2.9,1.7){$C_2$}\put(-3.4,0.3){$C_2$}\put(-2.7,1){$C^{(t)}$}
\end{picture}
\caption{Proof of Lemma \ref{al8}, part (4)} \label{af1}
\end{figure}

\smallskip
{\bf(4)} Suppose that $\xi$ is as in Lemma \ref{al4}(2iii) and $\bn(p)\in\Tor_\R(P_\Delta)^\times_+$. We closely follow the proof of \cite[Theorem 5]{Mir}, providing here details for the reader's convenience.
Consider the path $\{\widehat\bw(t)\}_{0\le t\le1}$ locally obtained from $\widehat\bw(t^*)$ by moving one point of $\widehat\bw(t^*)\cap C_1$ and one point of $\widehat\bw(t^*)\cap C_2$, while fixing the position of the other points of $\widehat\bw(t^*)$. Then,
we obtain a smooth one-dimensional germ in $\overline{\mathcal M}_{0,n}(\Delta)$, which locally induces a deformation of the hyperbolic node $\bn(\widehat C,p)$ equivalent in suitable local coordinates in a neighborhood of $z=\bn(p)$ to
$uv=\lambda(t)$ with $\lambda$ a smooth function in a neighborhood of $t^*$ such that $\lambda(t^*)=0$ and $\lambda'(t^*)>0$ (see Figure \ref{af1}(a)).
The smoothness and the local deformation claims are straightforward from the (standard) considerations in the proof of \cite[Lemma 11(2)]{IKS3}, which reduce both claims to the following $h^0$-vanishing (cf. \cite[Formula (16) and computations in the first paragraph in page 251]{IKS3}):
\begin{equation}h^0(\widehat C_r,{\mathcal O}_{\widehat C_r}(\bd_r))=0,\quad r=1,2\ ,\label{ae11}\end{equation} where
$$\deg\bd_r=C_r^2-(C_r^2+C_rK_{\Tor(P_\Delta)}+2)-(C_rK_{\Tor(P_\Delta)}-1)=-1>-2\ ,$$ and hence (\ref{ae11}) follows.

Now we analyze the change of the refined invariants, closely following the argument from \cite[Section 6.4]{Mir}. To this end, we pass to the consideration of complex oriented curves. Denote by $\widehat C_1^\pm$ and $\widehat C_2^\pm$ the connected components of $\widehat C_1\setminus\R\widehat C_1$ and $\widehat C_2\setminus\R\widehat C_2$, respectively.
When $t$ varies around $t^*$, the curve $\widehat C_1\cup\widehat C_2$ turns into $\PP^1$, and we encounter two types of deformations: either $\widehat C_1^+,\widehat C_1^-$ glue up with
$\widehat C_2^+,\widehat C_2^-$, respectively, or $\widehat C_1^+,\widehat C_1^-$ glue up with $\widehat C_2^-,\widehat C_2^+$, respectively. The type of the deformation agrees with the local complex orientation of the central curve as shown
in Figure \ref{af1}(b).
This means
that if for $t<t^*$ one encounters a deformation of the first type,
then for $t>t^*$ it is of the second type, and {\it vice versa}.
Put (cf. (\ref{ae13}) and (\ref{ae9i}))
$$\eps_j^\pm=\prod_{\renewcommand{\arraystretch}{0.6}
\begin{array}{c}
\scriptstyle{\sigma\in P^1_\Delta,\ 1\le i\le n^\sigma}\\
\scriptstyle{w_i^\sigma\in C_j,\ k_i^\sigma \ \equiv \ 0  \mod 2}\end{array}}\sign(\widehat C_j^\pm,w_i^\sigma),
\quad\widetilde\eps_j^{\;\pm}=\prod_{\sigma\in P^1_\Delta,\ 1\le i\le n^\sigma}\sign(\widehat C_j^\pm,w_i^\sigma)\quad j=1,2.$$
Depending on the point of $C_1\cap C_2$ chosen for the smoothing, the triple $(\kappa,\eps,\widetilde\eps)$, where $\kappa$ is the quantum index and $\eps,\widetilde\eps$ denote the second product in (\ref{ae13}), (\ref{ae9i}), respectively, changes as follows (see Figure \ref{af1}(c)): either the real tangent vectors to $C_1,C_2$ at the smoothed out intersection points form a negative pair and
$$\begin{cases}(\kappa_1+\kappa_2,\eps_1^+\eps_2^+,\widetilde\eps_1^{\;+}\widetilde\eps_2^{\;+}),&\\ (-\kappa_1-\kappa_2,\eps_1^-\eps_2^-,\widetilde\eps_1^{\;-}\widetilde\eps_2^{\;-}),&\end{cases}\ \Longrightarrow\ \begin{cases}
(\kappa_1-\kappa_2,\eps_1^+\eps_2^-,\widetilde\eps_1^{\;+}\widetilde\eps_2^{\;-}),&\\ (-\kappa_1+\kappa_2,\eps_1^-\eps_2^+,\widetilde\eps_1^{\;-}\widetilde\eps_2^{\;+}),&\end{cases}$$ or the real tangent vectors to $C_1,C_2$ form a positive pair and
$$\begin{cases}(\kappa_1-\kappa_2,\eps_1^+\eps_2^-,\widetilde\eps_1^{\;+}\widetilde\eps_2^{\;-}),&\\ (-\kappa_1+\kappa_2,\eps_1^-\eps_2^+,\widetilde\eps_1^{\;-}\widetilde\eps_2^{\;+}),&\end{cases}\ \Longrightarrow\ \begin{cases}
(\kappa_1+\kappa_2,\eps_1^+\eps_2^+,\widetilde\eps_1^{\;+}\widetilde\eps_2^{\;+}),&\\ (-\kappa_1-\kappa_2,\eps_1^-\eps_2^-,\widetilde\eps_1^{\;-}\widetilde\eps_2^{\;-}).&\end{cases}$$
For a simple geometric reason, the number of intersection points of the first kind equals the number of the intersection points of the second type.
Hence, summing up over smoothings fo all real intersection points of $C_1,C_2$, we derive the constancy of the numbers
$W_0^\kappa(t)$ and $\widetilde W_0^\kappa(t)$ along the path $\{\widehat\bw(t)\}_{|t-t^*|<\zeta}$ ($0<\zeta\ll1$) in case of $\bn:\widehat C\to\Tor(P_\Delta)$ a birational immersion.

The case of multiple coverings can be treated in the same manner. We leave details to the reader.

\smallskip
{\bf(5)} Suppose that $\xi$ is as in Lemma \ref{al4}(2iii) and $\bn(p)=w_i^\sigma(t^*)=:w$. Observe that
$(\bn\big|_{\widehat C_r})^*(w)=2l_rp$, $r=1,2$, where $l_1+l_2=k_i^\sigma$.

\smallskip({\it5a}) Suppose that $\bn\big|_{\widehat C_1},\bn\big|_{\widehat C_2}$ both are immersions and that $l_1\le l_2$.

Define the path $\{\widehat\bw(t)\}_{|t-t^*|<\zeta}$, $0<\zeta\ll1$, by picking one mobile point in $C_1\cap\widehat\bw(t^*)\setminus\{w\}$ and the other in $C_2\cap\widehat\bw(t^*)\setminus\{w\}$. It defines a one-dimensional subvariety of $\overline{\mathcal M}_{0,n}(\Delta)$, and
let ${\mathcal V}=\{[\bn_t:(\widehat C^{(t)},\bp_t)\to\Tor(P_\Delta)]\}_{|t-t^*|<\zeta}$ be one of the irreducible germs of that variety at $\xi$, where
$\widehat C^{(t)}\simeq\PP^1$ as $t\ne0$. The one-dimensional family of curves $C^{(t)}=\bn_t(\widehat C^{(t)})\in|{\mathcal L}_{P_\Delta}|$, $|t-t^*|<\zeta$, has a tangent cone at $C_1\cup C_2$ spanned by $C_1\cup C_2$ and some curve $C^*\in|{\mathcal L}_{P_\Delta}|\setminus\{C_1\cup C_2\}$.

\begin{figure}
\setlength{\unitlength}{1cm}
\begin{picture}(12,10)(-3,-0.5)
\includegraphics[width=0.6\textwidth, angle=0]{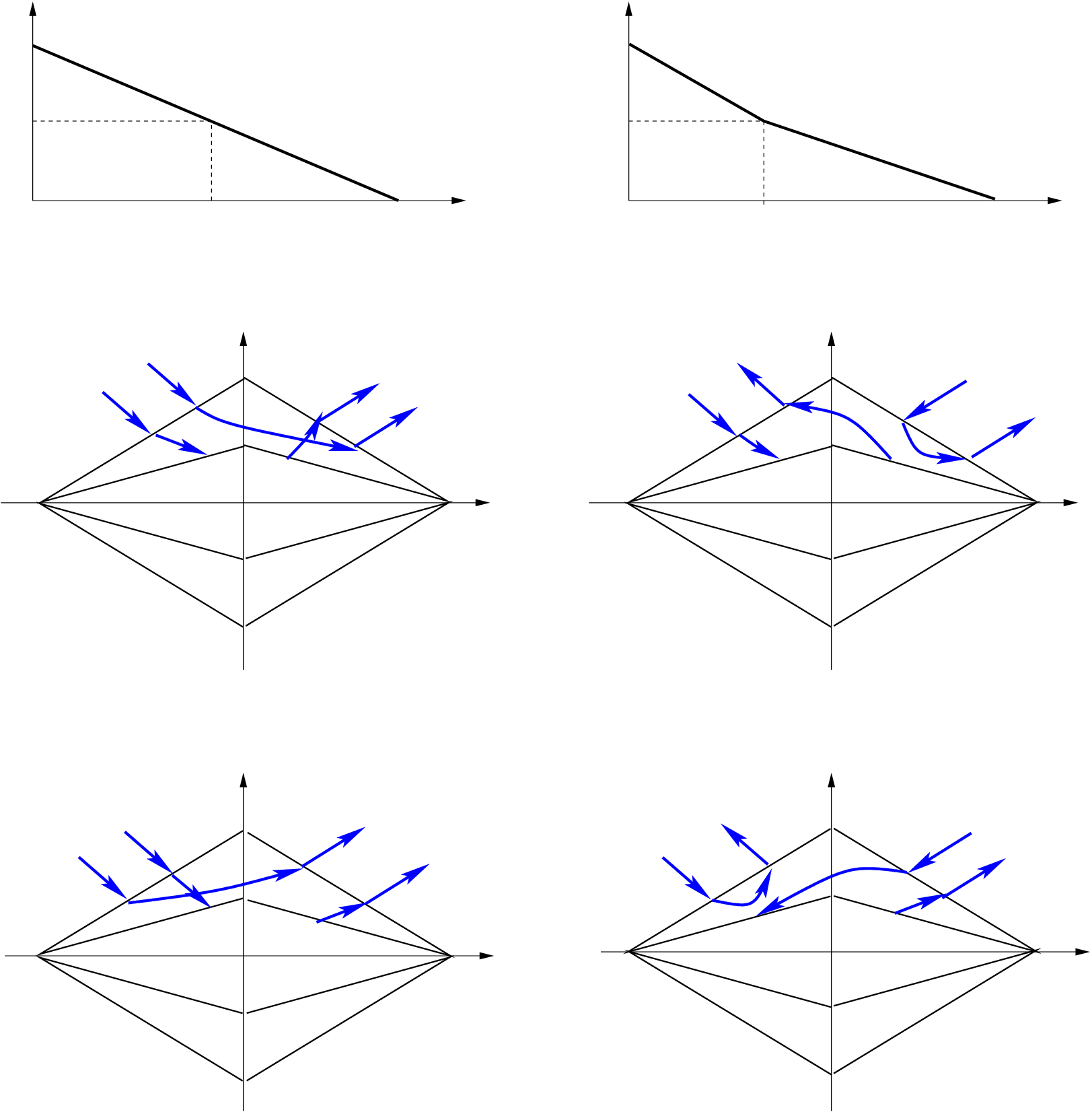}
\put(-8.7,6.7){(a)}\put(-3.8,6.7){(b)}
\put(-8.7,-0.3){(c)}\put(-3.8,-0.3){(d)}
\put(-7.4,7.1){$2l$}\put(-5.9,7.1){$4l$}
\put(-9,8){$1$}\put(-9.05,8.7){$2$}\put(-4.1,8){$1$}\put(-4.15,8.7){$2$}
\put(-2.9,7.1){$2l_1$}\put(-1.4,7.1){$2l_1+2l_2$}
\put(-7.4,8.6){$l_1=l_2=l$}\put(-2.4,8.6){$l_1<l_2$}
\put(-7.1,3.1){$\Downarrow$}\put(-2.3,3.1){$\Downarrow$}
\put(-5.9,4.1){$t<t^*$}\put(-1.1,4.1){$t<t^*$}\put(-5.9,0.3){$t>t^*$}\put(-1.1,0.3){$t>t^*$}
\end{picture}
\caption{Proof of Lemma \ref{al8}, part (5), I} \label{af2}
\end{figure}

We now describe the local behavior of $C^*$ at the point $w$. Choose local conjugation-invariant coordinates $x,y$ in a neighborhood of $w$ so that $w=(0,0)$,
$\Tor(\sigma)=\{y=0\}$, $\Tor_\R^+(P_\Delta)=\{y\ge0\}$, $C_s=\{y=\eta_sx^{l_s}+\text{h.o.t.}\}$ with $\eta_s>0$, $s=1,2$,
where without loss of generality we assume $\eta_1\ne\eta_2$, and h.o.t. denotes terms above the Newton diagram. Thus, the equation of $C_1\cup C_2$ has the Newton diagram $\Gamma$ at $w$ as shown in Figures \ref{af2}(a,b) by fat lines. We claim that $C^*$ has under $\Gamma$ the monomial $x^{2l_1-1}y$ with a nonzero coefficient and no other monomials. Since $C^*$ can be taken close to $C_1\cup C_2$, its Newton diagram either coincides or lies below $\Gamma$. The coincidence is not possible for the following reason. The curve $C^*$ passes through
$C_1\cap C_2\cap\Tor(P_\Delta)^\times$ (which consists of only nodes by Lemma \ref{al4}(2iii)) and, in case $C^*$ has Newton diagram $\Gamma$ at $w$, intersects $C_1$ at $w$ with multiplicity $\ge2l_1=(C_1\cdot C_2)_w$. Hence, by the Noether's fundamental theorem, $C^*$ lies in the subspace of $|{\mathcal L}_{P_\Delta}|$ spanned by the curves of type $C_1\cup C'_2$, $C'_2\in|C_2|$, and $C'_1\cup C_2$, $C'_1\in|C_1|$. However, the B\'ezout type restriction dictates that
$C'_1=C_1$ and $C'_2=C_2$: for example, at each point $z\in\Sing(C_1)$, the curve $C'_1$ must induce in ${\mathcal O}_{C_1,z}$ an element of the conductor ideal (cf. \cite[Theorem 4.15]{DH}), and hence $(C'_1\cdot C_1)_z\ge2\delta(C_1,z)$, at each fixed point $w_{i'}^{\sigma'}(t^*)$ of $\widehat\bw(t^*)\cap C_1$ (including $w$), we have $(C'_1\cdot C_1)_{w_{i'}^{\sigma'}(t^*)}\ge 2k_{i'}^{\sigma'}$, and at last, at the mobile point $w_{i'}^{\sigma'}(t^*)\in\widehat\bw(t^*)\cap C_1$, we have
$(C'_1\cdot C_1)_{w_{i'}^{\sigma'}(t^*)}\ge2k_{i'}^{\sigma'}-1$, which altogether amounts to
$$\sum_{z\in\Sing(C_1)}\delta(C_1,z)-C_1K_{\Tor(P_\Delta)}-1=C_1^2+1\ .$$
Further on, a monomial in an equation of $C^*$ below $\Gamma$ cannot be $x^s$, $0\le s<2(l_1+l_2)$, since $(C^*\cdot \Tor(\sigma))_w\ge2(l_1+l_2)$, and it cannot be $x^sy$ with $0\le s\le 2l_1-2$. For the latter claim, we observe that each curve $C^{(t)}$, $t\ne 0$, has $2l_1-1$ nodes in a neighborhood of $w$ and satisfies $(C^{(t)}\cdot\Tor(\sigma))_w=2(l_1+l_2)$. Thus, if the tangent line at $C^{(t)}$ to the considered family of curves is spanned by $C^{(t)}$ and $(C^*)^{(t)}$, then
$C^{(t)}$ and $(C^*)^{(t)}$ intersect in a neighborhood of $w$ with multiplicity $\ge2(2l_1-1)+2(l_1+l_2)=4l_1+2l_2-2$. Hence, $C^*$ and $C_1\cup C_2$ intersect at $w$ with multiplicity $\ge4l_1+2l_2-2$, which leaves the only possibility of the extra monomial $x^{2l_1-1}y$.

If $l_1=l_2=l$, we can write down an equation of $C^{(t)}$ in a neighborhood of $w$ in the form
$$y^2(1+O(t_1^{>0}))-(\eta_1+\eta_2+O(t_1^{>0}))x^{2l}y+x^{4l}(\eta_1\eta_2+O(t_1^{>0}))+\text{h.o.t.}
+\sum_{r=0}^{2l-1}a_{r1}(t_1)x^ry=0,$$ where h.o.t. denotes the terms above the Newton diagram,
and
$$t_1=t-t^*,\quad a_{r1}(0)=0\quad \text{for all}\quad r=0,...,2l-1\ .$$ Consider the tropical limit at $t_1\to0$.
It includes a subdivision of the triangle $T=\conv\{(0,1),(4l,0),(0,2)\}$
and certain limit curves in the toric surfaces associated with the pieces of the subdivision. Since locally $C^{(t)}$ is an immersed cylinder, we obtain that the union of the limit curves is a curve of arithmetic genus zero, which finally allows the only following tropical limit: the triangle $T$ is the unique piece of the subdivision, and
$$a_{r1}(t_1)=t_1^{\lambda(2l-r)}(a_{r1}^0+O(t_1^{>0})),\quad r=0,...,2l-1\ ,$$ where $a_{2l-1,1}^0\ne0$ in view of the above conclusion on the Newton diagram of $C^*$, and in addition,
the limit curve with Newton triangle $T$ given by
\begin{equation}y^2-(\eta_1+\eta_2)x^{2l}y+\eta_1\eta_2x^{4l}+\sum_{r=0}^{2r-1}a_{r1}^0x^ry=0\label{ae12}
\end{equation} is rational. It also follows, that ${\mathcal V}$ is smooth, regularly parameterized by $t_1$ ({\it i.e.}, $\lambda=1$),
and
its real part submersively projects onto the path $\{\widehat\bw(t)\}_{|t-t^*|<\zeta}$. By the patchworking theorem \cite[Theorems 3.1 and 4.2]{Sh3}, we uniquely restore ${\mathcal V}$ as long as we compute the coefficients $a_{r1}^0$, $0\le r\le 2l_1-2$ (here $a_{2l-1,1}^0$ is determined by $C^*$). The rationality of the curve (\ref{ae12}) can be expressed as follows. Write equation (\ref{ae12}) in the form
$y^2-2P(x)y+\eta_1\eta_2x^{4l}=0$, then resolve with respect to $y$:
$$y=P(x)\pm\sqrt{P(x)^2-\eta_1\eta_2x^{4l}}=P(x)\pm x^{2l}\sqrt{Q(1/x)^2-\eta_1\eta_2},\quad\deg Q=2l\ .$$ The rationality means that the expression
under the radical has $2l-1$ double roots (corresponding to the nodes of the curve). This means that $Q$ is a modified Chebyshev polynomial:
$$Q(u)=\lambda_1\Cheb_{2l}(u+\lambda_2)+\lambda_3,\quad\Cheb_{2l}(u)=\cos(2l\cdot\arccos u)\ ,$$
where $\lambda_1,\lambda_2,\lambda_3$ can be computed out of $\eta_1,\eta_2$, and $a_{2l-1,1}^0$. As noticed in \cite[Proof of Proposition 6.1]{Sh0}, there are exactly two real solutions, one corresponding to a curve with the real part, shown in Figure \ref{af2}(c) by lines inside the four Newton triangles designating four real quadrants, and this curve has one hyperbolic node and $2l-2$ non-real nodes, while the other curve, shown in Figure \ref{af2}(d), has $2l-1$ elliptic nodes.
Bringing the complex orientations to the play, we see that the former solution appears when the local real branches $(\R C_1,w)$ and $\R C_2,w)$ are
oriented coherently,
while the latter solution corresponds to the case when these real branches have opposite orientations,
see Figures \ref{af2}(c,d). That is, on both sides of the path
$\{\bw(t)\}_{|t-t^*|<\zeta}$ we encounter four elements of $\overrightarrow{\mathcal M}_{0,n}^{\R,+}(\Delta,\widehat\bw(t))$ having respective quantum
indices
$$\kappa_1+\kappa_2,\ -\kappa_1-\kappa_2,\ \kappa_1-\kappa_2,\ \kappa_2-\kappa_1.$$ It remains to notice that the Welschinger signs (\ref{ae13}) and (\ref{ae9i}) of the elements
with quantum index $\kappa_1+\kappa_2$ or $-\kappa_1-\kappa_2$ are the same on both sides, since the modified limit curves 
shown in Figure \ref{af2}(c) for $t>t^*$ and $t<t^*$
are obtained from each other by reflection with respect to the vertical axis. The same holds for the elements with quantum index $\kappa_1-\kappa_2$ or $\kappa_2-\kappa_1$ (see Figure \ref{af2}(d)). This completes the proof of the constancy of the functions $W_0^\kappa(t)$ and $\widetilde W_0^\kappa(t)$ in the considered case $l_1=l_2$.

\begin{figure}
\setlength{\unitlength}{1cm}
\begin{picture}(14,7.5)(-3,-0.2)
\includegraphics[width=0.7\textwidth, angle=0]{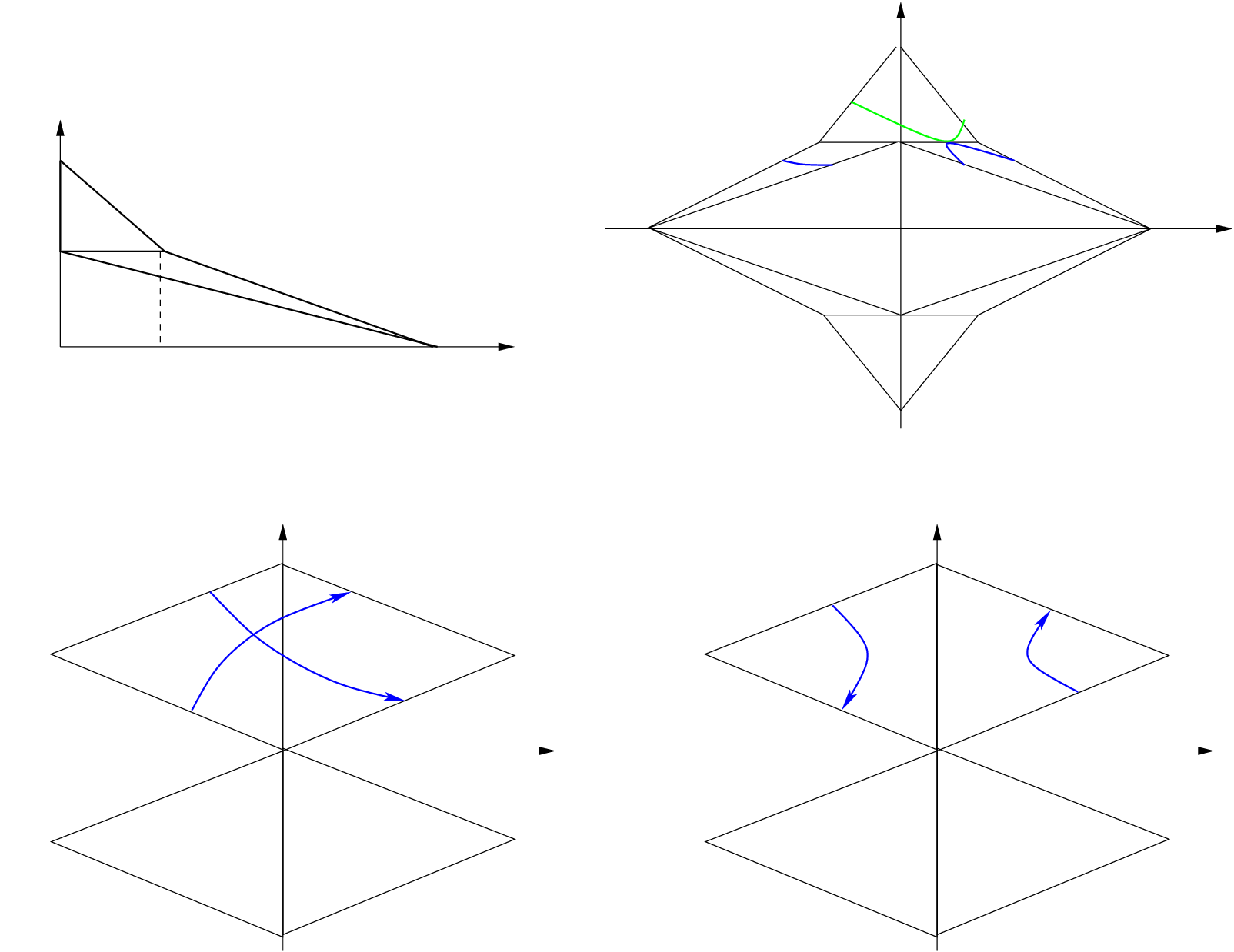}
\put(-10,4.2){(a)}\put(-4.5,4.2){(b)}
\put(-10,-0.3){(c)}\put(-4.5,-0.3){(d)}
\put(-9.2,4.7){$2l_1$}\put(-7.5,4.7){$2l_1+2l_2$}
\put(-10.3,5.8){$1$}\put(-10.3,6.6){$2$}
\end{picture}
\caption{Proof of Lemma \ref{al8}, part (5), II} \label{af3}
\end{figure}

In the case $l_1<l_2$, an equation of $C^{(t)}$ in a neighborhood of $w$ takes the form
$$(1+O(t_1^{>0}))y^2-(\eta_1+O(t_1^{>0}))x^{2l_1}y+(\eta_1\eta_2+O(t_1^{>0}))x^{2(l_1+l_2)}
+\sum_{r=0}^{2l-1}a_{r1}(t_1)x^ry+\text{h.o.t.}=0,$$
where $t_1=t-t^*$, $a_{r1}(0)=0$ for all $r=0,...,2l-1$, and h.o.t. stands for the terms above the Newton diagram. In this situation, due to the condition of arithmetic genus zero, the tropical limit is defined uniquely: the area under the Newton diagram is divided into two triangles
$$\conv\{(0,1),(2l_1,1),(2(l_1+l_2),0)\}\quad\text{and}\quad\conv\{(0,1),(2l_1,0),(0,2)\}$$ (see Figure \ref{af3}(a)), while the monomials
on the segment $[(0,1),(2l_1,1)]$ sum up to $\eta_1(x+\lambda t_1)^{2l_1}y$, where $\lambda\ne0$ is uniquely determined by the coefficient $a_{2l_1-1,1}^0$ coming from $C^*$. The two limit curves have branches intersecting each other with multiplicity $2l_1$ along the toric divisor $\Tor([(0,1),(2l_1,0)])$, see Figure \ref{af3}(b).
The genuine geometry of $C^{(t)}$ can be recovered when we deform that intersection point into $2l_1-1$ nodes. The
modification
describes such a deformation as a replacement of the intersection point by one of the $2l_1$ modified limit curves,
among which exactly two are real,
and their real parts are shown in Figure \ref{af3}(c,d) (cf. Figures \ref{af2}(c,d)): one of them has
a hyperbolic node and $2l_1-2$ non-real nodes, and the other has
$2l_1-1$ elliptic nodes. As in the preceding paragraph, the former modified limit curve fits the case
when the local real branches
$(\R C_1,w)$, $(\R C_2,w)$ are
oriented coherently,
which means that on each side of the path $\{\widehat\bw(t)\}_{|t-t^*|<\zeta}$ we have two elements of
$\overrightarrow{\mathcal M}_{0,n}^{\R,+}(\Delta,\widehat\bw(t))$
with quantum indices $\kappa_1+\kappa_2,-\kappa_1-\kappa_2$,
while all four elements have the same Welschinger sign
(\ref{ae13}) or (\ref{ae9i}) since the constructions for $t>t^*$ and $t<t^*$
are symmetric with respect to the vertical axis.
Hence, we obtain
the constancy of the functions $W_0^\kappa(t)$, $\widetilde W_0^\kappa(t)$. The same holds in the case of the opposite orientation of the branches $(\R C_1,w)$, $(\R C_2,w)$ with the use of the second modified limit curve. We also point out that the final deformation including the modification still depends on $t_1=t-t^*$.

The realizability of the deformations described in this item by variation of the constraints requires an independent variation of the coefficients of the monomials under the Newton diagram $\Gamma$ while keeping the total $\delta$-invariant outside the point $w$ and keeping the given intersection number $2k_{i'}^{\sigma'}$ with the toric divisor $\Tor(\sigma')$ at a point close to $w_{i'}^{\sigma'}$, for each $(i',\sigma')\ne(i,\sigma)$. This can be expressed as
\begin{equation}h^1(C,{\mathcal J}_{Z/C}\otimes_{{\mathcal O}_C}\ina^*{\mathcal L}_{P_\Delta})=0,\label{ehorace2}\end{equation}
where $C=C_1\cup C_2$, $C_r=\bn(\widehat C_r)$, $r=1,2$, the zero-dimensional scheme $Z\subset C$ is supported at $\Sing(C)\cup\widehat\bw(t^*)$ and defined by the ideal $J^{cond}(C,z)$ at each $z\in\Sing(C)\cap\Tor(P_\Delta)^\times$, by the ideal
$\{\varphi\in{\mathcal O}_{C,w_{i'}^{\sigma'}}\ :\ (\varphi\cdot C)_{w_{i'}^{\sigma'}}\ge 2k_{i'}^{\sigma'}-1\}$ at each $w_{i'}^{\sigma'}\in\widehat\bw(t^*)\setminus\{w\}$, and by the ideal generated by the monomials on the Newton diagram $\Gamma$ and above it at the point $w$. This $h^1$-vanishing lifts to the disjoin union of $\widehat C_1$ and $\widehat C_2$ in the form
\begin{equation}h^1(\widehat C_r,{\mathcal O}_{\widehat C_r}(\bd_r))=0,\quad r=1,2,\label{ecert}\end{equation}
where
$$\deg\bd_r=C^2_r-(C^2_r+C_rK_{\Tor(P_\Delta)}+2)-(-C_rK_{\Tor(P_\Delta)}-|C_r\cap\Tor(\partial P_\Delta)|+1)$$
\begin{equation}=|C_r\cap\Tor(\partial P_\Delta)|-3\ge-1>-2,\quad r=1,2,\label{ehorace3}\end{equation}
which certifies (\ref{ecert}).

\smallskip({\it 5b}) The remaining case is as follows: $\bn:\widehat C_1\to\Tor(P_\Delta)$ is an $2l_1$-multiple covering of a line $C_1\simeq\PP^1$ through $w$ and
$\bn:\widehat C_2\to\Tor(P_\Delta)$ is a (birational) immersion onto a rational curve $C_2$ having a smooth branch at $w$ intersecting $\Tor(\sigma)$ with multiplicity $2l_2$. Here,
we remove the restriction $l_1\le l_2$. We extend the range of the coordinate system $x,y$ to a neighborhood of the line $C_1$, assuming that $C_1=\{x=0\}$.
Then, the Newton diagram of $\bn_*(\widehat C)$
in this neighborhood of $C_1$
is as shown in Figure \ref{af4}(a) by fat lines, where $N$ is the intersection multiplicity of the curves $C_1,C_2$.
Following the recipe used in the preceding cases
and taking into account that the curves $C^{(t)}$, $0<|t-t^*|\ll1$, are rational, we obtain that
\begin{itemize}\item the combinatorial part of the tropical limit consists of the triangle
$\conv\{(0,1),(2l_1,1),(2(l_1+l_2),0)\}$ and the rectangle $\conv\{(0,1),(0,N),(2l_1,1),(2l_1,N)\}$ (see Figure \ref{af4}(b));
\item the limit curve associated with the triangle is a rational curvewhich intersect the toric divisor $\Tor([(0,1),(2l_1,1)])$ at one point, where it is smooth, while the limit curve associated with the rectangle is split into $2l_1$-multiple line intersecting the two horizontal toric divisors and $N-1$ distinct lines intersecting the two vertical toric divisors; here, the position of $N-1$ horizontal lines is determined by the intersection points of $C_2$ with $\Tor([(2l_1,1),(2l_1,N)])$, and the position of the vertical line is determined by the intersection point of he limit curve from the preceding item with $\Tor([(0,1),(2l_1,1)])$.
\end{itemize}
Note that the $2l_1$-multiple component of the limit curve associated with the rectangle is, in fact, the image of a $2l_1$-multiple covering $\bn:\widehat C_1\to C_1$, and the point $p$ is mapped to the point $x=\tau$ on the toric divisor $\Tor([(0,1),(2l_1,1)])$. The deformation can be written in the form
$$ay(x-\tau)^{2l_1}\prod_{j=1}^{N-1}(y-y_j)+f(x,y)+\text{h.o.t.}=0\ ,$$
where $\tau=t-t^*$, the non-zero constants $a,y_1,...,y_{N-1}$ are determined by $C_2$, and $f(x,y)$ is the sum of monomials strictly right to the Newton diagram with coefficients determined by $C_2$.
The higher order terms (h.o.t.) yield the appearance of $2l_1(N-1)$ nodes of $C^{(t)} (t\ne t^*$) in a neighborhood of $C_1$, local equigeneric deformation of each singular point of $C_2$, the fixed intersection multiplicity $2l_1+2l_2$ of $C^{(t)}$ and $\Tor(\sigma)$ at $w$, and intersection multiplicities of $C^{(t)}$ with $\Tor(\partial P_\Delta)\setminus\{w\}$ prescribed by the condition $C^{(t)}\in{\mathcal M}_{0,n}(\Delta,\bw(t))$.
Since this deformation does not produce elliptic nodes in a neighborhood of $C_1$, the functions $W_0^\kappa(\widehat\bw(t))$ and $\widetilde W_0^\kappa(\widehat\bw(t)$ remain constant for $t\in(t^*-\eps,t^*+\eps)$.

\begin{figure}
\setlength{\unitlength}{1cm}
\begin{picture}(12,6)(0,0)
\thinlines
\put(1,1){\vector(0,1){5}}\put(1,1){\vector(1,0){5}}
\put(7,1){\vector(0,1){5}}\put(7,1){\vector(1,0){5}}

\dottedline{0.1}(1,2)(3,2)\dottedline{0.1}(3,1)(3,2)
\dottedline{0.1}(9,1)(9,2)\dottedline{0.1}(1,5)(3,5)
\dashline{0.2}(3,5)(5,4.5)

\thicklines
\put(3,2){\line(2,-1){2}}\put(3,2){\line(0,1){3}}
\put(7,2){\line(4,-1){4}}\put(7,2){\line(1,0){2}}
\put(9,2){\line(2,-1){2}}\put(9,2){\line(0,1){3}}\put(7,5){\line(1,0){2}}

\put(2.8,0.6){$2l_1$}\put(4.4,0.6){$2l_1+2l_2$}
\put(0.6,1.9){$1$}\put(0.5,4.9){$N$}
\put(6.6,1.9){$1$}\put(6.5,4.9){$N$}
\put(8.8,0.6){$2l_1$}\put(10.4,0.6){$2l_1+2l_2$}
\put(1.2,0){(a)}\put(7.2,0){(b)}

\end{picture}
\caption{Proof of Lemma \ref{al8}, part (5), III}\label{af4}
\end{figure}

A sufficient condition for the existence of the required deformation can be written as
\begin{equation}h^1(C,{\mathcal J}_{Z/C}\otimes_{{\mathcal O}_C}\ina^*{\mathcal L}_{P_\Delta})=0,\label{ehorace}\end{equation}
where $C$ is the non-reduced union of a $2l_1$-multiple component $C_1$ and of $C_2$, and the zero-dimensional scheme $Z\subset C$ is supported at $\Sing(C_2)\cup(C\cap\Tor(\partial P_\Delta))\cup(C_1\cap C_2)$ and is defined by the ideal $I^{cond}(C_2(z)$at each $z\in\Sing(C_2)$, by the ideal $I_z=\langle x^{2l_1},y-y_i\rangle$ at each point $(0,y_i)\in C_1\cap C_2$, $i=1,...,N_1$, by the ideal
$$\{\varphi\in{\mathcal O}_{C_2,w_{i'}^{\sigma'}}\ :\ (\varphi\cdot C_2)_{w_{i'}^{\sigma'}}\ge 2k_{i'}^{\sigma'}-1\}$$
at some point $w_{i'}^{\sigma'}\in (C_2\cap\Tor(\partial P_\Delta)))\setminus\{w\}$, by the ideal
$$\{\varphi\in{\mathcal O}_{C_2,w_{i''}^{\sigma''}}\ :\ (\varphi\cdot C_2)_{w_{i''}^{\sigma''}}\ge 2k_{i''}^{\sigma''}\}$$at each point $w_{i''}^{\sigma''}\in C_2\cap\Tor(\partial P_\Delta)\setminus\{w,w_{i'}^{\sigma'}\}$,
and finally at the point $w$, by the ideal $I_w$ determined by the vanishing of the coefficients of the monomials $x^r$, $0\le r<2l_1+2l_2$, $x^ry$, $0\le r<2l_1$, see Figure \ref{af4}. We establish (\ref{ehorace}) using the basic Horace method \cite{Hir} (see also \cite[Section 3.4.1]{GLS1}): the vanishing relation (\ref{ehorace}) follows from the series of vanishing relations
\begin{equation}\begin{cases} & h^1(C_1,{\mathcal J}_{(Z:jC_1)\cap C_1/C_1}\otimes_{{\mathcal O}_{C_1}}\ina^*{\mathcal L}_{P_\Delta}(-jC_1))=0,\quad j=0,...,2l_1-1,\\
& h^1(C_2,{\mathcal J}_{Z:(2l_1C_1)/C_2}\otimes{{\mathcal O}_{C_2}}\ina^*{\mathcal L}_{P_\Delta}(-2l_1C_1))=0.\end{cases}\label{ehorace1}\end{equation}
Each of the $2l_1$ relations in the upper line of (\ref{ehorace1}) reads
$$h^1(C_1,{\mathcal O}_{C_1}(-w))=0,$$
which follows by the Riemann-Roch theorem.
The residue scheme $Z:(2l_1C_1)$
in the lower line of (\ref{ehorace1}) coincides with $Z$ at
$\Sing(C_2)\cup(C_2\cap\Tor(\partial P_\Delta))\setminus\{w\}$
and is defined at $w$ by the ideal $I_w=\langle x^{2l_2},y\rangle$.
Then, it lifts to $\widehat C_2\simeq\PP^1$ in the form
$$h^1(\widehat C_2,{\mathcal O}_{\widehat C_2}(\bd))=0,$$
where
$$\deg \bd=C_2^2-(C^2_2+C_2K_{\Tor(P_\Delta)}+2)-(-C_2K_{\Tor(P_\Delta)}-1)$$ 
$$=-1 
>-2,$$
and hence (\ref{ehorace1}) together with (\ref{ehorace}) follow.

\smallskip
{\bf(6)} Suppose that $\xi$ is as in Lemma iv). Without loss of generality we can suppose that
$$P=\conv\{(0,2m),\ (2p,0),\ (2q,0)\},\quad 0\le p<q,\ m<q\ ,$$ and that $\bn$ is ramified at the points
of $\widehat\bw(t^*)$ on the toric divisors associated with the segments $[(0,2m),(2p,0)]$ and
$[(0,2m),(2q,0)]$, while the two points $w_1^\sigma(t),w_2^\sigma(t)\in\bw(t)$ on the toric divisor $\Tor(\sigma)$, where $\sigma=[(2p,0),(2q,0)]$, merge to one point as $t\to t^*$.
Then, $C=\bn(\PP^1)$ admits a parametrization
$$x=a\theta^m,\quad y=b\theta^p(\theta-1)^{q-p},\quad\theta\in\C\ ,$$ and correspondingly $\bn$ is given by
$$x=a\theta^{2m},\quad y=b\theta^{2p}(\theta^2-1)^{q-p},\quad\theta\in\C\ .$$
Assuming that the path $\{\widehat\bw(t)\}_{|t-t^*|<\eps}$, is defined by fixing the point $w_1^\sigma$,
associated with $\theta=1$, and
the point of $\widehat\bw(t^*)$ on $\Tor([(0,2m),(2q,0)])$, we obtain a one-parameter deformation of $\bn$
\begin{equation}x=a\theta^{2m},\quad y=b\theta^{2p}(\theta-1)^{q-p}(\theta+1+\lambda)^{q-p},
\quad|\lambda|\ll1\ .\label{ae98}\end{equation} This deformation is regularly parameterized by the difference of the $x$-coordinates of $w_2^\sigma(t)$ and $w_1^\sigma(t)$ equal to $2ma\lambda+O(\lambda^2)$.
It follows from formulas (\ref{ae98}) that the corresponding element $\xi(t)\in\overline{\mathcal M}_{0,n}^{\R,+}(\Delta,\widehat\bw(t))$
has the same orientations at a point on the toric divisor $\Tor([(0,2m),(2p,0)])$ for $\lambda<0$ and for $\lambda>0$, and so does for the point on the toric divisor
$\Tor([(0,2m),(2q,0)])$, while the local branches at $w_1^\sigma(t)$, $w_2^\sigma(t)$ have opposite orientations with respect to the orientation of
$\Tor_\R(\sigma)$. Furthermore, $\xi(t)$ has exactly two elliptic nodes in a neighborhood of each elliptic node of $C$.
Thus, the constancy of $W_0^\kappa(t)$ and $\widetilde W_0^\kappa(t)$ follows.
\proofend

\section{Invariance in genus one: proof of Theorem \ref{at2}}\label{asec5}

Take
two sequences $\widehat\bw(0)=(\widehat\bw_\partial(0),w_0(0))$ and $\widehat\bw(1)=(\widehat\bw_\partial(1),w_0(1))$,
satisfying the conditions of Theorem \ref{at2}. We may assume that
$\widehat\bw_\partial(0),\widehat\bw_\partial(1)\in \Men_\R^\rho(\Delta)$ for some $\rho>0$.
Then,
we join these sequences by a generic path $\{\widehat\bw_\partial(t)\}_{0\le t\le 1}$, in $\Men_\R^\rho(\Delta)$ and
the points $w_0(0)$, $w_0(1)$ by a generic path $\{w_0(t)\}_{0\le t\le1}$ in ${\mathcal Q}$, and verify the constancy
of $W_1^\kappa(\Delta,\widehat\bw(t))$ and $\widetilde W_1^\kappa(\Delta,\widehat\bw(t))$
in all possible wall-crossing events, for all quantum indices $\kappa$.

\subsection{Invariance inside chambers}\label{chamb1}
Similarly to the proof of Theorem \ref{at1}, we start with the verification
of the constancy of $W_1^\kappa(\Delta,\widehat\bw(t))$
along intervals $t'<t<t''$ such that $\{\widehat\bw(t)\}_{t'<t<t''}$ lies in the same chamber of $\Men_1(\Delta)$.
To this end, it is enough to show that the projection ${\mathcal M}_{1,n}^{\;\R,+}(\Delta,\widehat\bw(t))_{t'<t<t''}\to(t',t'')$ has no critical points (cf. the proof of Theorem \ref{at1}).
This requirement amounts to the following statement (cf., Lemma \ref{lis2} and \cite[Lemma 2.3]{Sh}).

\begin{lemma}\label{lis3}
Let $\xi=[\bn:(\widehat C,\bp)\to\Tor(P_\Delta)]\in{\mathcal M}_{1,n}^{\;\R,+}(\Delta,\widehat\bw)$, where $\bw\in\Men_1(\Delta)$. Then,
\begin{equation}H^0\left(\widehat C,{\mathcal N}_\bn\left(-p_0-\sum_{\sigma\ne\sigma_0,i\ne i_0}2k^\sigma_ip^\sigma_i-(2k^{\sigma_0}_{i_0}-1)p^{\sigma_0}_{i_0}\right)\right)=0.\label{eis3}\end{equation}
\end{lemma}

{\bf Proof.}
We proceed along the lines of the proof of \cite[Lemma 2.3]{Sh}. Assume that $H^0\ne0$ in (\ref{eis3}).
Then, there exists a real curve $C'\in|{\mathcal L}_{P_\Delta}|\setminus\{C\}$, where $C=\bn_*(\widehat C)$,
which intersects $C$
\begin{itemize}\item at each point $q\in\Sing(C)$ with multiplicity $\ge2\delta(C,q)$,
\item at each point $w^\sigma_i$, $(\sigma,i)\ne(\sigma_0,i_0)$ with multiplicity $\ge 2k^\sigma_i$,
\item at $w^{\sigma_0}_{i_0}$ with multiplicity $\ge 2k^{\sigma_0}_{i_0}-1$,
\item and at $w_0$.
\end{itemize}
Since $C$ has in $\overline {\mathcal Q}$ a null-homologous immersed circle $S$, there must
be an additional intersection point $w'\in C'\cap S$, not mentioned in the above list. However, then
$$C'C\ge2\sum_{q\in\Sing(C)}\delta(C,q)+\sum_{\sigma\in P^1_\Delta}\sum_{i=1}^{n^\sigma}2k^\sigma_i-1+2$$
$$=(C^2-c_1(\Tor(P_\Delta)[C]))+c_1(\Tor(P_\Delta)[C]-1+2=C^2+1,$$
which is a contradiction.
\proofend

\subsection{Walls in the space of constraints}\label{walls1}
The set of sequences $\widehat\bw\in \Men_\R^\rho(\Delta)\times {\mathcal Q}$ satisfying
$$\overline{\mathcal M}^\R_{1,n}(\Delta,\widehat\bw)={\mathcal M}^\R_{1,n}(\Delta,\widehat\bw)\quad\text{and}\quad \widehat\bw\in \Men_1(\Delta),$$
is a dense semialgebraic subset of full dimension
$\dim \Men_\R^\rho(\Delta)\times {\mathcal Q}=n+1$.
The complement is the union of finitely many semialgebraic strata of codimension $\ge1$. Since the path $\{(\widehat\bw_\partial(t),w_0(t))\}_{0\le t\le 1}$ is generic, it avoids strata of $\Men_\R^\rho(\Delta)\times {\mathcal Q}$ of codimension $\ge2$ and intersects strata of codimension one only in their generic points.
Note that the image $\pr(S)\subset|{\mathcal L}_{P_\Delta)}|$ of $\overline{\mathcal M}^\R_{1,n}(\Delta,\widehat\bw)$ in the linear system $|{\mathcal L}_{P_\Delta}|$ when $\widehat\bw$ ranges over a stratum $S$ of codimension one, has dimension $n-1$.

Now, we study the strata $S$ of the dimension $n$ in $(\Men_\R^\rho\times {\mathcal Q})\setminus \Men_1(\Delta)$
in terms of the geometry of elements
$\xi\in\overline{\mathcal M}_{1,n}^{\;\R}(\Delta,\widehat\bw(t))\setminus{\mathcal M}_{1,n}^{\;\R}(\Delta,\widehat\bw(t))$, $0\le t\le1$.
In view of the results of Section \ref{chamb1}, we can ignore the strata along which
$\overline{\mathcal M}^\R_{1,n}(\Delta,\widehat\bw)={\mathcal M}^\R_{1,n}(\Delta,\widehat\bw)$.
The following lemma classifies other strata.

\begin{lemma}\label{al5}
(1) The following elements $\xi=[\bn:(\widehat C,\bp)\to\Tor(P_\Delta)]$ cannot occur in $\overline{\mathcal M}_{1,n}^{\;\R}(\Delta,\widehat\bw(t))\setminus{\mathcal M}_{1,n}^{\;\R}(\Delta,\widehat\bw(t))$, $0\le t\le1$, as $\widehat\bw\in S$:
\begin{enumerate}\item[(1i)]
$\widehat C$ is a reducible curve of arithmetic genus $1$ having a component mapped onto a toric divisor;
\item[(1ii)]
$\widehat C$ is a connected curve of arithmetic genus $1$ either having at least three irreducible components, or with two rational irreducible components
    and no other components.
\end{enumerate}

(2) If $\xi=[\bn:(\widehat C,\widehat\bp)\to\PP^2]\in\overline{\mathcal M}_{1,n}(\Delta,\widehat\bw(t^*))$ when $\widehat\bw(t^*)$ does not satisfy the conclusions of Lemma \ref{l3-1a}, then $\xi$ is of one of the following types:

(2a) either $w_0(t^*)$ is a singular point of the curve $C=\bn(\widehat C)$, and the following holds:
\begin{enumerate}
\item[(2i)] $\widehat\bw_\partial(t^*)$ consists of $n$ distinct points, the curve $\widehat C$ is smooth elliptic,
$\bn$ is an immersion onto a curve $C$ which is smooth along
the toric divisors, while $w_0(t^*)$ is a center of at least one real local branch;
\end{enumerate}

(2b) or $w_0(t^*)$ is a smooth point of the curve $C=\bn(\widehat C)$, and one of the following holds:
\begin{enumerate}
\item[(2ii)] $\widehat\bw_\partial(t^*)$ consists of $n$ distinct points, $\widehat C$ is a smooth elliptic curve,
the map $\bn$ is birational onto its image $C$ but not an immersion:
it may have singular local branches centered in $\Tor(P_\Delta)^\times$ as well as singularities of type $A_{2m}$ at some points $w_i^\sigma\in\widehat\bw_\partial(t^*)$, where $m\ge k_i^\sigma$;
\item[(2iii)] $\widehat\bw_\partial(t^*)$ consists of $n-1$ distinct points, in particular,
$w_i^\sigma(t^*)=w_j^\sigma(t^*)=:w$ for some $\sigma\in P^1_\Delta$ and $i\ne j$,
the curve $\widehat C$ is smooth elliptic, and the map $\bn$
is birational onto its image $C$; the curve $C$ is immersed outside $\widehat\bw_\partial(t^*)\setminus\{w\}$ and may have singularities of type
$A_{2m}$ at some points $w_{i'}^{\sigma'}\in\widehat\bw_\partial(t^*)\setminus\{w\}$, while at $w$,
the curve $C$ either is smooth or has two smooth branches;
\item[(2iv)] $\widehat\bw_\partial(t^*)$ consists of $n$ distinct points, $\widehat C=\widehat C_1\cup\widehat C_2$, where $\widehat C_1\simeq\PP^1$ and
$\widehat C_2$ is a smooth elliptic curve,
the intersection
$\widehat C_1\cap\widehat C_2$
consists of one point $p$;
the map $\bn:\widehat C_1\to\Tor(P_\Delta)$, $s=1,2$, is either an immersion, smooth along
$\Tor(\partial P_\Delta)$, or a multiple covering of a line intersecting only two toric divisors, while these divisors correspond to opposite parallel sides of $P_\Delta$ and the intersection points with these divisors are ramification points of the covering; the map $\bn:\widehat C_2\to\Tor(P_\Delta)$ is an immersion, smooth along $\Tor(\partial P_\Delta)$;
the point $p$ either is mapped to $\Tor(P_\Delta)^\times$, and then the curves $C_1=\bn(\widehat C_1)$, $C_2=\bn(\widehat C_2)$ intersect each other only in $\Tor(P_\Delta)^\times$ and each of their intersection points is an ordinary node, or $p$ is mapped to some point $w_i^\sigma(t^*)$, in which case
the curves $C_i=\bn(\widehat C_i)$, $i=1,2$, intersect $\Tor(\sigma)$ with even multiplicities, and they do not have other common point in $\widehat\bw_\partial(t^*)$;
\item[(2v)] $\widehat\bw_\partial(t^*)$ consists of $n$ distinct points, $\widehat C$ is an irreducible rational curve with one node $p$,
the map $\bn$ is an immersion
sending $p$ to some point $w_i^\sigma(t^*)$, center of two smooth branches
intersecting
with multiplicity $\min\{2l_1+1,2l_2+1\}$ and intersecting
$\Tor(\sigma)$ with odd multiplicities
$2l_1+1, 2l_2+1$,
where $l_1+l_2+1=k_i^\sigma$;
furthermore, each of the other points of $\widehat\bw_\partial(t^*)$ is a center of one smooth branch;
\item[(2vi)] $\widehat C$ is a smooth elliptic curve, $\widehat\bw_\partial(t^*)$ consists of $n=4$ points, the curve $C=\bn(\PP^1)$ is rational
and smooth at $\widehat\bw_\partial(t^*)$, and
$\bn:\widehat C\to C$ factors through the double covering $\widehat C\to C^\nu\simeq\PP^1$ ramified at four points projecting to $\widehat\bw_\partial(t^*)$, where $C^\nu\to C$ is the normalization.
\end{enumerate}
\end{lemma}

{\bf Proof.} (1i) In
this case,
we must have $\widehat C=\widehat C'\cup\widehat C''$, where $\widehat C'$ is a connected curve
of arithmetic genus one such that $\widehat C'$ is mapped onto the union of all toric divisors, and $\widehat C''$
is a non-empty (due to $w_0\in\Tor(P_\Delta)^\times$) union of connected curves of arithmetic genus zero, each one joined with $\widehat C'$ in one point (cf. the proof of Lemma \ref{al4}(1i)).
Let $\widehat C''_0$ be a connected component of $\widehat C''$ such that $w_0(t^*)\in\bn(\widehat C''_0)$.
Note that all local branches of
$\bn:\widehat C''_0\to\Tor(P_\Delta)$ centered on toric divisors, except for $\bn:(\widehat C''_0,p)\to\Tor(P_\Delta)$, where $p=\widehat C'\cap\widehat C''_0$,
are, in fact, centered at some points of $\widehat\bw_\partial(t^*)$. Observe that either $C''_0$ intersects at least three toric divisors of $\Tor(P_\Delta)$, or it intersects two toric divisors corresponding to parallel sides of $P_\Delta$. However, in the both cases we reach a contradiction with (AQC).

\smallskip
(1ii) In the case of at least $3$ irreducible components, we follow the lines of the proof of Lemma \ref{al4}(1ii). If $\widehat C$ contains an elliptic
component and at least $2$ rational components, then we similarly obtain at least three
independent Menelaus type conditions for $\widehat\bw_\partial(t^*)$, which bounds from above the dimension of the considered stratum by $(n-3)+2=n-1<n$, a contradiction.
If all irreducible components of $\widehat C$ are rational, then we encounter
at least two independent Menelaus conditions on $\widehat\bw_\partial(t^*)$ as well as a condition on the position of the point $w_0(t^*)$ due to the finiteness statement of Lemma \ref{l3-1a}, which altogether bounds from above the dimension of the considered stratum
by $(n-2)+1=n-1<n$.

\smallskip
(2a) Let $w_0(t^*)$ be a singular point of $C$. Then, for the dimension reason, $\widehat\bw_\partial(t^*)$ must be a generic element of $\Men_\R^\rho$.
This implies, in particular, that $\widehat C$ is irreducible and $\widehat\bw_\partial(t^*)$ consists of $n$ distinct points. Furthermore, $\widehat C$ cannot be a rational curve with a node.
Thus, $\widehat C$ is a smooth elliptic curve. Then (see Lemma \ref{l3-1a}),
we derive that $\bn:\widehat C\to\Tor(P_\Delta)$ is an immersion and $C$ is smooth along the toric divisors.
Since $w_0(t^*)$ turns into a real smooth point in the deformation along the path $\{(\widehat\bw_\partial(t),w_0(t))\}_{0\le t\le 1}$,
it must be a center of at least one real local branch.

\smallskip
(2b) From now on we can assume that $w_0(t^*)$ is a smooth point of $C=\bn_*(\widehat C)$.

\smallskip(2ii)
Suppose that $\widehat\bw_\partial(t^*)$ consists of $n$ distinct points, $\widehat C$ is a smooth elliptic curve, $\bn:\widehat C\to C=\bn(\widehat C)$ is birational, but not an immersion. To verify the conditions of item (2ii) of the lemma, it is enough to consider the case of $C$ having a singular branch at some point $w_i^\sigma\in\widehat\bw_\partial(t^*)$. Fixing the position of the point $w_i^\sigma$, we obtain a subfamily $\pr(S^{fix}(w_i^\sigma))\subset\pr(S)$ of dimension $\ge n-2\ge2$ and then apply \cite[Inequality (5) in Lemma 2.1]{IKS4} (cf. (\ref{ae9f})):
$$-CK_{\Tor(P_\Delta)}\ge\left(-CK_{\Tor(P_\Delta)}
-n+1\right)+(n-3)+(\mt Q-1)$$
\begin{equation}=-CK_{\Tor(P_\Delta)}-2+(\mt Q-1)\quad\Longrightarrow\quad\mt Q\le3\ ,\label{ae9fX}\end{equation}
where $Q$ is the local branch of $C$ at $w_i^\sigma$. It remains to show that $\mt Q=2$. Arguing by contradiction, assume that
$\mt Q=3$, {\it i.e.}, $Q$ is topologically equivalent to a singularity $y^3+x^m$,
where $\gcd(3,m)=1$ and $m\ge 2k_i^\sigma\ge4$. Without loss of generality, we can assume that the germ of $\pr(S)\subset|{\mathcal L}_{P_\Delta}|$ at $C$ is an equisingular family.
According to items (a) and (b) in Section \ref{sec-def},
we obtain the following upper bound completely analogous to (\ref{edop-alg2}):
\begin{equation}\dim\pr(S^{fix}(w_i^\sigma))\le h^0(C,{\mathcal J}_{Z/C}\otimes_{{\mathcal O}_C}\ina^*{\mathcal L}_{P_\Delta}=h^0\left(\widehat C,{\mathcal O}_{\widehat C}(-\bd)\otimes_{{\mathcal O}_{\widehat C}}\bn^*{\mathcal L}_{P_\Delta}\right),\label{ae9l}\end{equation}
where
$$\deg\bd=2\sum_{z\in\Sing(C)}\delta(C,z)+(-CK_{\Tor(P_\Delta)}-n+1)+2$$
$$\Longrightarrow\quad\deg({\mathcal O}_{\widehat C}(-\bd)\otimes_{{\mathcal O}_{\widehat C}}\bn^*{\mathcal L}_{P_\Delta})=C^2-\deg\bd=n-3\ge1.$$
By the Riemann-Roch formula, this yields
$$\dim\pr(S^{fix}(w_i^\sigma))\le n-3$$
contradicting the lower bound $\dim\pr(S^{fix}(w_i^\sigma))\ge n-2$ pointed above.

\smallskip(2iii)
Suppose that some of the points of the sequence $\widehat\bw_\partial(t^*)$ coincide.
For the dimension reason, we immediately get that $\widehat\bw_\partial(t^*)$ consists of $n-1$ distinct points, and all these points must be
in general position subject to the unique Menelaus relation (\ref{ae1}). Hence, $\widehat C$ must be irreducible. Moreover, $\widehat C$ cannot be rational, since otherwise, by Lemma \ref{l3-1a}, one would get the total dimension of the considered stratum $\le(n-2)+1=n-1<n$, a contradiction.
Thus, $\widehat C$ is a smooth elliptic curve. We have $w_i^\sigma(t^*)=w_j^\sigma(t^*)=w$ for some $\sigma\subset\partial P$ and $i\ne j$.
If $C$ is unibranch at each point of $\widehat\bw_\partial(t^*)$, then, by Lemma \ref{l3-1a},
the curve $C$ is immersed and smooth along the toric divisors.
Otherwise, the curve $C$
has two local branches at $w$.
Moreover, it must be unibranch at each point of $\widehat\bw_\partial(t^*)\setminus\{w\}$.

Let us analyze possible singularities of $C$.
Fixing the position of $w$ and some other point $w'\in\widehat\bw_\partial(t^*)\setminus\{w\}$, we obtain a family
$\pr(S^{fix}(w,w'))\subset\pr(S)$ of dimension $\ge n-3\ge1$; hence, \cite[Inequality (5) in Lemma 2.1]{IKS4} applies and similarly to (\ref{ae9f}) yields
$$-CK_{\Tor(P_\Delta)}\ge\left(-CK_{\Tor(P_\Delta)}
-n+3\right)+(n-4)+\sum_B(\mt Q-1)$$
\begin{equation}=-CK_{\Tor(P_\Delta)}-1+\sum_B(\mt Q-1)\quad\Longrightarrow\quad\sum_B(\mt Q-1)\le1\ ,\label{e2iii1}\end{equation}
where $B$ is the set of all singular local branches of $C$ in $\Tor(P_\Delta)^\times\cup\{w,w'\}$. That means a singular branch of $C$ (if any) must be of order $2$.
We show that $C$ is immersed outside $\widehat\bw_\partial(t^*)\setminus\{w\}$, and hence fits the claim of item (2iii).

Let $C$ have a singular local branch $Q$ centered in $\Tor(P_\Delta)^\times$. Note that $\dim\pr(S^{fix}(w))\ge n-2$.
Similarly to (\ref{ae9l}) we get
\begin{equation}\dim\pr(S^{fix}(w))\le h^0(\widehat C,{\mathcal O}_{\widehat C}(-\bd_1)\otimes_{{\mathcal O}_{\widehat C}}\bn^*{\mathcal L}_{P_\Delta}),\label{e2iii2}\end{equation}
where
$$\deg\bd_1=2\sum_{z\in\Sing(C)}\delta(C,z)+(-CK_{\Tor(P_\Delta)}-n+2)+(\mt Q-1)$$
\begin{equation}\Longrightarrow\quad\deg(\widehat C,{\mathcal O}_{\widehat C}(-\bd_1)\otimes_{{\mathcal O}_{\widehat C}}\bn^*{\mathcal L}_{P_\Delta})=C^2-\deg\bd_1=C^2-(C^2+CK_{\Tor(P_\Delta)})$$ $$-(-CK_{\Tor(P_\Delta)}-n+2)
-(\mt Q-1)=n-3\ge1,\label{e2iii3}\end{equation}
which yields $\dim\pr(S^{fix}(w))\le n-3$, and hence a contradiction
with the lower bound $\dim\pr(S^{fix}(w))\ge n-2$.

Let $C$ have a singular local branch $Q$ centered at $w$. Then $\dim\pr(S^{fix}(w))\ge n-2$, whereas the computation literally coinciding with that in the preceding paragraph yields $\dim\pr(S^{fix}(w))\le n-3$.

\smallskip(2iv) We are left with the case of $\widehat\bw_\partial(t^*)$ consisting of $n$ distinct points
and $\widehat C$ either consisting of two components, or being a rational curve with a node. Observe that the case of two rational components is not possible, since a generic element of ${\mathcal M}^\R_{0,n}(\Delta,\widehat\bw(t))$ will appear after smoothing out two intersection points of the given rational curves, and hence the pairs of rational curves would sweep a locus of codimension $2$.

Now suppose that $\widehat\bw_\partial(t^*)$ consists of $n$ distinct points, $\widehat C=\widehat C_1\cup\widehat C_2$, where $\widehat C_1\simeq\PP^1$ and $\widehat C_2$ is a smooth elliptic curve, and the intersection
$\widehat C_1\cap\widehat C_2$ consists of one point $p$. Then, for the dimension reason, the sequence $\widehat\bw_\partial(t^*)$ is in general position subject to exactly two Menelaus type relations, while $w_0(t^*)$ is in general position in ${\mathcal Q}$. Thus, by Lemma \ref{l3-1a}, we obtain the immersion and smoothness statements as required. The rest of the argument literally coincides with the corresponding part of the proof of Lemma \ref{al4}(2iii).

\smallskip(2v) Suppose that $\widehat C$ is a rational curve with a node $p$. Then,
this node cannot be mapped to $\Tor(P_\Delta)^\times$, since otherwise
one would encounter a real rational curve with two one-dimensional branches, one in $\Tor_\R^+(P_\Delta)$ and the other in $\overline {\mathcal Q}$. Hence the node is mapped to $w\in\widehat\bw_\partial(t^*)$,
and for the above reason, the intersection multiplicities of the
two branches of $C=\bn(\widehat C)$ at $w$ are odd. Observe that $\bn':\PP^1\to\widehat C\to\Tor(P_\Delta)$ is not a multiple covering of the image. Indeed,
in such a case, one would have that $P$ is a triangle, $l_1=l_2=l$,
and $C=\bn(\widehat C)$ is a rational curve intersecting at least two of the toric divisors with odd multiplicity,
while the covering is double
ramified at the two points of $\widehat\bw_\partial(t^*)\setminus\{w\}$. With an automorphism of $\Z^2$,
we can turn the Newton triangle of $C$ (which is $\frac{1}{2}P_\Delta$) into the triangle
$$\conv\{(0,s),(r,0),(r+2l+1,0)\},\quad \gcd(s,r)\equiv1\mod2\ ,$$ with $w\in\Tor([(r,0),(r+2l+1,0)])$ and $C$ given by a parametrization (in some affine coordinates $x,y$)
$$x=a\theta^s,\quad y=b\theta^r(\theta-1)^{2l+1},\quad a,b\in\R^\times\ .$$ Since $x(1)>0$, we get $a>0$, and hence $x>0$ as $0<\theta\ll1$. Then $y>0$ for $0<\theta\ll1$, and hence $b<0$. Thus, if $r$ and $s$ are odd, we obtain $x>0,y<0$ as $\theta\to+\infty$, and $x<0,y<0$ as $\theta\to-\infty$.
If $r$ is even and $s$ is odd, we get $x>0,y<0$ as $\theta\to+\infty$, and $x<0,y>0$ as $\theta\to-\infty$. In both cases,
it follows that
the intersection point with the toric divisor $\Tor([(0,s),(r+2l+1,0)])$ is out of $\partial\Tor_\R^+(P_\Delta)$ against the initial assumption. If $r$ is odd and $s$ is even, we get that the segment of $\R C$ in the quadrant $\overline{\{x>0,y<0\}}$ touches all three toric divisors, which contradicts the condition (AQC).
Thus, a multiple covering is excluded.
The rest of claims can be derived in the same way as
the statement of Lemma \ref{al4}(2ii) with one exception: we
show that the local branches of $C$ at $w$ intersect each other with multiplicity
$\min\{2l_1+1,2l_2+1\}$. For $l_1\ne l_2$, this is immediate. Suppose that $l_1=l_2=l$. If $n=3$, then $P$ is a triangle and, in the above setting, we get
the parametrization of $C$ in the form
$$x=a\theta^{2s},\quad y=b\theta^{2r}(\theta-1)^{2l+1}(\theta-\theta_0)^{2l+1},\quad\theta_0\ne1,\quad a,b\in\R^\times\ .$$ Since $x(1)=x(\theta_0)$, we obtain the
only real solution $\theta_0=-1$, and hence $y=b\theta^{2r}(\theta^2-1)^{2l-1}$ is a function of $\theta^2$ as well as $x$ is, which means that
this is a double covering, forbidden above.

\smallskip(2vi)
Suppose now that $\widehat C$ is a smooth elliptic curve and
$\bn:\widehat C\to\Tor(P_\Delta)$ is a multiple covering of its image.
Since the preimages of at least $n-2$ points of $\widehat\bw_\partial(t^*)$ in $\widehat C$ are irreducible, the map $\bn:\widehat C\to C\subset\Tor(P_\Delta)$ is ramified, and hence $C$ cannot be elliptic by the Riemann-Hurwitz formula.
Thus, $C$ is rational, which for dimension reason implies that $\widehat\bw_\partial(t^*)$ consists of $n$ points
in general position subject to one Menelaus condition. In particular, $C$ is smooth at each point of $\widehat\bw_\partial(t^*)$,
and $\bn$ has ramification index $s$ at all points of $\widehat\bw_\partial(t^*)$. By the Riemann-Hurwitz formula,
$$0\le2s-(s-1)n.$$
Hence, $n\le\frac{2s}{s-1}$,which in view of the assumption $n\ge4$, yields $s=2$ and $n=4$ as stated in item (vi).
\proofend

\subsection{Wall-crossing}\label{crossing1}
We complete the proof of Theorem \ref{at2} with the following lemma.

\begin{lemma}\label{al9}
Let $\{(\widehat\bw_\partial(t),w_0(t))\}_{0\le t\le1}$ be a generic path in $\Men_\R^\rho\times {\mathcal Q}$, and let $t^*\in(0,1)$ be such that $\overline{\mathcal M}_{1,n}^\R(\Delta,\widehat\bw(t^*))$ contains an element
$\xi$ as described in one of the items of Lemma \ref{al5}(2).
Then,
for each $\kappa \in \frac{1}{2}\Z$ such that $|\kappa| \leq \cA(\Delta)$,
the numbers $W_1^\kappa(t):=W_1^\kappa(\Delta,\widehat\bw(t))$ and $\widetilde W_1^\kappa(t):=\widetilde W_1^\kappa(\Delta,\widehat\bw(t))$
do not change as $t$ varies in a neighborhood of $t^*$.
\end{lemma}

{\bf Proof.}
We always can assume that, in a neighborhood of $t^*$, the path $\{(\widehat\bw_\partial(t),w_0(t))\}_{0\le t\le 1}$ is defined by fixing the position of $w_0(t^*)$ and some $n-2$ points of $\widehat\bw_\partial(t^*)$, while the other two points remain mobile. Except for the case of Lemma \ref{al5}(2iv) describing reducible degenerations, we work with families of curves which are trivially covered by families of complex oriented curves so that the quantum index persists along each component of the family of oriented curves.

\smallskip
{\bf(1)} Suppose that $\xi$ is as in Lemma \ref{al5}(2i).
Consider the lifts of $w_0(t^*)$ upon $\widehat C$ that correspond to real branches of $C$ centered at $w_0(t^*)$. Since $C$ is immersed, we only have to show that each of the lifts induces a germ of a smooth one-dimensional subfamily in ${\mathcal M}_{1,n}(\Tor(P_\Delta),{\mathcal L}_{P_\Delta})$. This follows from the smoothness statements
in \cite[Proposition 4.17(2)]{DH} and in \cite[Lemma 3(1)]{Sh2}, \cite[Lemma 2.4(2)]{IKS4}.

\smallskip
{\bf(2)} Suppose that $\xi$ is as in Lemma \ref{al5}(2ii). We follow the argument used in Steps (1) and (2) of the proof of Lemma \ref{al8}. Since the local treatment of the occurring bifurcations is the same as in the proof of Lemma \ref{al8}, we derive the constancy of $W_1^\kappa(t)$ and $\widetilde W_1^\kappa(t)$ along the interval $(t^*-\eps,t^*+\eps)$ from the fact that the image of $H^0(C,\ina*{\mathcal L}_{P_\Delta})$ in $\prod_{\sigma,i}\Vers(C,w_i^\sigma)\times\prod_{z\in\Sing(C)\cap\Tor(P_\Delta)^\times}\Vers(C,z)$ transversally intersects with the tangent cone to the locus (\ref{elocus}), where $w_i^\sigma$ is a center of a local branch of $C$ of type $A_{2m}$, $m\ge k_i^\sigma$. This amounts to the relation of the form (\ref{ae9bb}) which we rewrite as
\begin{equation}h^1(\widehat C,{\mathcal O}_{\widehat C}(\bd_1))=0,\label{ae10}\end{equation} where the divisor $\bd_1$ satisfies
$$\deg\bd_1=C^2-(C^2+CK_{\Tor(P_\Delta)})-(-CK_{\Tor(P_\Delta)}-n+1)-(\mt(C,w_i^\sigma)-1)=n-2>0,$$
and $\mt(C,w_i^\sigma)$ stands for the multiplicity of $C$ at $w_i^\sigma$.
Thus, the required transversality follows.

\smallskip
{\bf(3)} The treatment of the case where $\xi$ is as in Lemma \ref{al5}(2iii)
is the same as that in step (3) of the proof of Lemma \ref{al8}. The required transvesality condition reads now as (cf. (\ref{etrans1}))
$$h^1(\widehat C,{\mathcal O}_{\widehat C}(\bd))=0,$$
where
$$\deg\bd=C^2-((C^2+CK_{\Tor(P_\Delta)})-(-CK_{\Tor(P_\Delta)}-2)=2>0,$$
which immediately yields the above $h^1$-vanishing.

\smallskip
{\bf(4)} Suppose that $\xi$ is as in Lemma \ref{al5}(2iv) and $\bn(p)\in\Tor(P_\Delta)^\times$.
The real part of the elliptic curve $C_2=\bn(\widehat C_2)$ has two one-dimensional components,
one in $\Tor_\R^+(P_\Delta)$ and the other in $\overline {\mathcal Q}$. The real part of the rational curve $C_1$ lies either in $\Tor_\R^+(P_\Delta)$, or in $\overline {\mathcal Q}$. Again the treatment of the wall-crossing is very similar to that in step (4) of the proof of Lemma \ref{al8}. We only comment on the transversality conditions ensuring the realizability of this scenario of the wall-crossing ina suitable variation of the constraints. Thus, we have a combination of two $h^1$-vanishinhg relations: (i) one is for $\widehat C_1$, and it coincides with (\ref{ae11}), (ii) the other reads (cf. (\ref{ae11}))
$$h^1(\widehat C_2,{\mathcal O}_{C_2}(\bd_2))=0,$$
where
$$\deg\bd_2=C_2^2-(C_2^2+CK_{\Tor(P_\Delta)})-(-CK_{\Tor(P_\Delta)}-1)=1>0,$$
which yields the desired $h^1$-vanishing.

Suppose that $\xi$ is as in Lemma \ref{al5}(2iv) and $\bn(p)=w_i^\sigma\in\widehat\bw(t^*)$ (to relax notation, below we write $w$ instead of $w^\sigma_i$).
Again, the local description of the corresponding wall-crossing event can be borrowed from part (5) of the proof of Lemma \ref{al8} and correspondingly, the signs $W_1^\kappa(\widehat\bw(t))$ and $\widetilde W_1^\kappa(\widehat\bw(t))$ do not change as $t$ transints through $t^*$. Thus, we should confirm the required transversality conditions.

If $\bn:\widehat C_1\to C_1$ and $\bn:\widehat C_2\to C_2$ both are immersions, then we need the relation (\ref{ehorace2}) which reduces to two $h^1$-vanishong conditions (\ref{ecert}) on $\widehat C_1\simeq\PP^1$ and $\widehat C_2$, a smooth elliptic curve. For $\widehat C_1$, that $h^1$-vanishing holds due to (\ref{ehorace3}). For $\widehat C_2$, it reads
\begin{equation}h^1(\widehat C_2,{\mathcal O}_{C_2}(\bd_2))=0,\label{ehorace4}\end{equation}
where the divisor $\bd_2$ satisfies
$$\deg\bd_2=C_2^2-(C_2^2+C_2K_{\Tor(P_\Delta)})-(-C_2K_{\Tor(P_\Delta)}-|C_2\cap\Tor(\partial P_\Delta|+1)$$
$$=|C_2\cap\Tor(\partial P_\Delta|-1\ge2>0,$$ and hence (\ref{ehorace4}) holds too.

In the case of $\bn:\widehat C_1\to C_1$ a $2l_1$-multiple ramified covering and $\bn:\widehat C_2\to C_2$ an immersion, the required transversality condition takes the form of (\ref{ehorace}). Using the Horace method in the same manner as in part (5b) of the proof of Lemma \ref{al8}, we come to the system of $h^1$-vanishing relations (\ref{ehorace1}). The upper line relations in (\ref{ehorace1}) are the same in the current setting, and hence hold for the same reason, while the second line relation turns into
$$h^1(\widehat C_2,{\mathcal O}_{\widehat C_2}(\bd))=0,$$
where
$$\deg\bd=C_2^2-(C_2^2+C_2K_{\Tor(P_\Delta)})-(-C_2K_{\Tor(P_\Delta)}-1)=1>0,$$
and hence it holds too.

\begin{figure}
\setlength{\unitlength}{1cm}
\begin{picture}(12,9)(-3,-0.5)
\includegraphics[width=0.6\textwidth, angle=0]{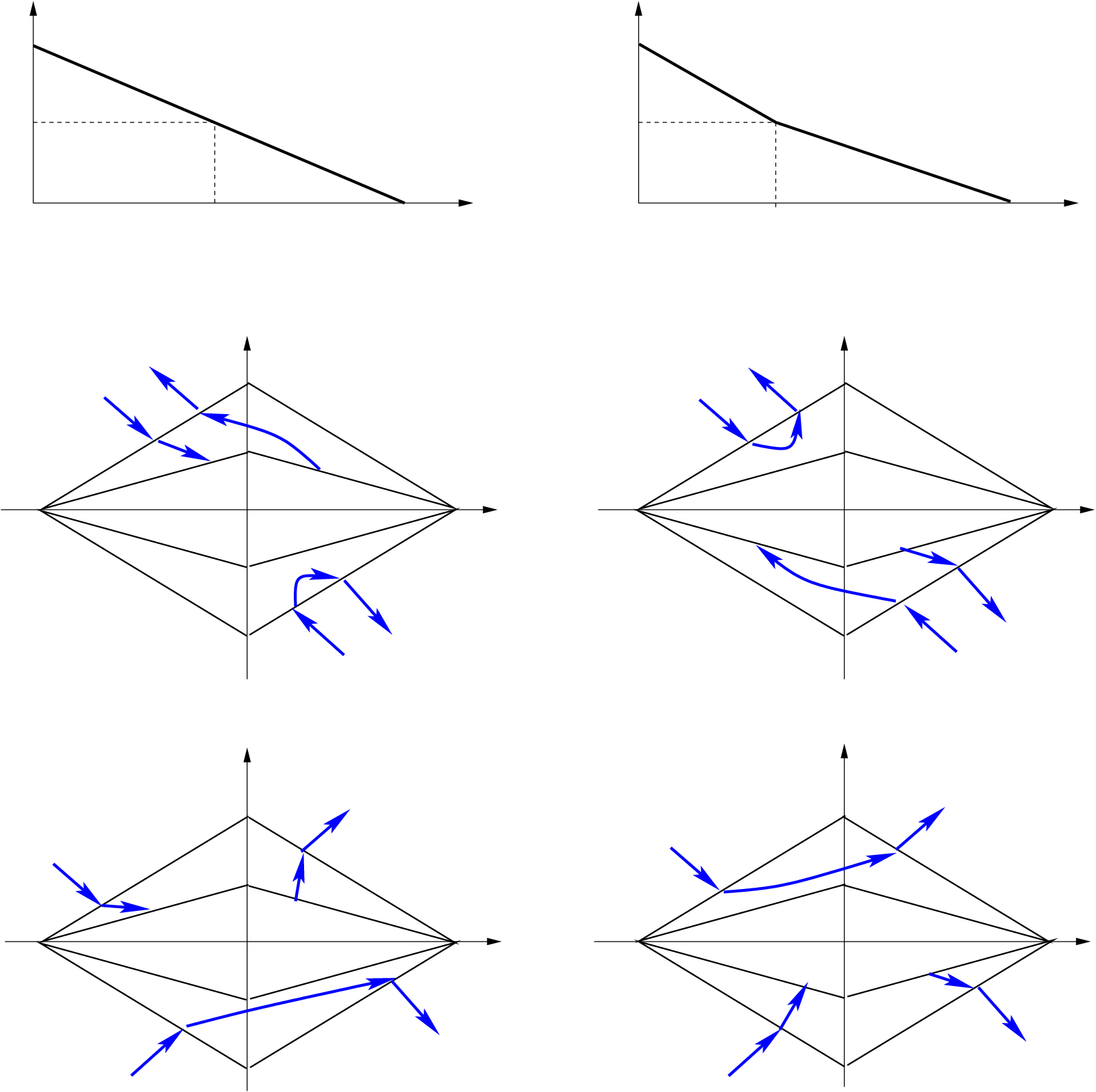}
\put(-8.9,6.6){(a)}\put(-3.9,6.6){(b)}
\put(-8.6,3.1){(c)}\put(-8.6,-0.3){(d)}\put(-4.7,4.7){$\Rightarrow$}\put(-4.7,1.1){$\Rightarrow$}
\put(-7.8,6.9){$2l+1$}\put(-6.3,6.9){$4l+2$}
\put(-6.3,5.8){$t<t^*$}\put(-1.3,5.8){$t>t^*$}\put(-6.3,2.5){$t<t^*$}\put(-1.3,2.5){$t>t^*$}
\put(-9,7.9){$1$}\put(-9.05,8.5){$2$}\put(-4,7.9){$1$}\put(-4.05,8.5){$2$}
\put(-3.2,6.9){$2l_1+1$}\put(-1.6,6.9){$2l_1+2l_2+2$}
\put(-7.4,8.4){$l_1=l_2=l$}\put(-2.4,8.4){$l_1<l_2$}
\end{picture}
\caption{Proof of Lemma \ref{al9}, part (5)} \label{af8}
\end{figure}

\smallskip
{\bf(5)} Suppose that $\xi$ is as in Lemma \ref{al5}(2v). Since $\widehat C$ is rational, both local branches $B_1,B_2$ of $C=\bn(\widehat C)$ at $w=w_i^\sigma(t^*)$ are real and they intersect $\Tor(\sigma)$ with odd multiplicities $2l_1+1,2l_2+1$, respectively, where
$2l_1+2l_2+2=2k_i^\sigma$, $l_1\le l_2$.
Choose two mobile points $w_{i'}^{\sigma'},w_{i''}^{\sigma''}\in\widehat\bw_\partial(t^*)\setminus\{w\}$ and fix the other points of $\widehat\bw(t^*)$. The one-dimensional stratum $\{\overline{\mathcal M}_{1,n}^{\R,+}(\Delta,
\bw(t))\}_{|t-t^*|<\zeta}$, $0<\zeta\ll1$, projects to the germ at $C$ of a one-dimensional variety in $|{\mathcal L}_{P_\Delta}|$,
and we denote by ${\mathcal V}$
one of its real irreducible components.
Let the tangent line to ${\mathcal V}$ at $C$ be spanned by $C$ and $C^*\in|{\mathcal L}_{P_\Delta}|\setminus\{C\}$.
Fix conjugation-invariant coordinates $x,y$ in a neighborhood of $w$ so that $w=(0,0)$ and $\Tor(\sigma)=\{y=0\}$.
Then, in a neighborhood of $w$, the curve $C$ is given by an equation
$$(y+\eta_1x^{2l_1+1})(y+\eta_2x^{2l_2+1})+\text{h.o.t.}=0,\quad\eta_1,\eta_2\in\R^\times\ ,$$ with the Newton diagram (see Figure \ref{af8}(a,b))
$$\Gamma=[(0,2),(2l_1+1,1)]\cup[(2l+1,1),(2l_1+2l+2+2,0)]$$ (cf. Figures \ref{af3}(a,b)). Exactly as in step (5a) of the proof of Lemma \ref{al8}, an equation of $C^*$ has the monomial $x^{2l_1}y$ with a non-zero coefficient and no other monomials below the Newton diagram $\Gamma$.

If $l_1=l_2=l$, then as in step (5a) of the proof of Lemma \ref{al8}, we obtain that the equation of the germ of $C^{(t)}\in{\mathcal V}$ at
$w=(0,0)$ is as follows:
\begin{equation}y^2+P(x)y+\eta_1\eta_2x^{4l+2}+\text{h.o.t.}=0,\quad P(x)=(\eta_1+\eta_2)x^{2l+1}+\sum_{j=2l}^0a_j\tau^{2l+1-j}x^j\ ,
\label{ae20}\end{equation}
$$\text{where}\quad \tau=t-t^*, \quad\text{and}\quad x^{2l+1}P(1/x)=\lambda_1\Cheb_{2l+1}(\tau x+\lambda_2)+\lambda_3$$
with real numbers $\lambda_1,\lambda_2,\lambda_3$
uniquely determined by $\eta_1,\eta_2$ and by the coefficient of $x^{2l}y$ in the equation of $C^*$,
and h.o.t. includes higher powers of $\tau$ in the present monomials and the sum of monomials above $\Gamma$.
If $\eta_1\eta_2>0$, the latter formula yields $2l$ elliptic nodes in a neighborhood of $w$, all of them lying either in the domain $\{y>0\}$, or in the domain $\{y<0\}$. Observe that the change of sign of $\tau$
yields the symmetry with respect to the origin for the real part of $C^{(t)}$ (up to h.o.t.), while the orientation of the local branch of $C^{(t)}$ at $w$ persists, and hence 
the quantum index and the signs $W_1^\kappa(t)$ and $\widetilde W_1^\kappa(t)$ remain constant along the path
$\{(\widehat\bw(t))\}_{|t-t^*|<\zeta}$.
If $\eta_1\eta_2<0$, then equation (\ref{ae20}) implies
$$y=\frac{1}{2}\left(-P(x)\pm\sqrt{P(x)^2-4\eta_1\eta_2x^{4l+2}}\right)+\text{h.o.t.}\ ,$$ which means that $C^{(t)}$ has only non-real nodes
in a neighborhood of $w$. Again, the orientation of the local branch of $C^{(t)}$ at $w$ persists (see Figure \ref{af8}(d)), and we conclude that the signs $W_1^\kappa(t)$ and $\widetilde W_1^\kappa(t)$ do not change as $t$ transits through $t^*$. 

Basically, the same local bifurcation occurs when If $l_1<l_2$ (cf. step (5a) of the proof of Lemma \ref{al8}). The family ${\mathcal V}$ appears to be smooth, regularly parameterized by $t_1=t-t^*$ in the form
\begin{equation}y^2+\eta_1(x-\eta_0t_1)^{2l_1+1}y+\eta_1\eta_2x^{2l_1+2l_2+2}+\text{h.o.t.}=0\ ,
\label{ae22}\end{equation} where $\eta_0$ is determined by $C^*$, and h.o.t. includes higher powers of $t_1$ in the present monomials and the sum of monomials above $\Gamma$.
We skip the rest of details.

\begin{remark}\label{ar3}
In equations (\ref{ae20}) and (\ref{ae22}), the sign of the coefficient of $x^{2l_1+2l_2+2}$ remains constant,
and the sign of the coefficient of $y$ changes
as $t_1$ changes its sign. Geometrically, this means that, for $t_1<0$,
the real branch of $C^{(t)}$, which is tangent to $\Tor(\sigma)$ at $w$,
lies in
$\Tor_\R^+(P_\Delta)$, while, for $t_1>0$, this real branch lies in $\overline {\mathcal Q}$, or vice versa (see Figure \ref{af8}(c,d)).
\end{remark}

At last, we establish the transversality condition sufficient for the realizability of the wall-crossing described above by a variation of constraints. Namely, we
prove that representing the local deformation as above, and crossing the wall transversally,
we need the following transversality condition:
the Zariski tangent space to the germ at $C^*$ of the family constrained by the given tangency conditions at fixed points of $\widehat\bw_\partial(t^*)$ and at two mobile points of $\widehat\bw_\partial(t^*)$ as well as by fixing the point $w_0$ has dimension zero. The routine computations, performed in preceding steps, reduce this condition to
\begin{equation}
h^0\left(\PP^1,{\mathcal N}_{\bn'}\left(-\sum_{p_i^\sigma\ne p}2k_i^\sigma p_i^\sigma+p_{i'}^{\sigma'}+p_{i''}^{\sigma''}-(2l_1+1)p'-(2l_2+1)p''-p_0\right)\right)=0,\label{tce1}\end{equation}
where $\bn':\PP^1\to\widehat C\overset{\bn}{\to}\Tor(P_\Delta)$ is $\bn$ combined with the normalization, $p',p''$ are the preimages of $p\in\widehat C$. This $h^1$-vanishing follows from the fact that the degree of the line bundle is
$$(-CK_{\Tor(P_\Delta)}-2)+(CK_{\Tor(P_\Delta)}-2k_i^\sigma+2)-((2l_1+1)+(2l_2+1)+1)=-1.$$

\smallskip
{\bf(6)} Suppose that $\xi$ is as in Lemma \ref{al5}(2vi). Since the rational curve $C=\bn(\widehat C)$ has one global one-dimensional real branch and contains the point $w_0\in Q$, it intersects the toric divisor $\Tor(\sigma)$ at two points $w_1,w_2$ with odd multiplicities $k_1,k_2$, respectively, where $\Tor^+_\R(\sigma)$
is the (unique) common side of $\Tor^+_\R(P_\Delta)$ and $\overline Q$.
Correspondingly, $P_\Delta$ is a lattice triangle, and each of
the two other toric divisors intersects $C$ with an even multiplicity,
see Figure \ref{fhorace}(a) (the two fragments of $\R C$ shown in bold are doubly covered by the two ovals of $\R\widehat C$).

\begin{figure}
\setlength{\unitlength}{1cm}
\begin{picture}(14,6)(0,0)
\thinlines
\put(2.5,1){\vector(0,1){5}}\put(0,3.5){\vector(1,0){5}}
\put(7,1){\vector(0,1){5}}\put(7,1){\vector(1,0){7}}

\thicklines
\put(0.5,3.5){\line(1,1){2}}\put(0.5,3.5){\line(1,-1){2}}
\put(0.5,3.5){\line(1,0){4}}\put(2.5,1){\line(0,1){4}}
\put(2.5,5.5){\line(1,-1){2}}\put(2.5,1.5){\line(1,1){2}}

\put(7,3){\line(0,1){2}}\put(7,3){\line(3,-1){6}}\put(7,5){\line(3,-2){6}}

\qbezier(3.5,4.5)(4.5,3.5)(3,3.5)
\qbezier(2.2,4.2)(2.2,3.5)(3,3.5)
\qbezier(2.2,4.2)(2.2,5.8)(3.5,4.5)

\put(1.6,4){$Q$}\put(6.7,2.9){$1$}\put(6.7,4.9){$2$}\put(12.7,0.5){$2k_m$}\put(8.4,3){$T_m$}
\put(2.3,0){(a)}\put(9.3,0){(b)}

\end{picture}
\caption{Proof of Lemma \ref{al8}, part (6)}\label{fhorace}
\end{figure}

We define the path $\{(\widehat\bw_\partial(t),w_0(t))\}_{|t-t^*|<\eps}$ by fixing the position of $w_1$, $w_2$, $w_3$, and $w_4$ and leaving $w_0$ mobile. Note that $w_0$ can be chosen in general position on $C$. We suppose that $\{w_0(t)\}_{|t-t^*|<\eps}$ sweeps the germ of a real line $L\subset Q$ transversally intersecting $\R C$ at $w_0=w_0(t^*)$.

Note also that $(C, (w_0,w_1,w_2,w_3,w_4))$ has two preimages: two
marked curves $(\widehat C,(p'_0,p_1,p_2,p_3,p_4))$ and $(\widehat C,(p''_0,p_1,p_2,p_3,p_4))$. It follows from (\ref{eis1}) that the quantum index of the two latter curves equipped with any complex orientation orientation is zero. The signs (\ref{ae9i}) and (\ref{ae21}) are preserved too when we reverse the orientation.

The family $\{\overline{\mathcal M}_{1,n}^{\R,+}(\Delta,\widehat\bw(t))\}_{|t-t^*|<\eps}$ naturally projects to a bunch of real one-dimensional subvarieties in the linear system $|{\mathcal L}_{P_\Delta}|$. Denote by ${\mathcal V}$ the union of those components of the image that pass through $2C$ and have the same tangent line (${\mathcal V}$ is determined up to a finite choice). Let this common tangent line $\Lambda$ be spanned by $2C$ and some $C^*\in|{\mathcal L}_{P_\Delta}|\setminus\{2C\}$. Two cases are possible: either $C^*$ contains $C$ as a component, or not.

Let $C\not\subset C^*$. Introduce local conjugation-invariant coordinates $x,y$ in a neighborhood of $w_0$ so that $$w_0=(0,0),\ L=\{y=0\},\ C=\{x=0\},\ C^*=\left\{a+\sum_{i+j>0}a_{ij}x^iy^j=0\right\},$$
where $a\ne0$. Let us parameterize $\Lambda$ as follows:
$$C(\theta)=\left\{x^2+\theta\left(a+\sum_{i+j>0}a_{ij}x^iy^j\right)=0\right\},\quad \theta\in(\R,0).$$
Then, for $a\theta<0$, the curve $C(\theta)$ intersects $L$ at two points, $x=\sqrt{-a\theta}+o(\sqrt{\theta})$ and $x=-\sqrt{-a\theta}+o(\sqrt{\theta})$. The same holds for any curve $C'\in{\mathcal V}$ projecting to $C(\theta)$; hence, $C'$ contributes to the count of $W_1^0(t)$ and $\widetilde W_1^0(t)$ both, for $t<t^*$ and for $t>t^*$. This yields the constancy of $W_1^0(t)$ and $\widetilde W_1^0(t)$ over the segment $t\in(t^*-\eps,t^*+\eps)$.

Let $C^*=C\cup C'$, $C'\ne C$. In the above coordinates $x,y$, we have
$$C^*=\left\{x\left(a+\sum_{i+j>0}a_{ij}x^iy^j\right)=0\right\},\quad a\ne0.$$
Thus,
$$\Lambda=\left\{C(\theta):=\left\{x\left(x+\theta\left(a+\sum_{i+j>0}a_{ij}x^iy^j\right)\right)=0
\right\}\right\}_{\theta\in\R,0)},$$
\begin{equation}C(\theta)\cap L=\{w_0,w_0(t)\},\quad w_0(t)=(-a\theta+O(\theta^2),0).\label{esecond}\end{equation}
Without loss of generality, we can assume that ${\mathcal V}$ consists of one component. If ${\mathcal V}$ homeomorphically projects onto the germ $(\Lambda,2C)$, and hence onto the segment $(t^*-\eps,t^*+\eps)$, we establish the constancy of $W_1^0(t)$ and $\widetilde W_1^0(t)$ by showing that a curve $C^{(t)}$ makes the same contribution to $W_1^0(t)$ both for $t<t^*$ and for $t>t^*$, and so for $\widetilde W_1^0(t)$.

To this end, we analyze the singularities of  curve $C^{(t)}$, $0<|t-t^*|<\eps$. The curve $C$ has $(C^2+CK_{\Tor(P_\Delta)}+2)/2$ nodes. In a neighborhood of each of these nodes, the curve $C^{9t)}$ has four nodes:
four hyperbolic nodes in a neighborhood of a hyperbolic node of $C$, two elliptic and two complex conjugate nodes in a neighborhood of an elliptic node of $C$, four non-real nodes in a neighborhood of a non-real node of $C$.
 There remain
$$\frac{4C^2+2CK_{\Tor(P_\Delta)}}{2}-4\cdot\frac{C^2+CK_{\Tor(P_\Delta)}+2}{2}=k_1+k_2+k_3+k_4-4$$
nodes of $C^{(t)}$. We
show that they are located in a neighborhood of the points $w_1,w_2,w_3,w_4$.
Consider the tropical limit of the family of the parameterized curves $\{\bn_t:\widehat C^{(t)}\to C^{(t)}\hookrightarrow{\mathfrak X}_t=\Tor(P_\Delta)\}_{C^{(t)}\in{\mathcal V},0<|t-t^*|<\eps}$ as $t\to t^*$
(see details in \cite[Section 1]{GaSh}, where the setting is most suitable for our argument).
The limit of the maps $\bn_t:\widehat C^{(t)}\to {\mathfrak X}_t=\Tor(P_\Delta)$, $t\to t^*$, is $\bn_0:\widehat C^{(0)}\to {\mathfrak X}_0$, where $\widehat C^{(0)}$ is the union of the elliptic curve $\widehat C$ with four marked points $p_1,p_2,p_3,p_4$ mapped to $w_1,w_2,w_3,w_4$, respectively, and of connected curves $\widehat C_1,\widehat C_2,\widehat C_3,\widehat C_4$ of arithmetic genus zero attached to $p_1,p_2,p_3,p_4$, respectively. Pick $m\in\{1,2,3,4\}$ and introduce local conjugation-invariant
coordinates $x_m,y_m$ in a neighborhood of $w_m$ so that $w_m=(0,0)$, the toric divisor passing through $w_m$ is given by $\{y=0\}$, and the germ of $C$ at $w_m$ is given by $\{y-x^{k_m}=0\}$. In particular, a local equation of $2C$ reads
$$(y_m-x^{k_m})^2=0.$$ Components of $\widehat C_m$ are mapped to toric components of ${\mathfrak X}_0$ that correspond to a certain subdivision of the region
under the segment $[(0,2),(2k_m,0)]$. In view of the given data, the only possible subdivision contains just one triangle $T_m=\conv\{(0,1),(0,2),(2k_m,0)\}$ (see Figure \ref{fhorace}(b)). The limit curve $C_m\subset\Tor(T_m)$ is rational, quadratically tangent to the toric divisor $\Tor([(0,2),(2k_m,0)])$, and it has $k_m-1$ nodes (cf. \cite[Lemma 3.5]{Sh0}). An equation of $C_m$ reads
$$(y_m-x_m^{k_m})^2+\sum_{i=0}^{k_m-1}a_i^{(m)}x_m^iy_m=0,\quad a_0^{(m)}\ne0.$$
Correspondingly, the curves $C^{(t)}\in{\mathcal V}$ are given in a neighborhood of $w_m$ by
\begin{equation}(y_m-x_m^{k_m})^2+\sum_{i=0}^{k_m-1}\tau^{(k_m-i)r_m}a_i^{(m)}x_m^iy_m 
+\text{h.o.t.}=0,\label{esecond1}\end{equation} where $\tau=t-t^*$, and $\text{h.o.t.}$ denotes terms with higher powers of $\tau$ in the preceding monomials as well as monomials above the segment $[(0,2),(2k_m,0)]$.

If $m=1$ or $2$, then $k_m$ is odd. In this case, all exponents of $\tau$ in the second summand of (\ref{esecond1}) must be even since otherwise the coefficient of $y_m$ will change sign when $\tau$ does, but geometrically this would mean that the local real branch of $C^{(t)}$ at $w_m$ would jump from the quadrant $Q$ to the positive quadrant when $t$ transits through $t^*$ which would contradict the fact that $C^{(t)}$ must have a global real branch in $Q$. This observation yields that the distribution among the quadrants of the elliptic nodes in a neighborhood of $w_m$ does not change when $\tau$ changes its sign. If $m=3$ or $4$,
then $k_m$ is even.
If all exponents of $\tau$ in the second summand of (\ref{esecond1}) are even,
we obtain the preceding conclusion. If there are odd exponens,
then the change of sign of $\tau$ geometrically is equivalent to the reflection of the limit curve $C_m$
with respect to the $y_m$-axis, which again does not affect the distribution of elliptic nodes among quadrants.
In all the cases, we immediately derive the constancy of $W_1^0(t)$ and $\widetilde W_1^0(t)$
along the segment $t\in(t^*-\eps,t^*+\eps)$.

It remains to show that the case of ${\mathcal V}$ doubly covering a half of the interval $(t^*-\eps,t^*+\eps)$
with ramification at $t^*$ is not possible. We argue by contradiction based on the B\'ezout theorem.
Let $C^{(t')},C^{(t'')}\in{\mathcal V}$ be two curves with $t'<t^*<t''$ that intersect the line $L$
at the same point $w_0(t')=w_0(t'')$ given by the second formula in (\ref{esecond}).
In a neighborhood of $\Sing(C)$,
the curves $C^{(t')}$ and $C^{(t'')}$ have at least
$$8\cdot\frac{C^2+CK_{\Tor(P_\Delta)}+2}{2}=4C^2+4CK_{\Tor(P_\Delta)}+8=4C^2-\sum_{m=1}^4(4k_m-2)$$
intersection points. To estimate the intersection multiplicity of $C^{(t')}$ and $C^{(t'')}$ in a neighborhood of $w_m$, $m=1,2,3,4$, we rescale the formula (\ref{esecond1}) bringing it to the common for $t'$ and $t''$ form (since $|t'-t^*|$ and $|t''-t^*|$ are asymptotically equivalent as $t',t''\to t^*$)
$$(y_m-x_m^{k_m})^2+\sum_{i=0}^{k_m-1}(\sign\tau)^{(k_m-i)r_m}a_i^{(m)}x_m^iy_m+O(\tau)=0.$$
In a neighborhood of $w_m$, $m=1,2,3,4$, the curves $C^{(t')},C^{(t'')}$ intersect with multiplicity at least
$$4k_m-2.$$ Indeed, if the change of sign of $\tau$ in (\ref{esecond1}) yields the same limit curve $C_m$, then $C^{(t')}$ and $C^{(t'')}$ intersect with multiplicity at least $2k_m$ at $w_m$ and intersect with multiplicity at least $2(k_m-1)$ near the nodes of $C_m$. If the change of sign of $\tau$ induces the reflection of $C_m$ against the $y_m$-axis, then $C^{(t')}$ and $C^{(t'')}$ intersect in the big torus of $\Tor(T_m)$ with multiplicity $2k_m-2$ (here we subtract the intersection multiplicity at the common tangency point on the toric divisor $\Tor([(0,2),(2k_m,0])$) and they intersect at $w_m$ with multiplicity $\ge2k_m$. Taking into account the extra intersection point $w_0(t')=w_0(t'')\in L$, we arrive to a contradiction with the B\'ezout theorem.
\proofend

\section{Invariance in genus two: proof of Theorem \ref{g2t1}}\label{asec6}
We closely follow the argument in the proof of Theorems \ref{at1} and \ref{at2}.
Take two sequences $\bw(0)=(\bw_\partial(0),w_1(0),w_2(0))$ and $\bw(1)=(\bw_\partial(1),w_1(1),w_2(t))$,
satisfying the conditions of Theorem \ref{at2}. We may assume that
$\widehat\bw_\partial(0),\widehat\bw_\partial(1)\in \Men_\R^\rho(\Delta)$ for some
$\rho > 0$.
Then,
we join these sequences by a generic path $\{\widehat\bw_\partial(t)\}_{0\le t\le 1}$, in $\Men_\R^\rho(\Delta)$ and
the points $w_i(0)$, $w_i(1)$ by a generic path $\{w_i(t)\}_{0\le t\le1}$ in the respective open quadrant, $i=1,2$. Then we verify the constancy
of $W_2^\kappa(\Delta,\widehat\bw(t))$ inside the chambers ({\it i.e.}, connected components of $\Men_2(\Delta)$)
and in all possible wall-crossing events, for all $\kappa$.

\subsection{Invariance inside chambers}\label{chamb2}
Let us verify
the constancy of $W_2^\kappa(\Delta,\widehat\bw(t))$
along intervals $t'<t<t''$ such that $\bw(t)\in\Men_2(\Delta)$ for all $t'<t<t''$.
To this end, it is enough to show that the projection ${\mathcal M}_{2,n}^{\;\R}(\Delta,\widehat\bw(t))_{t'<t<t''}\to(t',t'')$ has no critical points (cf. the proof of Theorems \ref{at1} and \ref{at2}).
This requirement is a consequence of the following statement similar to Lemmas \ref{lis2} and \ref{lis3} (cf. \cite[Lemma 2.3]{Sh}).

\begin{lemma}\label{lis5}
Let $\xi=[\bn:(\widehat C,\bp)\to\Tor(P_\Delta)]\in{\mathcal M}_{2,n}^{\;\R}(\Delta,\widehat\bw)$,
where $\widehat\bw\in\Men_2(\Delta)$.
Then,
\begin{equation}h^0\left(\widehat C,{\mathcal N}_\bn\left(-p_1-p_2-\sum_{s=0}^2\sum_{i=1}^{n_s}2k_i^{(s)}p_i^{(s)}+p_{n_0}^{(0)}\right)\right)=0,
\label{g2e1}\end{equation}
where $[\bn:(\widehat C,\bp)\to\PP^2]\in{\mathcal M}_{2,n}^\R(\PP^2,\Delta,\widehat\bw(t))$.
\end{lemma}

{\bf Proof.}
The non-vanishing in (\ref{g2e1}) would mean the existence of a real plane curve $C'\ne C=\bn(\widehat C)$ of degree $2d$ passing through $w_1$ and $w_2$, intersecting $D_s$, $s=0,1,2$, at each point $n_i^{(s)}$, $1\le i\le n_s$, $(s,i)\ne(0,n_0)$, with multiplicity $2k_i^{(s)}$, intersecting $D_0$ at $w_{n_0}^{(0)}$ with multiplicity $2k_{n_0}^{(0)}-1$, and intersecting $C$ at each point $z\in\Sing(C)$ with multiplicity $\ge\delta(C,z)$. Note that $C'$ must additionally intersect $C$ in an extra point on the real branch in $(\PP_\R^2)^{(-,+)}$ and in an extra point on the real branch in $(\PP_\R^2)^{(+,-)}$. However, this leads to a contradiction with the B\'ezout theorem:
$$4d^2\ge 4+\sum_{s=0}^2\sum_{i=1}^{n_s}2k_i^{(s)}-1+\sum_{z\in\Sing(C)}2\delta(C,z)$$
$$=4+(6d-1)+((2d-1)(2d-2)-4)=4d^2+1.$$
\proofend

\subsection{Walls in the space of constraints}\label{walls2}
Let $[\bn:(\widehat C,\bp)\to\PP^2]\in\overline{\mathcal M}_{2,n}^\R(\Delta,\widehat\bw(t))_{0\le t\le1}$.
In this section, we assume that $\widehat C$
does not contain contractible components, {\it i.e.}, we just contract them (if any)
and we call the singularities of $\widehat C$ ``nodes".
Note that such a contraction does not affect the arithmetic genus, since a generic curve $\widehat C$ has three real ovals. We also allow the points of $\widehat\bp$ to collate. We
use the notations of Section \ref{walls1}: namely, $S$ denotes a locus in the space of constraints
$\Men^\rho_\R(\Delta)\times(\PP^2_\R)^{(+,-)}\times(\PP^2_\R)^{(-,+)}$ and $\pr(S)$ denotes the image of
$\overline{\mathcal M}_{2,n}^\R(\Delta,\widehat\bw)$ in the linear system $|{\mathcal L}_{P_\Delta}|$
when $\widehat\bw$ ranges over $S$.

In view of the results of Section \ref{chamb2}, we can ignore the walls along which
$\overline{\mathcal M}^\R_{2,n}(\Delta,\widehat\bw)={\mathcal M}^\R_{2,n}(\Delta,\widehat\bw)$.
The following lemma classifies other walls.

\begin{lemma}\label{g2l1} (1) The following elements $\xi=[\bn:(\widehat C,\bp)\to\PP^2]$ cannot occur in $\overline{\mathcal M}_{2,n}^\R(\Delta,\widehat\bw(t))_{0\le t\le1}$:
\begin{enumerate}\item[(1i)] $\widehat C$ is a reducible connected curve of arithmetic genus $2$ with a component mapped onto a toric divisor;
\item[(1ii)] $\widehat C$ is a connected curve of arithmetic genus $2$ with at least two nodes;
    \item[(1iii)] a component of $\widehat C$ multiply covers its image.
\end{enumerate}

(2) If $\xi=[\bn:(\widehat C,\widehat\bp)\to\PP^2]\in\overline{\mathcal M}_{2,n}(\Delta,\widehat\bw(t^*))$ when $\widehat\bw(t^*)$ does not satisfy the conclusions of Lemma \ref{l3-1a}, then $\xi$ is of one of the following types:

(2a) either exactly one of the points $w_1(t^*)$, $w_2(t^*)$ is a singular point of the curve $C=\bn(\widehat C)$, and the following holds:
\begin{enumerate}
\item[(2i)] $\widehat\bw_\partial(t^*)$ consists of $n$ distinct points, the curve $\widehat C$ is smooth of genus $2$, $\bn$ is an immersion onto a curve $C$ which is smooth along
the toric divisors, while $w_i(t^*)$ is a center of at least one real local branch, $i=1,2$;
\end{enumerate}

(2b) or $w_1(t^*)$, $w_2(t^*)$ are smooth points of the curve $C=\bn(\widehat C)$, and one of the following holds:
\begin{enumerate}
\item[(2ii)] $\widehat C$ is a smooth curve of genus $2$, $\widehat\bw_\partial(t^*)$ consists of $n-1$ distinct points, $\bn$ is birational onto its image, but not an immersion: it may have singular local branches centered in $\Tor(P_\Delta)^\times$ as well as singularities of type $A_{2m}$ at some points $w_i^\sigma\in\widehat\bw_\partial(t^*)$, where $m\ge k_i^\sigma$;
\item[(2iii)]
$\widehat C$ is a smooth curve of genus $2$, $\widehat\bw_\partial$ contains $n-2$ distinct points so that $w_i^{(s)}(t^*)=w_j^{(s)}(t^*)=:w$ for some $s\in\{0,1,2\}$ and $i\ne j$; the map $\bn$
is birational onto its image $C$, where the curve $C$ is immersed outside $\widehat\bw_\partial(t^*)\setminus\{w\}$ and may have singularities of type
$A_{2m}$ at some points of $\widehat\bw_\partial(t^*)\setminus\{w\}$, while at $w$, $C$ either is smooth or has two smooth branches;
\item[(2iv)] $\widehat C=\widehat C_1\cup\widehat C_2$, where $g(\widehat C_1)+g(\widehat C_2)=2$,  $\widehat\bw_\partial(t^*)$ consists of $n-1$ distinct points, and the maps $\bn:\widehat C_i\to\PP^2$, $i=1,2$, are immersions; furthermore, the only node $p$ joining $\widehat C_1$ and $\widehat C_2$ either is mapped to $(\R^\times)^2\subset\PP_\R^2$, or
    is mapped to some point $w_i^{(s)}(t^*)\in\widehat\bw(t^*)$;
the germ $\bn(\widehat C_j,p)$ intersects $D_s$ with multiplicity $2l_j$, $j=1,2$, where $l_1+l_2=k_i^{(s)}$ and these two germs intersect each other with multiplicity $\min\{2l_1,2l_2\}$;
\item[(2v)] $\widehat\bw_\partial(t^*)$ consists of $n-1$ distinct points, $\widehat C$ is an irreducible elliptic curve with one node $p$,
the map $\bn$ is an immersion that sends $p$ to some point $w_i^{(s)}(t^*)$, center of two smooth branches intersecting $D_s$ with odd multiplicities
    $2l_1+1,2l_2+1$, where $l_1+l_2+1=k_i^{(s)}$, and intersecting each other with multiplicity $\min\{2l_1+1,2l_2+1\}$; furthermore, each of the other points of $\widehat\bw_\partial(t^*)$ is a center of one smooth branch.
\end{enumerate}
\end{lemma}

{\bf Proof.}
(1i)
We can assume that the sequences $\{\bw(t)\}_{0\le t\le1}$ avoid $\rho_0$-neighborhoods
of the intersection points of the coordinate axes, for some
$\rho_0 > 0$.

Suppose that $\bn(\widehat C)$
contains some coordinate axis.
In this case, $\bn(\widehat C)$
contains all of them.
Indeed, otherwise, we would have $\bw(t)$ arbitrarily close to the intersection points of the coordinate axes, contrary to the above assumption. Furthermore, the intersection points of the coordinate axes must lift to nodes of $\widehat C$, and they smooth out along the path $\{\bw(t)\}_{0\le t\le1}$. In particular, the union $\widehat C_1$ of the components of $\widehat C$ mapped onto the coordinate axes is a curve of arithmetic genus $\ge1$.

The union $\widehat C_2$ of the other components of $\widehat C$ cannot be mapped onto a plane curve of odd degree. Indeed, otherwise, one would get an odd number of real branches of $\bn(\widehat C_2)$ intersecting $D_0$ with odd multiplicity yielding a one-dimensional real arc in the quadrant $\{x_1/x_0<0,\ x_2/x_0<0\}$ which is not possible. Thus, $\bn_*(\widehat C_2)$ has an even degree. If $\widehat C_2$ has an arithmetic genus $0$, then the real part of $\bn(\widehat C_2)$ contains a connected graph passing through $w_1$ and $w_2$, and then there should be four real local branches of $\bn(\widehat C_2)$ centered on
$D_1^+\cup D_2^+$ (see Section \ref{sec-g2} for notations)
that intersect $D_1\cup D_2$ with odd multiplicity. Hence, the centers of these branches must lift to the nodes of $\widehat C$, which would yield that the arithmetic genus of $\widehat C$ is at least $4$;
a contradiction. Thus, $\widehat C_1$ and $\widehat C_2$ both have arithmetic genus $1$, and are joined by one node.
It is easy to see, that this node must be mapped by $\bn$ to the center of the unique local branch of $\bn(\widehat C_2)$
on $D_0^+$. This local branch must lie in the closure of $(\PP_\R^2)^+$. Hence the image of one connected component
of the real part of $\widehat C_2$ must pass through, say, $w_1$ and contain the real local branch centered on $D_0^+$.
It follows that $\bn(\widehat C_2)$ has two real local branches centered on $D_1^+$
and intersecting $D_1$ with odd multiplicity. Thus, their centers must lift to two nodes of $\widehat C$
joining $\widehat C_1$ and $\widehat C_2$ making the arithmetic genus of $\widehat C$
at least $3$ which is a contradiction.

\smallskip (1ii) We argue on contradiction assuming that $\widehat C$ has $k\ge2$ nodes and no component of $\widehat C$ is mapped onto a coordinate axis.

First, $\widehat C$ cannot have $\ge3$ components, since it would impose two extra Menelaus conditions on the disposition of $\bw_\partial$ defining a stratum of codimension $2$ in ${\mathcal W}$ which must be avoided by the path $\{\bw(t)\}_{0\le t\le1}$.

Second, $\widehat C$ cannot be an irreducible rational curve. Indeed, if $\bw_\partial\in{\mathcal W}_\partial$ is a general position, then $\bw(\widehat C)$ belongs to a finite set, see Lemma \ref{l3-1a}. Hence, if $\bw_\partial\in{\mathcal W}_\partial$ varies in a stratum of codimension $1$, then we obtain a one-dimensional family of rational curves, and this imposes an extra condition on the position of the points $w_1,w_2$, altogether defining a stratum of codimension $\ge2$ in ${\mathcal W}$.

Third, $\widehat C$ cannot be the union of a rational and an elliptic component. Then,
we obtain an extra Menelaus condition on $\bw_\partial$, while the parts of $\bw_\partial$ lying on the image
of the rational and
the elliptic components are in general position.
Thus, $\bw_\partial$ and one of the points $w_1,w_2$ determine a finite set of pairs
of a rational and an elliptic components of $\widehat C$ (Lemma \ref{l3-1a}),
which imposes a condition on the location of the second point in $\{w_1,w_2\}$,
hence, defining a stratum of codimension $\ge2$ in ${\mathcal W}$.

\smallskip (1iii) Suppose that $\widehat C$ is irreducible
and $\bn:\widehat C\twoheadrightarrow C\hookrightarrow\PP^2$ is a multiple covering.
By the Riemann-Hurwitz formula,
the curve $C$ cannot be of geometric genus $2$. Also, $C$ cannot be rational since otherwise, the fixed position of $\bw_\partial$ determines a finite set of rational curves, and hence two conditions are imposed on the position of $w_1,w_2$. This defines a stratum of codimension $\ge2$ in ${\mathcal W}$ avoided by the path $\{\bw(t)\}_{0\le t\le1}$. The only remaining option is an elliptic curve $C$. As in the case of the rational curve $C$,
we derive that $\bw_\partial$ varies in an open dense subset of ${\mathcal W}_\partial$,
while the position of $w_1,w_2$ is constrained by one condition.
It follows that all the points of the sequence $\widehat\bw_\partial$ are distinct,
and hence they all must be ramification points of the covering.
However, $\#\widehat\bw_\partial\ge3$ which yields a contradiction in the Riemann-Hurwitz formula.

In view of (1ii), it remains to study the case of two-component $\widehat C=\widehat C_1\cup\widehat C_2$ such that
either $g(\widehat C_1)=2$, $g(\widehat C_2)=0$, or $g(\widehat C_1)=g(\widehat C_2)=1$. It follows that $\bw_\partial$ is constrained by an extra Menelaus condition, and this is the only condition imposed. Hence, all the points of $\widehat\bw_\partial$ differ from each other, and the genus of the image of any component must be equal to the genus of that component. The argument of the preceding paragraph leaves the only option of a multiple covering $\bn:\widehat C_2\twoheadrightarrow C_2\hookrightarrow\PP^2$ with rational $\widehat C_2$ and $C_2$. If $m$ is the degree of the covering, then each point $w_i^{(s)}\in C_2$ is a ramification point of the covering of index $m-1$, and hence we get a contradiction in the Riemann-Hurwitz formula:
$$2m-(m-1)\cdot\#(\widehat\bw_\partial\cap C_2)\ge2,\quad\#(\widehat\bw_\partial\cap C_2)\ge3.$$

\smallskip (2) We mainly follow the argument in parts (2ii)-(2v) of the proof of Lemma \ref{al5}. The items (2i), (2iv), (2v) of the present lemma can be treated in the same manner as the items (2i), (2iv), (2v) in Lemma \ref{al5}.
The case of (2v) appears to be even much simpler, since multiple coverings are not possible as we established in the preceding part of the proof. The cases (2ii) and (2iii) require additional arguments.

\smallskip (2ii)
Suppose that $\widehat\bw_\partial(t^*)$ consists of $n$ distinct points, $\widehat C$ is a smooth curve of genus $2$, $\bn:\widehat C\to C=\bn(\widehat C)$ is birational, but not an immersion. To verify the conditions of item (2ii) of the lemma, it is enough to consider the case of $C$ having a singular branch at some point $w_i^\sigma\in\widehat\bw_\partial(t^*)$. Fixing the position of the point $w_i^\sigma$, we obtain a subfamily $\pr(S^{fix}(w_i^\sigma))\subset\pr(S)$
of dimension at least $n-2\ge3$ and then apply \cite[Inequality (5) in Lemma 2.1]{IKS4} (cf. (\ref{ae9fX})):
$$-CK_{\Tor(P_\Delta)}\ge-2+\left(-CK_{\Tor(P_\Delta)}
-n+2\right)+(n-3)+(\mt Q-1)$$
$$=-CK_{\Tor(P_\Delta)}-3+(\mt Q-1)\quad\Longrightarrow\quad\mt Q\le4\ ,$$
where $Q$ is the local branch of $C$ at $w_i^\sigma$. It remains to show that $\mt Q=2$. Arguing by contradiction,
assume that $\mt Q=3$ or $4$.
Without loss of generality, we can assume that the germ of $\pr(S)\subset|{\mathcal L}_{P_\Delta}|$ at $C$ is an equisingular family. According to items (a) and (b) in Section \ref{sec-def}, we obtain the following upper bound completely analogous to (\ref{edop-alg2}):
\begin{equation}\dim\pr(S^{fix}(w_i^\sigma))\le h^0(C,{\mathcal J}_{Z/C}\otimes_{{\mathcal O}_C}\ina^*{\mathcal L}_{P_\Delta})=h^0\left(\widehat C,{\mathcal O}_{\widehat C}(-\bd)\otimes_{{\mathcal O}_{\widehat C}}\bn^*{\mathcal L}_{P_\Delta}\right),\label{ae9lX}\end{equation}
where
$$\deg\bd=2\sum_{z\in\Sing(C)}\delta(C,z)+(-CK_{\Tor(P_\Delta)}-n+2)+r,\quad 2\le r\le3,$$
$$\Longrightarrow\quad\deg({\mathcal O}_{\widehat C}(-\bd)\otimes_{{\mathcal O}_{\widehat C}}\bn^*{\mathcal L}_{P_\Delta})=C^2-\deg\bd=n-r\ge\begin{cases}3,\ &r=2,\\ 2,\ &r=3\end{cases}.$$

In case of $r=2$, we have $h^1({\mathcal O}_{\widehat C}(-\bd)\otimes_{{\mathcal O}_{\widehat C}}\bn^*{\mathcal L}_{P_\Delta})=0$ since $2g(\widehat  C)-2=2$. Thus,
by the Riemann-Roch formula, we have
$$\dim\pr(S^{fix}(w_i^\sigma))\le h^0({\mathcal O}_{\widehat C}(-\bd)\otimes_{{\mathcal O}_{\widehat C}}\bn^*{\mathcal L}_{P_\Delta})=n-3$$
contradicting the lower bound $\dim\pr(S^{fix}(w_i^\sigma))\ge n-2$ pointed above. In case of $r=3$, we have
$$h^0({\mathcal O}_{\widehat C}(-\bd)\otimes_{{\mathcal O}_{\widehat C}}\bn^*{\mathcal L}_{P_\Delta})\le h^8({\mathcal O}_{\widehat C}(-\bd')\otimes_{{\mathcal O}_{\widehat C}}\bn^*{\mathcal L}_{P_\Delta}),$$
where $\bd'=\bd-w$ with $w$ being one of the points of $\bd$. Thus,
$$\deg\bd'=\deg\bd-1\quad\Longrightarrow\quad\deg({\mathcal O}_{\widehat C}(-\bd')\otimes_{{\mathcal O}_{\widehat C}}\bn^*{\mathcal L}_{P_\Delta})=n-r+1\ge3,$$
and we derive the same contradiction.

\smallskip(2iii) Proceeding along the lines of part (2iii) of the proof of Lemma \ref{al5}, we need to
verify suitable analogs of the estimates (\ref{e2iii1}), (\ref{e2iii2}), and (\ref{e2iii3}) in the genus $2$ case.

An analog of (\ref{e2iii1}) reads
$$-CK_{\Tor(P_\Delta)}\ge-2+\left(-CK_{\Tor(P_\Delta)}
-n+4\right)+(n-4)+\sum_B(\mt Q-1)$$
$$=-CK_{\Tor(P_\Delta)}-2+\sum_B(\mt Q-1)\quad\Longrightarrow\quad\sum_B(\mt Q-1)\le2\ ,$$
where $B$ is the set of all singular local branches of $C$ in $\Tor(P_\Delta)^\times\cup\{w,w'\}$
(recall that $w'$ is any point in $\widehat\bw_\partial(t^*)\setminus\{w\}$).
We have analogs of (\ref{e2iii2}) and (\ref{e2iii3}):
$$\dim\pr(S^{fix}(w))\le h^0(\widehat C,{\mathcal O}_{\widehat C}(-\bd_1)\otimes_{{\mathcal O}_{\widehat C}}\bn^*{\mathcal L}_{P_\Delta}),$$
where
$$\deg\bd_1=2\sum_{z\in\Sing(C)}\delta(C,z)+(-CK_{\Tor(P_\Delta)}-n+3)+r,\quad r=1\ \text{or}\ 2,$$
and hence
$$\deg(\widehat C,{\mathcal O}_{\widehat C}(-\bd_1)\otimes_{{\mathcal O}_{\widehat C}}\bn^*{\mathcal L}_{P_\Delta})=C^2-\deg\bd_1=C^2-(C^2+CK_{\Tor(P_\Delta)}-2)$$ $$-(-CK_{\Tor(P_\Delta)}-n+3)-r=n-1-r\ge\begin{cases}3,\ & r=1,\\ 2,\ & r=2.\end{cases}$$
Then, we literally follow the last paragraph of part (2ii) of the proof of the lemma and derive that $\dim\pr(S^{fix}(w))\le n-3$ contrary to the initial lower bound $\dim\pr(S^{fix}(w))\ge n-2$. This contradiction completes the study of the case (2iii) of the lemma.
\proofend

\subsection{Wall-crossing}\label{crossing2}
In the analysis of wall-crossings, we
use arguments similar to the ones presented in Sections \ref{asec2} and \ref{asec5}.
The following lemma completes the proof of Theorem \ref{g2t1}.

\begin{lemma}\label{g2l2}
The numbers $W^\kappa_2(\Delta,\widehat\bw(t))$ and $\widetilde W^\kappa_2(\Delta,\widehat\bw(t))$ remain constant in each of the wall-crossing events described in Lemma \ref{g2l1}(2).
\end{lemma}

{\bf Proof.}
All local considerations in the proof of Theorems \ref{at1} and \ref{at2} can literally be transferred to our situation. Thus, we only confirm the global transversality relations which all take form of certain cohomology vanishing conditions.

\smallskip{\bf(1)} The wall-crossing described in Lemma \ref{g2l1}(2ii) coincides with that in Lemma \ref{al5}(2ii). The transversality condition in the treatment of this wall-crossing reads (cf. formula (\ref{ae10}))
\begin{equation}H^0(\widehat C,{\mathcal O}_{\widehat C}(\bd-p_1-p_2))=0,\quad\bd\in\Div(\widehat C),\label{g2e4}\end{equation}
where in the notation of the proof of Lemma \ref{al9}, part (2),
$$\deg\bd=C^2-\sum_{z\in\Sing(C)}2\delta(C,z)-\sum_{w\in\widehat\bw_\partial(t^*)}\dim{\mathcal O}_{C,w}/I_w$$
$$=4d^2-(4d^2-6d-2)-(6d-1)=3>(2g-2)\big|_{g=2}=2.$$
Thus, the Riemann-Roch formula yields $h^0(\widehat C,{\mathcal O}_{\widehat C}(\bd))=2$, and hence (\ref{g2e4}) holds due to the general position of $w_1,w_2$ on $C=\bn(\widehat C)$.

\smallskip{\bf(2)} The wall-crossing described in Lemma \ref{g2l1}(2iii) coincides with that in the proof of Lemma \ref{al5}(2iii). 
The required transversality condition coincides with relation (\ref{g2e4}) confirmed in the preceding paragraph.

\smallskip{\bf(3)} The wall-crossing described in Lemma \ref{g2l1}(2iv) coincides with that in Lemma \ref{al5}(2iv), the subcase of $\bn:\widehat C_i\to\Tor(\Delta)$, $i=1,2$, birational immersions. The local study is a word-for-word copy of the proof of Lemma \ref{al8}, part (5a) (restricted to the case of birational immersions). Similarly to the proof of Lemma \ref{al8}, part (5b), and the proof of Lemma \ref{al9}, part (4), the global transversality condition amounts to
$$h^0(\widehat C_1,{\mathcal O}_{\widehat C_1}(\bd_1-\{p_1,p_2\}\cap\widehat C_1))=h^0(\widehat C_2,{\mathcal O}_{\widehat C_2}(\bd_1-\{p_1,p_2\}\cap\widehat C_2))=0,$$
where
$$\deg\bd_1=2g(\widehat C_1)-1,\quad\deg\bd_2=2g(\widehat C_2)-1.$$ Thus, the required $h^1$-vanishing holds.

\smallskip{\bf(4)} The last wall-crossing to consider is associated with the wall described in Lemma \ref{g2l1}(2v). That wall is very similar to the wall appearing in Lemma \ref{al5}(2v). The local study of the wall-crossing is a word-for-word copy of the study performed in the proof of Lemma \ref{al9}, part (5).
Thus, we only address the global transversality condition analogous to that in formula (\ref{tce1}):
$$h^0(\widehat C,{\mathcal N}_{\bn'}(\bd))=0,$$
$$\bd=-\sum_{p_i^{(s)}\ne p}2k_i^{(s)} p_i^{(s)}+p_{i'}^{(s')}+p_{i''}^{(s'')}-(2l_1+1)p'-(2l_2+1)p''-p_1-p_2,$$
where $\bn':\widehat C\to\widehat C\overset{\bn}{\to}\PP^2$ is a composition of $\bn$ with the normalization, ${\mathcal N}_{\bn'}$ is the normal bundle, $p_{i'}^{(s')},p_{i''}^{(s'')}$ are mobile points in $\widehat\bw_\partial$ (while the other points are fixed), and $p',p''$ are the preimages of the node $p$ in $\widehat C$. It is easy to see that
$$\deg{\mathcal N}_{\bn'}(\bd)=0.$$
However, the points $w_1,w_2$ are in general position on $C=\bn(\widehat C)$, and hence
$p_1,p_2$ are in general position on $\widehat C$, which by the Riemann-Roch formula
implies the desired $h^0$-vanishing.
\proofend

{\ncsc Sorbonne Universit\'e and Universit\'e Paris Cit\'e, CNRS, IMJ-PRG \\[-21pt]

F-75005 Paris, France} \\[-21pt]

{\it E-mail address}: {\ntt     ilia.itenberg@imj-prg.fr}

\vskip10pt

{\ncsc School of Mathematical Sciences \\[-21pt]

Raymond and Beverly Sackler Faculty of Exact Sciences\\[-21pt]

Tel Aviv University \\[-21pt]

Ramat Aviv, 6997801 Tel Aviv, Israel} \\[-21pt]

{\it E-mail address}: {\ntt shustin@tauex.tau.ac.il}

\end{document}